\documentstyle[12pt,leqno]{article}
  \font\EUFMv	eufm5 scaled \magstephalf
  \font\EUFMvii	eufm7 scaled \magstephalf
  \font\EUFMx	eufm10 scaled \magstephalf
  \font\MSBMv	msbm5 scaled \magstep1
  \font\MSBMvii	msbm7 scaled \magstep1
  \font\MSBMx	msbm10 scaled \magstep1
  \font\EUFMxi   eufm10 scaled \magstep2

\newfam\German
\textfont\German=\EUFMx \scriptfont\German=\EUFMvii
\scriptscriptfont\German=\EUFMv
\def\frak{\fam=\German}

\newfam\Bgerman
\textfont\Bgerman=\EUFMxi
\def\bgerm{\fam=\Bgerman}

\newfam\Blackboad
\textfont\Blackboad=\MSBMx \scriptfont\Blackboad=\MSBMvii
\scriptscriptfont\Blackboad=\MSBMv
\def\Bbb{\fam=\Blackboad}

\def\newsection{\section}
\makeatletter
 \@addtoreset{equation}{subsection}
\makeatother
\renewcommand{\theequation}{\thesubsection.\arabic{equation}}

\def\GGg{{\bgerm{g}}}
\def\hb{\hfill\break}
\def\K{{\cal K}}
\newcommand{\CT}{{\cal T}}
\def\KR{{{\cal K}_{\rm reg}}}
\def\KRP{{{\cal K}^+_{\rm reg}}}

\def\Gsl{{\frak{sl}}}
\def\Ga{{\frak{a}}}
\def\Gb{{\frak{b}}}
\def\Gg{{\frak{g}}}
\def\Gh{{\frak{h}}}
\def\Gk{{\frak{k}}}

\def\Gn{{\frak{n}}}

\def\BA{{\Bbb A}}
\def\BC{{\Bbb C}}
\def\C{{\Bbb C}}
\def\BD{{\Bbb D}}

\def\BH{{\Bbb H}}
\def\BL{{\Bbb L}}
\def\BM{{\Bbb M}}
\def\BMA{{{\Bbb M}_{adm}(\Gg)}}
\def\MA{{{\Bbb M}_{adm}}}
\def\BMD{\HD}
\def\HD{{\Bbb H}}
\def\BN{{\Bbb N}}
\def\BO{{\Bbb O}}
\def\BP{{\Bbb P}}
\def\BQ{{\Bbb Q}}
\def\BR{{\Bbb R}}
\def\BZ{{\Bbb Z}}

\def\QC{{\overline{\Bbb Q}}}
\def\MH{{\rm MH}}
\def\DM{{\Bbb M}}
\def\Aut{{\rm Aut}}
\def\CA{{\cal A}}

\def\CB{{\cal B}}

\def\CD{{\cal D}}
\def\CE{{\cal E}}
\def\CH{{\cal H}}

\def\CL{{\cal L}}
\def\CM{{\cal M}}
\def\CN{{\cal N}}
\def\CO{{\cal O}}

\def\lam{\lambda}

\def\lan{\langle}
\def\ran{\rangle}

\newcommand{\ol}{\overline}
\def\on{\operatorname}
\def\operatorname#1{{\rm #1}}
\def\gge{>\kern-3pt>}

\def\mod{{\rm{mod}}}
\def\pt{{\rm{pt}}}
\def\im{{\rm{im}}}
\def\Supp{{\rm{Supp}}}
\def\ad{\mathop{\rm ad}\nolimits}
\def\pr{{\rm pr}}
\def\BUL{{\bullet}}
\def\ch{\mathop{\rm ch}\nolimits}
\def\Ch{{\rm{Ch}}}
\def\Coker{\mathop{\rm Coker}\nolimits}

\def\Cok{\mathop{\rm Coker}\nolimits}
\def\dual{{\bf D}}

\def\Ext{{\mathop{\rm Ext}\nolimits}}
\def\gr{{\mathop{\rm gr}\nolimits}}
\def\HT{{\rm ht}}
\def\ot{\otimes}
\def\htensor{{\hat{\otimes}}}
\def\Hom{\mathop{\rm Hom}\nolimits}
\def\id{\mathop{\rm id}\nolimits}

\def\Image{\mathop{\rm Im}\nolimits}
\def\Ker{\mathop{\rm Ker\hskip.5pt}\nolimits}
\def\limi{\mathop{\mathop{\rm lim}_{\longrightarrow}}}
\def\limp{\mathop{\displaystyle{\displaystyle\lim_{\longleftarrow}}}\limits}

\def\Ob{\mathop{\rm Ob}\nolimits}
\def\re{{\rm re}}
\def\sHom{\mathop{{\CH}om}\nolimits}
\def\sEnd{\mathop{{\CE}nd}\nolimits}

\def\Spec{\mathop{\rm Spec}\nolimits}

\newcommand{\var}{\on{var}}
\def\longhookrightarrow{%
\lhook\joinrel\relbar\joinrel\relbar\joinrel\rightarrow}

\def\TO{\mathop{\relbar\joinrel\longrightarrow}}

\def\tp{{\tilde p}}

\def\TX{{\tilde X}}

\def\TZ{\tilde Z}

\def\Ti{{\hbox{$\tilde \imath$}}}

\def\tH{\tilde{H}}
\def\tG{\tilde{\Gamma}}

\def\beq{\begin{eqnarray}}
\def\beqn{\begin{eqnarray*}}
\def\endeq{\end{eqnarray}}
\def\endeqn{\end{eqnarray*}}
\def\eq{\begin{eqnarray}}
\def\eqn{\begin{eqnarray*}}
\def\nn{\nonumber}
\def\proof{\noindent{\it Proof.}\quad}
\def\qed{\hspace*{\fill}\hbox{Q.E.D.}\hss\par\medskip}
\newtheorem{lemma}{Lemma}[subsection]
\newtheorem{definition}[lemma]{Definition}
\newtheorem{corollary}[lemma]{Corollary}
\newtheorem{sublemma}[lemma]{Sublemma}
\newtheorem{proposition}[lemma]{Proposition}
\newtheorem{theorem}[lemma]{Theorem}
\newtheorem{conjecture}{Conjecture}

\def\Conjecture{\begin{conjecture}}
\def\enconjecture{\end{conjecture}}

\def\Lemma{\begin{lemma}}
\def\enlemma{\end{lemma}}
\def\Definition{\begin{definition}}
\def\Def{\begin{definition}}
\def\endefinition{\end{definition}}
\def\Cor{\begin{corollary}}
\def\encor{\end{corollary}}
\def\Sublemma{\begin{sublemma}}
\def\ensublemma{\end{sublemma}}
\def\Proposition{\begin{proposition}}
\def\Prop{\begin{proposition}}
\def\enproposition{\end{proposition}}
\def\enprop{\end{proposition}}
\def\Theorem{\begin{theorem}}
\def\entheorem{\end{theorem}}

\def\ritem#1{\item{${\rm{#1}}$}}
\def\Remark{\medskip\noindent{\it Remark.}\quad}
\newcommand{\U}{{U(\Gg)}}
\newcommand{\overset}[2]{{\buildrel#1\over#2}}
\newcommand{\iso}{\mathrel{\longrightarrow{\kern-18pt\raise3.5pt%
\hbox{$\sim$}}}\enspace}

\newcommand{\TBO}{{\widetilde\BO}}

\newenvironment{tenumerate}{
  \begin{enumerate}
  \renewcommand{\labelenumi}
  {\rm
  {(\roman{enumi})}
  }
  }{\end{enumerate}}

\newcommand{\maprightu}[1]{%
\smash{\mathop{%
\hbox to 1cm{\rightarrowfill}}\limits^{#1}}}
\newcommand{\maprightd}[1]{%
\smash{\mathop{%
\hbox to 1cm{\rightarrowfill}}\limits_{#1}}}
\newcommand{\mapleftu}[1]{%
\smash{\mathop{%
\hbox to 1cm{\leftarrowfill}}\limits^{#1}}}
\newcommand{\mapleftd}[1]{%
\smash{\mathop{%
\hbox to 1cm{\leftarrowfill}}\limits_{#1}}}
\newcommand{\mapdownl}[1]{\Big\downarrow
\llap{$\vcenter{\hbox{$\scriptstyle#1\,$}}$ }}
\newcommand{\mapdownr}[1]{\Big\downarrow
\rlap{$\vcenter{\hbox{$\scriptstyle#1$}}$ }}

\newcommand{\mapr}[1]{%
\hbox{
\hbox to 28pt{\rightarrowfill}\kern-18pt\raise5pt\hbox{$\scriptstyle{#1}$}}}
\newcommand{\mapl}[1]{%
\smash{%
\hbox{\hbox to 25pt{\leftarrowfill}}\kern-17pt\raise.5pt\hbox{${}^{#1}$}}}

\newcommand{\downeq}{\vline\kern2pt\hbox{$\wr$}}
\newcommand{\e}{{\rm e}}
\newcommand{\ba}{\begin{array}}
\newcommand{\ea}{\end{array}}

\title{Kazhdan-Lusztig conjecture
for symmetrizable Kac-Moody Lie algebras. III\\
{\Large {\rule[.5ex]{25pt}{.7pt}}Positive rational case%
{\rule[.5ex]{25pt}{.7pt}}}}
\author{\it 
Dedicated to Professor Mikio Sato\\
\it in celebration of his seventieth birthday\\[10pt]
\sc Masaki Kashiwara%
\thanks{Research Institute for Mathematical Sciences, Kyoto University,
Kyoto, 606--8502, Japan} and
Toshiyuki Tanisaki%
\thanks{Department of Mathematics, Hiroshima University,
Higashi-Hiroshima, 739--8526, Japan}}
\begin{document}
\maketitle
\pagestyle{plain}

\tableofcontents

%
%
%
\newsection{Introduction}
\renewcommand{\theequation}{\thesection.\arabic{equation}}
\renewcommand{\thelemma}{\thesection.\arabic{lemma}}
The aim of this paper is to prove 
the Kazhdan-Lusztig type character formula 
for irreducible highest weight modules 
with positive rational highest weights over 
symmetrizable Kac-Moody Lie algebras.

Let us formulate our results precisely.
Let $\Gg$ be a symmetrizable Kac-Moody Lie algebra 
over the complex number field $\BC$ with Cartan subalgebra $\Gh$.
We denote by $W$ the Weyl group 
and by $\{\alpha_i\}_{i\in I}$  the set of simple roots.
For a real root $\alpha$, 
we define the corresponding coroot by 
$\alpha^\vee=2\alpha/(\alpha,\alpha)$, 
where $(\;,\;)$ denotes a standard non-degenerate 
symmetric bilinear form on $\Gh^*$.
For $\lam\in\Gh^*$, let $\Delta^+(\lam)$
denote the set of positive real roots 
$\alpha$ satisfying $(\alpha^\vee,\lam)\in\BZ$, 
and let $\Pi(\lam)$ denote the  set of 
$\alpha\in\Delta^+(\lam)$ such that 
$s_\alpha\big(\Delta^+(\lam)\setminus\{\alpha\}\big)
=\Delta^+(\lam)\setminus\{\alpha\}$.
Here $s_\alpha\in W$ denotes the reflection with respect to $\alpha$.
Then the subgroup $W(\lam)$ of $W$ generated by 
$\{s_\alpha\;;\;\alpha\in\Delta^+(\lam)\}$ 
is a Coxeter group with the canonical generator system 
$\{s_\alpha\;;\;\alpha\in\Pi(\lam)\}$.
Fix $\rho\in\Gh^*$ satisfying $(\rho,\alpha_i^\vee)=1$ for any $i\in I$ 
and define a shifted action of $W$ on $\Gh^*$ by 
\[
w\circ\lam=w(\lambda+\rho)-\rho\quad
\mbox{
for $w\in W$ and $\lambda\in\Gh^*$.}
\]

For $\lam\in\Gh^*$ let $M(\lam)$ (resp.\ $M^*(\lam)$, $L(\lam)$) 
be the Verma module (dual Verma module, irreducible module) 
with highest weight $\lam$. 
We denote their characters by 
$\ch(M(\lam))$, $\ch(M^*(\lam))$, $\ch(L(\lam))$ respectively.
We have $\ch(M(\lam))=\ch(M^*(\lam))$, and $\ch(M(\lam))$ is easily described.

The main result of this paper is the following.
\begin{theorem}
\label{introduction:Theorem1}
Assume that $\lam\in\Gh^*$ satisfies the following conditions.
\eq
\label{introduction:eq1}
&&\mbox{
$2(\alpha,\lam+\rho)\ne(\alpha,\alpha)$ 
for any positive imaginary root $\alpha$.
}\\
\label{introduction:eq2}
&&\mbox{
$(\alpha^\vee,\lam+\rho)\notin\BZ_{\leq0}$ 
for any positive real root $\alpha$.
}\\
\label{introduction:eq3}
&&\mbox{
If $w\in W$ satisfies $w\circ\lam=\lam$, then $w=1$.
}\\
\label{introduction:eq4}
&&\mbox{
$(\alpha^\vee,\lam)\in\BQ$ for any real root $\alpha$.
}
\endeq
Then for any $w\in W(\lam)$ we have 
\eq
&&
\ch(M(w\circ\lam))
=\sum_{y\ge_\lam w}
P^\lam_{w,y}(1)
\ch(L(y\circ\lam)),
\\
&&
\ch(L(w\circ\lam))
=\sum_{y\ge_\lam w}
(-1)^{\ell_\lam(y)-\ell_\lam(w)}
Q^\lam_{w,y}(1)
\ch(M(y\circ\lam)).\label{eq:1.6}
\endeq
Here, $\ge_\lam$, $\ell_\lam$, $P^\lam_{w,y}$, $Q^\lam_{w,y}$  
denote the Bruhat ordering, the length function, 
the Kazhdan-Lusztig polynomial, and 
the inverse Kazhdan-Lusztig polynomial for the Coxeter group 
$W(\lam)$, respectively.
\end{theorem}

When $\Gg$ is a finite-dimensional semisimple Lie algebra, 
this result for integral weights was conjectured by 
Kazhdan-Lusztig~\cite{KL1}, and proved  by 
Beilinson-Bernstein~\cite{BB} and Brylinski-Kashiwara~\cite{BK} independently.
Later its generalization to rational weights was obtained 
by combining the results by Beilinson-Bernstein 
(unpublished) and Lusztig~\cite{L0}.

As for the symmetrizable Kac-Moody Lie algebra, 
Theorem~\ref{introduction:Theorem1} 
for integral weights was obtained by Kashiwara(-Tanisaki) in 
\cite{Kpos1} and \cite{KTpos2} (see also Casian~\cite{C1}).
We note that a generalization of 
the original Kazhdan-Lusztig conjecture to affine Lie algebras 
in the negative level case was obtained by 
Kashiwara-Tanisaki~\cite{KTneg1}, 
Casian~\cite{C2} (integral weights), 
Kashiwara-Tanisaki~\cite{KTneg2} (rational weights).
We finally point out that 
(\ref{eq:1.6}) for $w=1$ was 
proved by Kac-Wakimoto~\cite{WK}.

Let us give a sketch of the proof of our theorem.

Let $X=G/B^-$ be the flag manifold introduced in Kashiwara~\cite{Kflag}, 
which is an infinite-dimensional scheme.
We have a stratification $X=\bigsqcup_{w\in W}X_w$ 
by finite-codimensional Schubert cells $X_w=BwB^-/B^-$.
For $\lam\in\Gh^*$ let $D_\lam$ be the TDO-ring 
(ring of twisted differential operators) on $X$ 
corresponding to the parameter $\lam$.
For $w\in W$ define $D_\lam$-modules $\CB_w(\lam)$ 
(resp.\  $\CM_w(\lam)$, $\CL_w(\lam)$) as the meromorphic extension 
(resp.\ dual meromorphic extension, minimal extension) 
of the $D_{X_w}$-module $\CO_{X_w}$ to a $D_\lam$-module.
They are objects of the category $\BH(\lam)$ 
consisting of $N^+$-equivariant holonomic $D_\lam$-modules.

For $\lam\in\Gh^*$ satisfying the conditions 
(\ref{introduction:eq1}), (\ref{introduction:eq2}) and
(\ref{introduction:eq3}), 
we define a modified global section functor $\tG$ 
from $\BH(\lam)$ to the category $\BM(\Gg)$ of $\Gg$-modules.
Then Theorem~\ref{introduction:Theorem1} 
is a consequence of the following results.
\begin{theorem}
\label{introduction:Theorem2}
Assume that $\lam\in\Gh^*$ satisfies 
the conditions {\rm (\ref{introduction:eq1}),
(\ref{introduction:eq2}) and (\ref{introduction:eq3})}.
\begin{tenumerate}
\item
The functor $\tG:\BH(\lam)\to\BM(\Gg)$ is exact.
\item
$\tG(\CB_w(\lam))=M^*(w\circ\lam)$ for any $w\in W$.
\item
$\tG(\CM_w(\lam))=M(w\circ\lam)$ for any $w\in W$.
\item
$\tG(\CL_w(\lam))=L(w\circ\lam)$ for any $w\in W$.
\end{tenumerate}
\end{theorem}
\begin{theorem}
\label{introduction:Theorem3}
Assume that $\lam\in\Gh^*$ satisfies the condition
$(\ref{introduction:eq4})$.
Then for any $w\in W$ which is the smallest element of $wW(\lam)$
and any $x\in W(\lam)$, we have
\eq
[\CL_{wx}(\lam)]
=\sum_{y\ge_\lam x}
(-1)^{\ell_\lam(y)-\ell_\lam(x)}
Q^\lam_{x,y}(1)
[\CM_{wy}(\lam)]
\endeq
in the $($modified$)$ Grothendieck group $K(\BH(\lam))$ of $\BH(\lam)$.
\end{theorem}

The proof of Theorem~\ref{introduction:Theorem2} 
is similar to the one in \cite{Kpos1}.
In the course of the proof we also use the modified localization functor 
$D_\lam\htensor\,\bullet\,$ from a category of certain $\Gg$-modules 
to a category of certain $D_\lam$-modules as in \cite{Kpos1}, 
and we prove simultaneously that 
$D_\lam\htensor\tG(\CM)\simeq\CM$ for any $\CM\in\Ob(\BH(\lam))$.
Aside from the technical complexity 
in dealing with non-integral weights,
the main new ingredients compared with the integral case \cite{Kpos1} 
are the embeddings of Verma modules (Theorem~\ref{KM:emb:thm:embed})
and the proof of the injectivity of
the canonical morphism $M(w\circ\lam)\to\tG(\CM_w(\lam))$
(Proposition~\ref{pro:inj}).

The proof of Theorem~\ref{introduction:Theorem3} is based on
the theory of Hodge modules by M. Saito~\cite{Saito} as in \cite{KTpos2}.
As for the combinatorics concerning the Kazhdan-Lusztig polynomials 
we use the dual version of the result in \cite{L}.

\medskip
In the affine case, we can deduce the non-regular highest weight case
from the above result by using the translation functors.

\Theorem
Let $\Gg$ be an affine Lie algebra, and assume that
$\lam\in\Gh^*$ satisfies
\eq
&&\mbox{
$(\delta,\lam+\rho)\not=0$, where $\delta$ is the imaginary root.
}\\
&&\mbox{
$(\alpha^\vee,\lam+\rho)\in\BQ\setminus\BZ_{<0}$ 
for any positive real root $\alpha$.
}
\endeq
Then $W_0(\lam)=\{w\in W;w\circ\lam=\lam\}$ is a finite group.
Let $w$ be an element of $W(\lam)$ which is
the longest element of $wW_0(\lam)$.
Then we have
\eqn
\ch(L(w\circ\lam))
=\sum_{y\ge_\lam w}
(-1)^{\ell_\lam(y)-\ell_\lam(w)}
Q^\lam_{w,y}(1)
\ch(M(y\circ\lam)).
\endeqn
\end{theorem}

A motivation of our study comes from a recent work of 
W. Soergel~\cite{Soergel} concerning tilting modules over affine Lie algebras.
We would like to thank H. H. Andersen for leading our attention 
to this problem.


\renewcommand{\theequation}{\thesubsection.\arabic{equation}}
\renewcommand{\thelemma}{\thesubsection.\arabic{lemma}}
\newsection{Highest weight modules} 
\subsection{Kac-Moody Lie algebras}
In this section, we shall review the definition of Kac-Moody Lie algebras, 
and fix notations employed in this paper.

Let $\Gh$ be a finite-dimensional vector space over $\BC$, and let
$\Pi=\{\alpha_i\}_{i\in I}$ and $\Pi^\vee=\{h_i\}_{i\in I}$ be subsets of
$\Gh^*$ and $\Gh$ respectively indexed by the same finite set $I$ subject to
\eq
&&
\mbox{$\Pi$ and $\Pi^\vee$ are linearly independent subsets of $\Gh^*$ and
$\Gh$ respectively,}
\\
&&
\mbox{$(\langle h_i,\alpha_j\rangle)_{i,j\in I}$ is a symmetrizable
generalized Cartan matrix.}
\endeq
Here $\lan\ ,\ \ran:\Gh\times\Gh^*\to\BC$ denotes the natural paring.
The elements of $\Pi$ and $\Pi^\vee$ are called simple roots and simple
coroots respectively.
We fix a non-degenerate symmetric bilinear form $(\ ,\ )$ on $\Gh^*$ such that
\begin{eqnarray}
&&(\alpha_i,\alpha_i)\in\BQ_{>0} \quad\mbox{ for any } i\in I,\\
&&\langle h_i,\lambda \rangle={2(\lambda,\alpha_i)}/{(\alpha_i,\alpha_i)}
\quad\mbox{for any $\lambda\in\Gh^*$ and $i\in I$.}
\end{eqnarray}

We denote the corresponding Kac-Moody Lie algebra by $\Gg$.
Recall that  $\Gg$ is the Lie algebra over $\BC$ generated by elements 
$e_i$, $f_i$ $(i\in I)$
and the vector space $\Gh$ satisfying the following defining relations
(see Kac~\cite{Kac}):
\eq
&&
\begin{array}{l}
[h,h']=0\quad \mbox{ for $h, h'\in\Gh$},
\\
{}[h,e_i]=\langle h,\alpha_i\rangle e_i, \quad
[h,f_i]=-\langle h,\alpha_i\rangle f_i
\quad\mbox{for $h\in\Gh$ and $i\in I$,}
\\
{}[e_i,f_j]=\delta_{ij}h_i  \quad\mbox{for $i$, $j\in I$,} \\
\ad(e_i)^{1-\langle h_i,\alpha_j\rangle}(e_j)
=\ad(f_i)^{1-\langle h_i,\alpha_j\rangle}(f_j)=0
\quad\mbox{for $i$, $j\in I$ with $i\ne j$.}
\end{array}
\endeq

Define the subalgebras $\Gn^+, \Gn^-, \Gb, \Gb^-$ of $\Gg$ by
\begin{eqnarray}
&&
\begin{array}{ll}
\Gn^+=\langle e_i\,;\,i\in I\rangle,\qquad\qquad
&\Gn^-=\langle f_i\,;\,i\in I\rangle,\\
\Gb=\langle h, e_i\,;\,h\in\Gh, i\in I\rangle,\qquad\qquad
&\Gb^-=\langle h, f_i\,;\,h\in\Gh, i\in I\rangle.
\end{array}
\end{eqnarray}
The vector space $\Gh$ is naturally regarded as an abelian subalgebra of
$\Gg$, and we have the decompositions
\begin{equation}\label{km:not:tri}
\Gg=\Gn^-\oplus\Gh\oplus\Gn^+,\quad
\Gb=\Gh\oplus\Gn^+,\quad
\Gb^-=\Gh\oplus\Gn^-.
\end{equation}

For $\lambda\in \Gh^*$ set
$\Gg_\lambda=
\{x\in\Gg\,;\,[h,x]=\langle h,\lambda\rangle x \mbox{ for } h\in \Gh\}$, and
define the root system $\Delta$ of $\Gg$ by
\begin{equation}
\Delta=\{\lambda\in\Gh^*\,;\, \Gg_\lambda\ne0\}\setminus\{0\}.
\end{equation}
Set
\begin{eqnarray}
&&Q=\sum_{i\in I}\BZ\alpha_i,\quad
Q^\pm=\pm\sum_{i\in I}\BZ_{\geq0}\alpha_i,\\
&&\Delta^\pm=\Delta\cap Q^\pm.
\end{eqnarray}
We have $\Pi\subset\Delta^+$ and  $\Delta=\Delta^+\sqcup\Delta^-$.
The elements of $\Delta^+$ and $\Delta^-$ are called positive and negative
roots respectively.

For a subset $\Theta$ of $\Delta$ such that
$(\Theta+\Theta)\cap\Delta\subset\Theta$ we define the subalgebra
$\Gn(\Theta)$ of $\Gg$ by
\eq
\label{eq:nTheta}
\Gn(\Theta)=\sum_{\alpha\in\Theta}\Gg_\alpha.
\endeq

For $\alpha=\sum_{i\in I}m_i\alpha_i\in Q$, 
its height $\HT(\alpha)$ is defined by
\eq\label{not:ht}
{\HT}(\alpha)=\sum_{i\in I}m_i.
\endeq

For $i\in I$ define the simple reflection $s_i\in GL(\Gh^*)$ by
\begin{equation}
s_i(\lambda)=\lambda-\langle h_i,\lambda\rangle\alpha_i.
\end{equation}
The subgroup $W$ of $GL(\Gh^*)$ generated by $S=\{s_i \,;\,i\in I\}$ is
called the Weyl group.

It is a Coxeter group with the canonical generator system $S$.
The length function  $\ell:W\to\BZ_{\geq0}$ of the Coxeter group $W$
satisfies
\eq
\ell(w)=\sharp(\Delta^-\cap w\Delta^+)
\qquad
\mbox{ for any }
 w\in W.
\endeq
We denote the Bruhat ordering on $W$ by $\ge$.
Note that we have
\begin{equation}
(w\lambda,w\mu)=(\lambda,\mu)\qquad
\mbox{ for any } \lambda, \mu\in\Gh^* \mbox{ and } w\in W.
\end{equation}

Set
\eq
&&
\Delta_{\rm re}=W\,\Pi\,,\quad
\Delta_{\rm im}=\Delta\setminus\Delta_{\rm re}\,,\quad
\Delta_{\rm re}^\pm=\Delta_{\rm re}\cap\Delta^\pm\,,\quad
\Delta_{\rm im}^\pm=\Delta_{\rm im}\cap\Delta^\pm\,.
\endeq
The elements of $\Delta_{\rm re}$ and $\Delta_{\rm im}$ are called real and
imaginary roots, respectively.
For $\alpha\in\Delta_{\rm re}$ set
\eq
\alpha^\vee=2\alpha/(\alpha,\alpha),
\endeq
and define the reflection $s_\alpha\in GL(\Gh^*)$ by
\eq
s_\alpha(\lambda)=\lambda-(\lambda,\alpha^\vee)\alpha.
\endeq
Then we have $s_\alpha\in W$ for any $\alpha\in\Delta_{\rm re}$.

We fix a vector $\rho\in\Gh^*$ such that
$\lan h_i,\rho\ran=1$ for any $i\in I$.
Then the shifted action of $W$ on $\Gh^*$ is defined by
\eq\label{def;shift}
w\circ \lam=w(\lam+\rho)-\rho\,.
\endeq
Note that $\rho-w\rho\in Q^+$ for any $w\in W$ and it does
not depend on the choice of $\rho$.

\subsection{Integral Weyl groups}\label{sec:weyl}

In this section, we study the properties of integral Weyl groups.
We start the study in a more general setting
\footnote{
After writing up this paper, the authors were informed 
by S. Naito
the existence of two papers,
R. Moody--A. Pianzola, Lie Algebras with Triangular Decompositions,
Canadian Mathematical Society series of monographs and advanced texts,
A Wiley-Interscience Publication, John Wiley \& Sons, 1995,
and
Jong-Min Ku, On the uniqueness of embeddings of Verma modules defined by the
Shapovalov elements, J. Algebra, Vol. 118, (1988) 85--101.
They showed results similar to those in this subsection by a different 
formulation and method. In the last paper, Ku also obtained
a result weaker than Theorem \ref{KM:emb:thm:embed}.}.

Let $\Delta_1$ be a subset of $\Delta_\re$
satisfying the following condition:
\eq\label{cond:root}
\hbox{ $s_\alpha\beta\in\Delta_1$ for any $\alpha,\beta\in\Delta_1$.}&&
\endeq
In particular, we have $-\Delta_1=\Delta_1$.
We set
\eq
\Delta_1^\pm&=&\Delta_1\cap\Delta^\pm,\\
\Pi_1&=&\{\alpha\in\Delta_1^+\,;\, 
s_\alpha(\Delta_1^+\setminus\{\alpha\})\subset\Delta_1^+\},\\
W_1&=&\lan s_\alpha;\alpha\in\Pi_1\ran\subset W.
\endeq
We call the elements of $\Delta_1^+$ (resp.\ $\Delta_1^-$, $\Pi_1$) positive
roots (resp.\ negative roots, simple roots) for $\Delta_1$, and $W_1$ the
Weyl group for $\Delta_1$.

Note that if $\Delta_1$ satisfies the condition (\ref{cond:root}),
then $w\Delta_1$ also satisfies (\ref{cond:root}) for any $w\in W$.

\Lemma \label{km:int:lem:1}
If $\Delta_1$ contains a simple root $\alpha_i$, then
$\alpha_i$ is in $\Pi_1$.
\enlemma

\Lemma\label{km:int:lem:2}
Assume that $\alpha_i\not\in \Delta_1$.
Set $\Delta_1'=s_i\Delta_1$.
Then $\Delta_1'$ satisfies the condition {\rm (\ref{cond:root})}.
Moreover $\Delta_1'\bigcap\Delta^+=s_i\Delta_1^+$,
the set of simple roots for $\Delta_1'$ is $s_i\Pi_1$,
and the Weyl group for $\Delta_1'$ is $s_iW_1s_i$.
\enlemma
The above two lemmas 
immediately follow from
$s_i\Delta^+=(\Delta^+\setminus\{\alpha_i\})\cup\{-\alpha_i\}$.

\Lemma\label{lemma:simple}
 If $\alpha\in\Pi_1$ and $i\in I$ satisfy
$(\alpha_i,\alpha)>0$, then
either $\alpha=\alpha_i$ or $\alpha_i\not\in\Delta_1$.
\enlemma
\proof
Assume $\alpha\not=\alpha_i$ and $\alpha_i\in\Delta_1$. 
Then $\beta=s_{\alpha}\alpha_i=\alpha_i-(\alpha^\vee,\alpha_i)\alpha$
is a positive root.
Then $\alpha_i=\beta+(\alpha^\vee,\alpha_i)\alpha$
contradicts $(\alpha^\vee,\alpha_i)\in\BZ_{>0}$ .
\qed

\Lemma\label{km:int:lem:8}
For any $\alpha\in\Pi_1$ there exist $w\in W$ and $i\in I$ such that
$w\alpha=\alpha_i$ and 
$w\Delta_1^+=w\Delta_1\bigcap\Delta^+$.
\enlemma
\proof
We shall show this by induction on $\HT(\alpha)$.
If $\HT(\alpha)=1$, then there is nothing to prove.
Assume $\HT(\alpha)>1$.
Write $\alpha=\sum_{j\in I}m_j\alpha_j$ with $m_j\ge0$.
Then we have
\[
0<(\alpha,\alpha)=\sum_{j\in I}m_j(\alpha,\alpha_j),
\]
and hence there exists some $j\in I$ such that $(\alpha,\alpha_j)>0$.
Since $\HT(\alpha)>1$, we have $\alpha\ne\alpha_j$ and hence
$\alpha_j\notin \Delta_1$ by Lemma \ref{lemma:simple}.
Set $\Delta_1'=s_j\Delta_1, (\Delta_1')^+=s_j\Delta_1^+, \Pi_1'=s_j\Pi_1$.
Then $(\Delta_1')^+$ and $\Pi_1'$ are the set of positive and simple roots
for $\Delta_1'$ respectively by Lemma \ref{km:int:lem:2}.
Set $\alpha'=s_j\alpha\in\Pi_1'$.
Since $\HT(\alpha')<\HT(\alpha)$, there exist some $w'\in W$
and $i\in I$ such that $w'\alpha'=\alpha_i$ and
$w'(\Delta_1')^+\subset\Delta^+$ by the hypothesis of induction.
Then setting $w=w's_j$ we have $w\alpha=\alpha_i$ and
$w\Delta_1^+=w'(\Delta_1')^+\subset \Delta^+$.
\qed

\medskip
The following lemma follows from the above lemma
by reducing to the case $\alpha=\alpha_i$ for $i\in I$.

\Lemma\label{lemma:simp2}
For any positive integer $n$ and $\alpha\in\Pi_1$, we have
$n\alpha\not\in\sum
\limits_{\beta\in(\Delta_1^+\setminus\{\alpha\})\cup\Delta_\im^+}
\BZ_{\ge0}\beta$.
\enlemma

\Lemma\label{km:int:lem:4}
For any $\beta\in\Delta_1^+$, there exists $\alpha\in\Pi_1$
such that $(\alpha,\beta)>0$.
\enlemma
\proof
We shall prove this by the induction on $\HT(\beta)$.
If $\HT(\beta)=1$, then  Lemma \ref{km:int:lem:1} implies
$\beta\in\Pi_1$, and we can take $\beta$ as $\alpha$.
Assume that $\HT(\beta)>1$.
Take $i$ such that $(\alpha_i,\beta)>0$.
If $\alpha_i\in\Delta_1$, then
it is enough to take $\alpha_i$ as $\alpha$.
Now assume that $\alpha_i\not\in \Delta_1$.
Set $\Delta_1'=s_i\Delta_1, (\Delta_1')^+=s_i\Delta_1^+, \Pi_1'=s_i\Pi_1$.
Then $(\Delta_1')^+$ and $\Pi_1'$ are the set of positive and simple roots
for $\Delta_1'$ respectively by Lemma \ref{km:int:lem:2}.
Set $\beta'=s_i\beta\in(\Delta_1')^+$.
We have $\HT(\beta')<\HT(\beta)$ by $(\beta,\alpha_i)>0$.
Hence by the hypothesis of induction there exists some $\alpha'\in\Pi_1'$
such that $(\alpha',\beta')>0$.
Then $\alpha=s_i\alpha'\in\Pi_1$ satisfies
$(\alpha,\beta)=(\alpha',\beta')>0$.
\qed

\Lemma \label{km:int:lem:5}
\begin{tenumerate}
\item
$\Delta_1=W_1\Pi_1$.
\item
$\Delta_1^+\subset \sum_{\alpha\in\Pi_1}\BZ_{\ge0}\alpha$.
\item
$W_1$ contains $s_\alpha$ for any $\alpha\in\Delta_1$.
\end{tenumerate}
\enlemma
\proof
Since (iii) follows from (i), it is enough to show that any $\beta\in
\Delta_1^+$ is contained in
$W_1\Pi_1\bigcap \sum_{\alpha\in\Pi_1}\BZ_{\ge0}\alpha$.
We shall prove this by the induction on $\HT(\beta)$.
By Lemma~\ref{km:int:lem:4}  there exists $\alpha_0\in\Pi_1$ such that
$(\alpha_0,\beta)>0$.
If $\beta=\alpha_0$, then there is nothing to prove.
If $\beta\not=\alpha_0$, then
$\gamma=s_{\alpha_0}\beta\in\Delta_1^+$ by 
the definition of $\Pi_1$
and
$\HT(\gamma)<\HT(\beta)$.
Now we can apply the hypothesis of induction to conclude
$\gamma\in W_1\Pi_1\bigcap \sum_{\alpha\in\Pi_1}\BZ_{\ge0}\alpha$,
which implies the desired result.
\qed





\Lemma\label{lemma:equiv}
For $\alpha\in\Delta_1^+$, the following conditions are equivalent.
\begin{tenumerate}
\item
$\alpha\in\Pi_1$.
\item
$s_\beta\alpha\in\Delta_1^-$ 
for any $\beta\in\Delta_1^+$ such that $(\alpha,\beta)>0$.
\item
$\alpha$ cannot be written as $\alpha=m_1\beta_1+m_2\beta_2$
for $\beta_\nu\in\Delta_1^+$ and $m_\nu\in\BZ_{>0}$
$($$\nu=1$, $2$$)$.
\item
$\alpha$ cannot be written as $\alpha=\sum_{\nu=1}^k\beta_\nu$
for $k>1$, $\beta_\nu\in\Delta_1^+$ $($$1\le\nu\le k$$)$.
\end{tenumerate}
\enlemma
\proof
(i)$\Rightarrow$(iv) follows from Lemma \ref{lemma:simp2}.
(iv)$\Rightarrow$(iii) is trivial.
(iii)$\Rightarrow$(ii) is also immediate.
Let us prove (ii)$\Rightarrow$(i).
By Lemma \ref{km:int:lem:4}, there exists $\beta\in\Pi_1$ such that
$(\beta,\alpha)>0$. 
Hence $\gamma=-s_\beta\alpha\in\Delta_1^+$.
Rewriting this, we have $(\beta^\vee,\alpha)\beta=\alpha+\gamma$.
Then Lemma \ref{lemma:simp2} implies $\alpha=\beta$
or $\gamma=\beta$. It is now enough to remark that
$\gamma=\beta$ implies $\alpha=\beta$.
\qed

The following proposition is proved
by a standard argument
(see e.g. \cite[\S 3.2]{KTneg2}).

\Prop \label{km:int:lCoxeter}
\begin{tenumerate}
\item
$W_1$ is a Coxeter group with a generator system
$S_1=\{s_\alpha\,;\,\alpha\in\Pi_1\}$.
\item
Its length function $\ell_1:W_1\to\BZ_{\ge0}$ is given by
$\ell_1(w)=\sharp(\Delta_1^-\cap w\Delta_1^+)$.
\item
For $x,y\in W$, $x\ge_1y$ 
with respect to the Bruhat order
$\ge_1$ for $(W_1,S_1)$ if and only if there exist
$\beta_1,\ldots,\beta_r\in\Delta_1^+$ $($$r\ge0$$)$
such that $x=ys_{\beta_1}\cdots s_{\beta_r}$ and
$ys_{\beta_1}\cdots s_{\beta_{j-1}}\beta_j\in\Delta_1^+$ for
$j=1,\ldots,r$.
\end{tenumerate}
\enprop

\Lemma\label{km:int:lem:6}
For $\alpha,\beta\in\Pi_1$ such that $\alpha\not=\beta$
we have $(\alpha,\beta)\le0$.
\enlemma
\proof
We have $s_\alpha\beta\in\Delta_1^+$ by the definition of $\Pi_1$.
Since $\beta\in\Pi_1$, Lemma \ref{lemma:equiv} implies the desired result.
\qed

\medskip
By this lemma,
$\Big((\beta,\alpha^\vee)\Big)_{\alpha,\beta\in\Pi_1}$ is a
symmetrizable generalized Cartan matrix.
Hence $W_1$ is isomorphic to the Weyl group for
the Kac-Moody Lie algebra with
$\Big((\beta,\alpha^\vee)\Big)_{\alpha,\beta\in\Pi_1}$ as
a generalized Cartan matrix.

\Prop\label{prop:cong}
For $w\in W$ the following conditions are equivalent.
\begin{tenumerate}
\item
$l(x)\ge l(w)$ for any $x\in wW_1$.
\item
$wx\ge wy$ for any $x,y\in W_1$ such that $x\ge_1y$.
\item
$w\Delta_1^+\subset \Delta^+$.
\end{tenumerate}
\enprop

\proof
Let us first prove (iii)$\Rightarrow$(ii).
We may assume without loss of generality that
$x=ys_\beta$ for some $\beta\in\Delta_1^+$.
Then $y\beta\in\Delta_1^+$ and hence
$wy\beta\in\Delta^+$.
This implies $wx=wys_{\beta}\geq wy$.

(ii) implies (i) by taking $y=1$ in (ii).
(i) implies (iii) because, for any $\alpha\in\Delta^+_1$,
$l(ws_\alpha)\ge l(w)$ implies $w\alpha\in\Delta^+$.
\qed

\bigskip
For $\lam\in \Gh^*$ set
\eq
&&
\ba{rcl}
\Delta(\lam)&=&\{\beta\in\Delta_\re;\,(\beta^\vee,\lam+\rho)\in\BZ\}\\
&=&\{\beta\in\Delta_\re;\,(\beta^\vee,\lam)\in\BZ\}.
\ea
\endeq
This satisfies the condition (\ref{cond:root}).
We set $\Delta^{\pm}(\lambda)=\Delta(\lambda)\cap\Delta^\pm$.
Let $\Pi(\lam)$ and $W(\lambda)$ be the set of simple roots
and the Weyl group for $\Delta(\lam)$, respectively.
We call $W(\lam)$ the {\em integral Weyl group} for $\lam$.
We denote by $\ell_\lambda:W(\lambda)\to\BZ_{\geq0}$ and $\geq_\lambda$ the
length function and the Bruhat order of the Coxeter group
$W(\lambda)$, respectively.

\Remark%
\begin{tenumerate}
\item
In \cite{KTneg2}, we introduced
$W(\lam)$ and $W'(\lam)$.
The integral Weyl group introduced here
is equal to $W'(\lam)$ loc.cit.
As a matter of fact, $W(\lam)$ and $W'(\lam)$ loc.cit. coincide.
The opposite statement in \cite[Remark 3.3.2]{KTneg2} should be corrected.

\item
The set $\Pi(\lam)$ is linearly independent 
when $\Gg$ is finite-dimensional. 
But it is not necessarily linearly independent in the affine case,
although we have assumed the linear independence of $\{\alpha_i\}_{i\in I}$.
For example, for $\Gg=A^{(1)}_3$ and $\lam=(\Lambda_1+\Lambda_3)/2$,
we have 
$\Pi(\lam)=\{\alpha_0,\,\alpha_2,\,-\alpha_0+\delta,\,-\alpha_2+\delta\}$.
\item
For $x,y\in W_1$, $x\ge_1y$ implies $x\ge y$ (Lemma \ref{prop:cong}).
However the converse is false in general.
For example for $\Gg=A_3$ and $\lam=(\Lambda_1+\Lambda_3)/2$,
we have $\Pi(\lam)=\{\alpha_2,\,\alpha_1+\alpha_2+\alpha_3\}$
and $s_{\alpha_1+\alpha_2+\alpha_3}\ge s_{\alpha_2}$.
\end{tenumerate}

\subsection{Category of highest weight modules}
In this subsection we shall recall some properties of the category $\BO$
of highest weight $\Gg$-modules.

In general, for a Lie algebra $\Ga$ we denote its enveloping algebra by
$U(\Ga)$ and the category of (left) $U(\Ga)$-modules by $\BM(\Ga)$.

For $k\in\BZ_{\ge0}$ set
\eq
\label{eq:nkpm}
\Gn_k^\pm=\Gn(\pm\Delta^+_k)
\;\mbox{ with }
\Delta^+_k=\{\alpha\in\Delta^+\,;\,\HT(\alpha)\ge k\}
\endeq
(see (\ref{eq:nTheta}) and (\ref{not:ht}) for the notation).
A $\U$-module $M$ is called {\em admissible} if, for any $m\in M$, there exists
some $k$ such that $\Gn^+_km=0$.
We denote by $\BMA$ 
the full subcategory of $\BM(\Gg)$ consisting of admissible $\U$-modules.
It is obviously an abelian category.

For $M\in\BM(\Gh)$ and $\xi\in\Gh^*$ we set
\[
M_\xi=\{u\in M\,;\,
\mbox{ $(h-\lan h,\xi\ran)^nu=0$ for any $h\in\Gh$ and $n\gge0$}
\}.
\]
It is called the generalized weight space of $M$ with weight $\xi$.
We denote by $\BO$ the full subcategory of $\BM(\Gg)$ consisting of
$U(\Gg)$-modules $M$ satisfying
\eq
&&M=\bigoplus_{\xi\in\Gh^*}M_\xi,\label{KM:cat:eq1}\\
&&\dim M_\xi<\infty\quad \mbox{ for any $\xi\in\Gh^*$},\label{KM:cat:eq2}\\
&&\mbox{for any $\xi\in\Gh^*$ there exist only finitely many
$\mu\in\xi+Q^+$ such that $M_\mu\ne0$.}\label{KM:cat:eq3}
\endeq
It is an abelian subcategory of $\BMA$.

For $M\in\Ob(\BO)$, or more generally for an $\Gh$-module $M$ satisfying
(\ref{KM:cat:eq1}) and (\ref{KM:cat:eq2}),
we define its character as the formal infinite sum
\[
\ch(M)=\sum_{\xi\in\Gh^*}(\dim M_\xi)\, \e^\xi.
\]
For a $U(\Gg)$-module $M$, the dual space $\Hom_\BC(M,\BC)$ is endowed with a
$U(\Gg)$-module structure by
\[
\langle xm^*, m\rangle
= \langle m^*, a(x)m\rangle\quad
\mbox{for $m^*\in\Hom_\BC(M,\BC)$, $m\in M$, $x\in\Gg$,}
\]
where $a:\Gg\to\Gg$ is the anti-automorphism of the Lie algebra $\Gg$ given by
\[
a(h)=h\quad\mbox{ for } h\in\Gh, \qquad
a(e_i)=f_i,\quad a(f_i)=e_i\quad\mbox{ for } i\in I.
\]
If $M\in\Ob(\BO)$, then
\[
M^*:=\bigoplus_{\xi\in\Gh^*}(M_\xi)^*\subset\Hom_\BC(M,\BC)
\]
is a $\U$-submodule of $\Hom_\BC(M,\BC)$ belonging to $\Ob(\BO)$.
Indeed we have
\eqn
&&(M^*)_\xi=(M_\xi)^*\,.
\endeqn
Moreover, it defines a contravariant exact functor $(\bullet)^*:\BO\to\BO$
such that $(\bullet)^{**}$ is naturally isomorphic to the identity functor
on $\BO$.
In particular, we have
\eq
\Hom_\Gg(M,N)\simeq\Hom_\Gg(N^*,M^*)
\quad\mbox{ for $M, N\in\Ob(\BO)$}.
\endeq
We also note
\eq
\mbox{$\ch(M^*)=\ch(M)$ for any $M\in\Ob(\BO)$.}
\endeq

An element $m$ of a $\U$-module $M$ is called a highest weight vector with
weight $\lambda$ if $m\in M_\lambda$ and $e_im=0$ for any
$i\in I$.
A $\U$-module $M$ is called a highest weight module with highest weight
$\lambda$ if it is generated by a highest weight vector with weight $\lam$.
Highest weight modules belong to the category $\BO$.

For $\lambda\in\Gh^*$ define  a highest weight module $M(\lambda)$ with
highest weight $\lambda$, called a Verma module, by
\[
M(\lambda)=\U/(\sum_{h\in\Gh}\U(h-\lambda(h))+
\sum_{i\in I}\U e_i).
\]
The element of $M(\lambda)$ corresponding to $1\in U(\Gg)$ will be denoted
by $u_\lam$.
Set $M^*(\lam)=(M(\lam))^*$.
There exists a unique (up to a constant multiple) non-zero homomorphism
$M(\lam)\to M^*(\lam)$.
Its image $L(\lam)$ is a unique irreducible quotient of $M(\lam)$ 
and a unique irreducible submodule of $M^*(\lam)$.
In particular, we have $(L(\lam))^*\simeq L(\lam)$.

We have the following lemma (see Lemma 9.6 of Kac~\cite{Kac}).
\begin{lemma}
\label{km:not:lemma:mult}
For any $M\in\Ob(\BO)$ and $\mu\in\Gh^*$, there exists a finite filtration
\[
0=M_0\subset M_1\subset\cdots\subset M_r=M
\]
of $M$ by $\U$-modules $M_k$ $($$k=0,\ldots, r$$)$
such that for any $k$ we have either
$(M_k/M_{k-1})_\mu=0$ or
$M_k/M_{k-1}\simeq L(\xi)$ for some $\xi\in\Gh^*$.

\end{lemma}

For $M\in\Ob(\BO)$ and $\mu\in\Gh^*$ we set
\[
[M:L(\mu)]=\sharp\{k\,;\,M_k/M_{k-1}\simeq L(\mu)\},
\]
for a filtration of $M$ as in Lemma~\ref{km:not:lemma:mult}.
%
It does not depend on the choice of a filtration.
Then we have the equalities
\eqn
&&[M:L(\mu)]=[M^*:L(\mu)]\,,\\[3pt]
&&\ch(M)=\sum_{\mu\in\Gh^*}[M:L(\mu)]\ch(L(\mu)).
\endeqn
We frequently use the following lemma later.

\begin{lemma}
\label{KM:cat:dimHom}
Let $M\in\Ob(\BO)$ and $\mu\in\Gh^*$.
\begin{tenumerate}
\item
$\dim\Hom_\Gg(M(\mu),M)$ and $\dim\Hom_\Gg(M,M^*(\mu))$
are less than or equal to $[M:L(\mu)]$.
\item
Assume that if $\xi\in\Gh^*$ satisfies $[M:L(\xi)]\ne0$ 
and $[M(\xi):L(\mu)]\ne0$, then $\xi=\mu$.
Assume further $[M:L(\mu)]\ne0$.
Then neither $\Hom_\Gg(M(\mu),M)$ nor $\Hom_\Gg(M,M^*(\mu))$
vanishes.
\end{tenumerate}
\end{lemma}
\proof
(i) is obvious.
Let us prove (ii).
Consider the set $\CA$ of submodules $R$ of $M$ satisfying $[R:L(\mu)]=0$.
There exists the largest element $K$ of $\CA$ with
respect to the inclusion relation.
Set $N=M/K$.\hb
We shall prove $N_{\mu+\gamma}=0$ for any $\gamma\in Q^+\setminus\{0\}$.
Assume that there exist some $\gamma\in Q^+\setminus\{0\}$ such that
$N_{\mu+\gamma}\ne0$.
Since $N$ is an object of $\BO$, there exists finitely many such $\gamma$.
Take $\gamma\in Q^+\setminus\{0\}$ such that
$N_{\mu+\gamma}\ne0$ and $N_{\mu+\gamma+\delta}=0$ for any $\delta\in
Q^+\setminus\{0\}$.
Then we have $[N:L(\mu+\gamma)]>0$.
Let $N'$ be the $\Gg$-submodule of $N$ generated by $N_{\mu+\gamma}$.
By the maximality of $K$ we have $[N':L(\mu)]\ne0$.
Hence, $[M(\mu+\gamma):L(\mu)]\ne0$ by the construction of $N'$.
This contradicts
\[
[M:L(\mu+\gamma)]\ge[N:L(\mu+\gamma)]>0.
\]
Hence $N_{\mu+\gamma}=0$ for any $\gamma\in Q^+\setminus\{0\}$,
which implies
\eqn
&&\Hom(M(\mu),N^*)=
\{u\in N^*_\mu;hu=\lam(h)u\quad\mbox{for any $h\in\Gh$}\}.
\endeqn
Since $\dim N^*_\mu\ge[N:L(\mu)]=[M:L(\mu)]>0$, 
$\Hom(M(\mu),N^*)$ does not vanish.
Hence $\Hom(M(\mu),M^*)$ which contains $\Hom(M(\mu),N^*)$,
does not vanish either.
By applying the same argument to $M^*$ we have
$\Hom(M(\mu),M)\ne0$.
\qed

\subsection{Enright functor for non-integral weights}
In order to obtain some results on Verma modules
(Proposition \ref{KM:bim:prop:HomMM}), 
we construct a version of
Enright functor with non-integral weights (see Enright~\cite{Enright},
Deodhar~\cite{Deodhar}).

Since the action of $\ad(f_i)$ on $\U$ is locally nilpotent,
the ring $\U[f_i^{-1}]$ ,
a localization of $\U$ by $f_i$,
is well-defined.
It contains $\U$ as a subring.
Similarly we can consider a $\U$-bimodule
$\U f_i^{a+\BZ}$ for any scalar $a\in\BC$.
As a left $\U$-module it is given by
\eq
\U f_i^{a+\BZ}=
\limi_n
\U f_i^{a-n},
\endeq
where $\U f_i^{a-n}$ is a rank one free $U(\Gg)$-module generated by the
symbol $f_i^{a-n}$ and the homomorphism $\U f_i^{a-n}\to\U f_i^{a-n-1}$ is
given by $f_i^{a-n}\mapsto f_if_i^{a-n-1}$.
The left module $\U f_i^{a-n}$ is naturally identified 
with a submodule of $\U f_i^{a+\BZ}$ 
and we have $\U f_i^{a+\BZ}=\bigcup_{n\in\BZ}\U f_i^{a-n}$.
Its right module structure is given by
\eq
\label{eq:Enright1}
\qquad
f_i^{a+m}P=\sum_{k=0}^{\infty} {a+m\choose k}(\ad(f_i)^kP)f_i^{a+m-k}
\mbox{ for any $m\in\BZ$ and any $P\in\U$}.
\endeq
As a right $\U$-module, we also have
\[
\U f_i^{a+\BZ}=
\limi_n
f_i^{a-n}\U
=\bigcup_n f_i^{a-n}\U.
\]
By (\ref{eq:Enright1}) we have
\eq
Pf_i^{a+m}=\sum_{k=0}^\infty (-1)^k{a+m\choose k}f_i^{a+m-k}(\ad(f_i)^kP).
\endeq
In particular, we have
\eq
e_i^n f_i^a=\sum_{k=0}^nf_i^{a-k}e_i^{n-k}(k!)^2{n\choose k}{a\choose k}
{h_i+n-a\choose k}.
\endeq
The $\U$-bimodule $\U f_i^{a+\BZ}$ depends only on $a$ modulo $\BZ$.

\Lemma\label{lem:comp}
For $a$, $b\in\BC$, the map $f_i^{a+b+n}\mapsto f_i^{a+n}\otimes f_i^b
=f_i^{a}\otimes f_i^{b+n}$
$($$n\in \BZ$$)$
defines an isomorphism of $\U$-bimodules
\[\U f_i^{a+b+\BZ}\to \U f_i^{a+\BZ}\otimes_\U\U f_i^{b+\BZ}\,.\]
\enlemma

Since the proof is straightforward, we omit it.
Hence $\bigoplus\limits_{a\in\BC/\BZ}\U f_i^{a+\BZ}$ has a structure of a ring
containing $\U$.

\medskip
For any $\Gg$-module $M$,
$\U f_i^{a+\BZ}\otimes_{\U} M$ is isomorphic
to the inductive limit
$$M\overset{f_i}{\longrightarrow}M\overset{f_i}{\longrightarrow}M\overset{f_
i}{\longrightarrow}\cdots$$
as a vector space.
%
Hence we obtain the following result.
\Prop
The functor $M\to \U f_i^{a+\BZ}\otimes_{\U} M$
is an exact functor from $\BM(\Gg)$ into itself.
\enprop

Let $\Ga$ be a Lie algebra, and let $\Gk$ be its subalgebra
such that $\Ga$ is locally $\Gk$-finite with respect to the adjoint action.
Then for any $\Ga$-module $M$, the subspace
$\{m\in M\,;\,\dim U(\Gk)m<\infty\}$
is an $\Ga$-submodule of $M$.
\hb
In particular, for a $\Gg$-module $M$ its subspace
$
\{m\in M\,;\,\dim \BC[e_i]m<\infty\}
$
is a $\Gg$-submodule of $M$.
For $a\in\BC$, we define a functor
\eq
T_i(a):\BM(\Gg)\to\BM(\Gg)
\endeq
by
\[
T_i(a)(M)=\{u\in \U f_i^{a+\BZ}\otimes_{\U} M\,;\, \dim
\BC[e_i]u<\infty\}
\]
for $M\in\Ob(\BM(\Gg))$.
It is obviously a left exact functor.

For $a\in\BC$, let
$\BM_i^a(\Gg)$ be the category of locally $\BC[e_i]$-finite $\U$-modules
$M$ such that $M$ has a weight decomposition
\[M=\bigoplus_{\lam\in\Gh^*}M_\lam,\]
the action of $\Gh$ on $M$ is semisimple,
and $M_\lam=0$ unless $\lan h_i,\lam\ran-a\in\BZ$.

\Lemma \label{lemma:nilp}
For $a\in\BC$, $M\in\Ob(\BM_i^a(\Gg))$ and
$u\in M$, we have $f_i^{a+n}\otimes u\in T_i(a)(M)$
for $n\gge0$.
\enlemma
\proof
We may assume that $u$ has weight $\lam$ such that
$\lan h_i,\lam\ran=a$ without loss of generality.
We have $e_i^{m_0}u=0$ for some $m_0>0$.
We have
\eqn
e_i^m(f_i^{a+n}\otimes u)
&=&
\sum_{k=0}^mf_i^{a+n-k}e_i^{m-k}(k!)^2{m\choose k}{a+n\choose k}
{h_i+m-n-a\choose k}\otimes u\\
&=&
\sum_{k=0}^mf_i^{a+n-k}(k!)^2{m\choose k}{a+n\choose k}
{m-n\choose k}\otimes e_i^{m-k}u.
\endeqn
Assume $m\ge n\ge m_0$.
Then each term survives only when $k\le m-n$ and 
$m-k<m_0$, or equivalently
$m-m_0<k\le m-n$, and there is no such $k$.
Hence $e_i^m(f_i^{a+n}\otimes u)=0$ for $m\ge n\ge m_0$.
\qed

\medskip
For $a\in\BC$, the functor $T_i(a)$ sends
$\BM_i^{a}(\Gg)$ to $\BM_i^{-a}(\Gg)$.

The morphism of $\U$-bimodules
\eq\label{mor:coa}
\U\to \U f_i^{-a+\BZ}\otimes_{\U}\U f_i^{a+\BZ}
\qquad(1\mapsto f_i^{-a}
\otimes f_i^{a})
\endeq
(see Lemma \ref{lem:comp})
induces a morphism of functors
\eq\label{mor:inv}
\id_{\BM_i^{a}(\Gg)}\to T_i(-a)\circ T_i(a).
\endeq
Indeed, for $M\in \Ob(\BM_i^{a}(\Gg))$,
(\ref{mor:coa}) gives a morphism
\eqn
M\to \U f_i^{-a+\BZ}\otimes_{\U}\U f_i^{a+\BZ}\otimes_\U M.
\endeqn
For any $u\in M$,
the image of $u$ by the above homomorphism
is equal to $f_i^{-a-n}\otimes 
f_i^{a+n}\otimes u$,
and Lemma \ref{lemma:nilp} implies that
$f_i^{a+n}\otimes u$ belongs to $T_i(a)(M)$
for $n\gge0$.
Hence the image of the above homomorphism is contained in 
$T_i(-a)\circ T_i(a)(M)\subset 
\U f_i^{-a+\BZ}\otimes_{\U}T_i(a)(M)$.

%
%
%
%
%
%
%

\medskip
Define the ideal $\Gn^-(i)$ of $\Gn^-$ by
\eq
\Gn^-(i)=\Gn(\Delta^-\setminus\{-\alpha_i\}).
\endeq
By the PBW theorem, we have $U(\Gn^-)=U(\Gn^-(i))\otimes\C[f_i]$,
which implies
$$U(\Gn^-(i))\otimes \C[f_i,f_i^{-1}]f_i^a\otimes U(\Gb)
\iso \U f_i^{a+\BZ}\,.$$
The following lemma follows immediately from this isomorphism.

\Lemma
\label{KM:bim:lem:trianular}
For any $\lam\in\Gh^*$, we have an isomorphism
$$U(\Gn^-(i))\otimes \C[f_i,f_i^{-1}]f_i^a\otimes\BC u_\lam
\iso \U f_i^{a+\BZ}\otimes_{\U}M(\lam).$$
\enlemma

\Lemma
For any $\lam\in\Gh^*$, the element 
$f_i^{\lan h_i,\lam\ran+1}\otimes u_\lam$
of $\;\U f_i^{\lan h_i,\lam\ran+\BZ}\otimes_{U(\Gg)}M(\lam)$
is a highest weight vector with weight
$s_i\circ\lam$.
Here $\circ$ is the shifted action defined in $(\ref{def;shift})$.
\enlemma

\proof
Set $\lam_i=\lan h_i,\lam\ran$.
We have for $j\not=i$
\eqn
e_j(f_i^{\lam_i+1}\otimes u_\lam)
=f_i^{\lam_i+1}e_j\otimes u_\lam=
f_i^{\lam_i+1}\otimes e_j u_\lam=0.
\endeqn
If $j=i$, then
\eqn
e_i(f_i^{\lam_i+1}\otimes u_\lam)
=(f_i^{\lam_i+1}e_i+(\lam_i+1)f_i^{\lam_i}(h_i-\lam_i))\otimes u_\lam
=0.
\endeqn
\qed


\Prop\label{prop:enr}
Assume that $\lam\in\Gh^*$ satisfies
$a=\lan h_i,\lam+\rho \ran\not\in \BZ_{>0}$.
Then we have
\eq
T_i(a)(M(\lam))
=\U(f_i^{\lan h_i,\lam+\rho\ran}\otimes u_\lam)\cong M(s_i\circ\lam).
\nn
\endeq
\enprop

\proof
We have
\[
T_i(a)(M(\lam))
=
\{u\in \U f_i^{a+\BZ}\otimes_{\U} M\,;\,
\mbox{
$e_i^mu=0$ for a sufficiently large $m$}\}.
\]
By the preceding lemma, $f_i^{a}\otimes u_\lam$
is a highest weight vector of $T_i(a)(M(\lam))$ .
It is enough to show that $T_i(a)(M(\lam))$ is generated by this vector.

By Lemma \ref{KM:bim:lem:trianular},
any $v\in T_i(a)(M(\lam))$
can be written in a unique way
$$v=\sum_{n\in\BZ}P_nf_i^{a+n}\otimes u_\lam$$
for $P_n\in U(\Gn^-(i))$. Here $P_n$ vanishes except for finitely many $n$.
\hb
Take a positive integer $m$ such that $e_i^mv=0$.
Then we have
\eqn
\quad0&=&e_i^mv\\
&=&\sum_n\sum_{k=0}^m{m\choose k}(\ad(e_i)^{m-k}P_n)
e_i^kf_i^{a+n}\otimes u_\lam\\
&=&\sum_n\sum_{k=0}^m{m\choose k}(\ad(e_i)^{m-k}P_n)
\sum_{\nu=0}^kf_i^{a+n-\nu}
e_i^{k-\nu}\\
&&\phantom{aaaaaaaaa}(\nu!)^2{k\choose\nu}{a+n\choose\nu}
{h_i+k-a-n\choose \nu}\otimes u_\lam\\
&=&\sum_n\sum_{k=0}^m{m\choose k}(\ad(e_i)^{m-k}P_n)f_i^{a+n-k}
{k-1-n\choose k}(k!)^2{a+n\choose k}
\otimes u_\lam.
\endeqn
Rewriting this equality, we have
\eqn
\quad0&=&\sum_n\sum_{k=0}^m{m\choose k}(\ad(e_i)^{m-k}P_n)f_i^{n-k}
{k-1-n\choose k}(k!)^2{a+n\choose k}\\
&=&
\sum_n\sum_{k=0}^m{m\choose k}(\ad(e_i)^{m-k}P_{n+k})f_i^{n}
{-1-n\choose k}(k!)^2{a+n+k\choose k}.
\endeqn
The vanishing of the coefficient of $f_i^n$ implies
\eq\label{KM:bim:eqn:1}
\sum_{k=0}^m{m\choose k}
{-1-n\choose k}(k!)^2{a+n+k\choose k}
(\ad(e_i)^{m-k}P_{n+k})=0
\endeq
for any $n$.

Now we shall prove that $P_n=0$ for $n<0$.
Assuming the contrary we take the largest $c>0$ such that
$P_{-c}\not=0$.
By taking $n=-c-m$ in (\ref{KM:bim:eqn:1}),
only $k=m$ survives, and we obtain
$${-1-n\choose m}(m!)^2{a+n+m\choose m}
P_{-c}=0.$$
Hence we obtain
$${-1-c-m\choose m}{a-c\choose m}=0.$$
Since $-1-c-m<0$, ${-1-c-m\choose m}$ does not vanish, and
${a-c\choose m}$ must vanish.
This means that $a-c$ is an integer and satisfies
$0\le a-c<m$.
This leads to the contradiction $a\ge c>0$.
Hence we have $P_n$=0 for $n<0$, and we conclude
$v\in U(\Gn^-)f_i^{a}\otimes u_\lam
= U(\Gn^-)(f_i^{a}\otimes u_\lam)$.
\qed

\medskip

Proposition \ref{prop:enr} implies the following proposition.
\begin{proposition}
Assume $a\equiv\lan h_i,\lam\ran\not\equiv0\ \mod\,\BZ$.
Then the morphism $(\ref{mor:inv})$ induces an isomorphism
\[
M(\lam)\iso
T_i(-a)\circ T_i(a)(M(\lam)).
\]
\end{proposition}

Now we are ready to prove the following proposition used later.
\Prop\label{KM:bim:prop:HomMM}
Assume that $\lam$, $\mu\in\Gh^*$ satisfy
$\lan h_i,\lam\ran\not\in\BZ$.
Then we have
$$\Hom(M(s_i\circ\mu),M(s_i\circ\lam))\simeq\Hom(M(\mu),M(\lam)).$$
\enprop
\proof If $\lam-\mu$ is not in the root lattice $Q$,
then the both sides vanish.
If $\lam-\mu\in Q$, then
$\lan h_i,\mu\ran\equiv\lan h_i,\lam\ran\not\equiv0\ \mod\,\BZ$.
Hence the assertion follows from the preceding proposition.
\qed

\subsection{Embeddings of Verma modules}
We shall use the following result of Kac-Kazhdan.

\Theorem[\cite{KK}]\label{thm:KK}
Let $\lam,\mu\in \Gh^*$.
Then the following three conditions are equivalent.
\begin{description}
\ritem{(i)}
The irreducible highest weight module $L(\mu)$
with highest weight $\mu$ appears as a subquotient of $M(\lam)$.
\ritem{(ii)}
There exist a sequence of
positive roots $\{\beta_k\}_{k=1}^l$,
a sequence of positive integers $\{n_k\}_{k=1}^l$
and a sequence of weights $\{\lam_k\}_{k=0}^l$
such that
$\lam_0=\lam$, $\lam_l=\mu$ and
$\lam_k=\lam_{k-1}-n_k\beta_k$,
$2(\beta_k,\lam_{k-1}+\rho)=n_k(\beta_k,\beta_k)$
for $k=1,\ldots,l$.
\ritem{(iii)}
There exists a non-zero homomorphism $M(\mu)\to M(\lam)$.
\end{description}
\entheorem
Note that any non-zero homomorphism from a Verma module 
to another Verma module must be a monomorphism.
The implication (ii)$\Rightarrow$(iii) is not explicitly stated in
Kac-Kazhdan~\cite{KK}.
But it easily follows from
Lemma 3.3 (b) in Kac-Kazhdan~\cite{KK} and (i)$\Leftrightarrow$(ii).

We use also the following result in Kashiwara~\cite{Kpos1}.

\Prop\label{KM:emb:prop:ext}
For $\lam,\mu\in\Gh^*$ and $i\in I$,
we assume $\lan h_i,\mu+\rho\ran\in\BZ_{\ge0}$
$($which implies $M(s_i\circ\mu)\subset M(\mu)$$)$
and $\lan h_i,\lam+\rho\ran\not\in\BZ_{<0}$.
Then
we have
$$\Ext^1_\Gg(M(\mu)/M(s_i\circ\mu),M(\lam))=0.$$
\enprop

\bigskip
Let $\K$ denote the set of $\lam\in\Gh^*$ satisfying 
the following two conditions.
\eq
&&\hbox{$2(\beta,\lam+\rho)\ne(\beta,\beta)$
for any positive imaginary root $\beta$,}
\label{cond:real}\\
&&\hbox{$\{\beta\in\Delta_\re^+; (\beta^\vee,\lam+\rho)\in \BZ_{<0}\}$
is a finite set.}\label{cond:real1}
\endeq

The condition
(\ref{cond:real1})
implies that there exists $w\in W(\lam)$
such that $w\circ\lam+\rho$ is integrally dominant
(i.e.\ $(\beta^\vee,w\circ\lam+\rho)\notin\BZ_{<0}$ for any
$\beta\in\Delta_\re^+$).

If $\lam$ in Theorem~\ref{thm:KK}
satisfies the condition (\ref{cond:real}), then
$\beta_k$ in (ii) must be a real positive root.
This easily follows from the fact that
$n\beta$ is an imaginary root for any positive integer $n$
and any imaginary root $\beta$.

Note that $\K$ is invariant by the shifted action of $W$.

\Theorem\label{KM:emb:thm:embed}
For $\lam\in\K$
we have
$$\dim \Hom_\Gg(M(\mu),M(\lam))\le1$$
for any $\mu\in\Gh^*$.
\entheorem

\proof
There exists an embedding $M(\lam)\hookrightarrow M(\lam')$
for some $\lam'\in W(\lam)\circ\lam$
such that $\lam'+\rho$ is integrally dominant.
Hence we may assume that $\lam+\rho$ is integrally dominant from the beginning.

We assume that $\Hom(M(\mu),M(\lam))$ is not zero.
Then by Theorem~\ref{thm:KK},
there exists $w\in W(\lam)$ such that $\mu=w\circ\lam$.

We shall argue by the induction on the lenght of $w$.

%

If $w=1$, then it is evident.
Assuming $w\not=1$, let us take
$\alpha\in\Pi(\lam)$ such that $l_\lam(s_\alpha w)<l_\lam(w)$,
which is equivalent to $w^{-1}\alpha\in\Delta^-(\lam)$.
Since $\lam+\rho$ is integrally dominant,
$(w^{-1}\alpha^\vee,\lam+\rho)\le0$.
Since we may assume $s_\alpha\circ\mu\not=\mu$, we have
\eq
(\alpha^\vee ,\mu+\rho)=(w^{-1}\alpha^\vee,\lam+\rho)\in\BZ_{<0}.
\endeq
Now we shall argue by the induction on $\HT(\alpha)$.
\par\medskip
\noindent
(1) Case $\HT(\alpha)=1$.\quad 
In this case, $\alpha=\alpha_i$ for some $i\in I$.
Then we have
$M(s_i\circ\mu)\supset M(\mu)$.
Since $\lan h_i,\lam+\rho\ran\in \BZ_{\ge0}$,
Proposition~\ref{KM:emb:prop:ext} implies
$$\Ext^1(M(s_i\circ\mu)/M(\mu), M(\lam))=0.$$
Therefore the following sequence is exact.
$$\Hom(M(s_i\circ\mu),M(\lam))\to\Hom(M(\mu),M(\lam))\to0.$$
Since $\dim \Hom(M(s_i\circ\mu),M(\lam))\le1$
by the induction hypothesis on the length of $w$,
we obtain  $\dim \Hom(M(\mu),M(\lam))\le1$.
\par\smallskip
\noindent
(2) Case $\HT(\alpha)>1$.\quad Take $i$ such that
$\lan h_i,\alpha\ran>0$. Then $\alpha_i\not\in\Delta(\lam)$ 
by Lemma \ref{lemma:simple}.
Hence we have
\eq\label{eq:hyp}
\lan h_i,\lam\ran\not\in \BZ.
\endeq
Set $\lam'=s_i\circ\lam$.
Then $s_i\Delta(\lam)=\Delta(\lam')$, 
and $s_i\Pi(\lam)=\Pi(\lam')$ by Lemma \ref{km:int:lem:2}.
Moreover $\lam'+\rho$ is also integrally dominant.
Then $w'=s_iws_i\in W(\lam')$ and
$l_{\lam'}(w')=l_\lam(w)$.
Set $\alpha'=s_i\alpha$.
Then $\alpha'\in\Pi(\lam')$, $\HT(\alpha')<\HT(\alpha)$ and
$l_{\lam'}(s_{\alpha'}w')<l_{\lam'}(w')$.
We have also $\mu'=s_i\circ\mu=w'\circ\lam'$.
Hence the induction hypothesis on $\HT(\alpha)$
implies
$$\dim\Hom(M(\mu'),M(\lam'))\le1.$$
By (\ref{eq:hyp}) we can apply Proposition~\ref{KM:bim:prop:HomMM}
to deduce
$$\Hom(M(\mu'),M(\lam'))\simeq\Hom(M(\mu),M(\lam)).$$
Thus we obtain the desired result
$\dim\Hom(M(\mu),M(\lam))\le1$.
\qed

\medskip
We denote by $\KR$ the set of $\lambda\in\K$ 
subject to the following condition:
\eq
&&\label{wt:st}
\mbox{If $w\in W$ satisfies $w\circ\lam=\lam$, then $w=1$.}
\endeq
In particular, this condition implies
\eqn
&&\mbox{
$(\lambda+\rho,\alpha^\vee)\ne0$ for any $\alpha\in\Delta_{\rm re}$.
}
\endeqn
Define a subset $\KRP$ of $\KR$ by
\eq&&
\KRP=
\{\lambda\in\KR\, ;\,(\lambda+\rho,\alpha^\vee)>0
\quad\mbox{for any $\alpha\in\Delta^+(\lam)$$\}$.}
\endeq

\begin{lemma}
We have $W\circ\KR=\KR$ and
$\KR=\bigsqcup\limits_{\lambda\in\KRP}W(\lambda)\circ\lambda$.
\end{lemma}
The proof is standard by using the results in \S \ref{sec:weyl}
and omitted.

\bigskip
By Theorem~\ref{thm:KK} and  Theorem~\ref{KM:emb:thm:embed} we have the
following proposition.
\begin{proposition}
\label{D:mod:prop:[M:L]}
Let $\lam\in \KRP$.
\begin{tenumerate}
\item
For $x\in W(\lam)$ and $\mu\in\Gh^*$ we have
$[M(x\circ\lam):L(\mu)]\ne0$ if and only if $\mu=y\circ\lam$ for some $y\in
W(\lam)$ satisfying $y\ge_\lam x$.
\item
For $x,y\in W(\lam)$ we have
$\dim\Hom(M(y\circ\lam),M(x\circ\lam))=1$ or $0$
according to whether $y\ge_{\lam}x$ or not.
\end{tenumerate}
\end{proposition}
\begin{corollary}
\label{D:mod:cor:[M:L]}
Let $\lam\in \KRP$,
$x\in W$ and $\mu\in\Gh^*$.
Then $[M(x\circ\lam):L(\mu)]\ne0$ implies
$\mu=y\circ\lam$ for some $y\in xW(\lam)$
satisfying $y\ge x$.
\end{corollary}
\proof
Assume $[M(x\circ\lam):L(\mu)]\ne0$.
Take $z_1\in W(x\circ\lam)$ such that $\lam'=z_1^{-1}x\circ\lam\in\KRP$.
Then we have $z_1\in W(\lam')$ and $x\circ\lam=z_1\circ\lam'$.
By Proposition~\ref{D:mod:prop:[M:L]} there exists some $z_2\in W(\lam')$
such that $\mu=z_2\circ\lam'$ and $z_2\geq_{\lam'}z_1$.
Setting $w=z_1^{-1}x,\; y=z_2w$ we have $x=z_1w$ and $\mu=y\circ\lam$.
Since $y=z_2z_1^{-1}x\in W(\lam')x=xW(\lam)$,
the assertion follows from the following lemma.
\begin{lemma}\label{lemma:ord}
Assume that $\lam, \lam'\in\KRP$ and $w\in W$ satisfy
$\lam'=w\circ\lam$.
Then for $z_1, z_2\in W(\lam')$ such that $z_2\geq_{\lam'} z_1$ we have
$z_2w\geq z_1w$.
\end{lemma}
\proof
For $\lam\in\KRP$, we have
$\Delta^+(\lam)=\{\alpha\in\Delta_\re;(\alpha^\vee,\lam+\rho)\in\BZ_{>0}\}$.
This implies
$w^{-1}\Delta^+(\lam')=\Delta^+(\lam)\subset\Delta^+$.
Then it is enough to apply Lemma \ref{prop:cong}.
\qed

For a subset $\Omega$ of $\KR$ we denote by $\BO\{\Omega\}$ the full
subcategory of $\BO$ consisting of $M\in\Ob(\BO)$ such that any irreducible
subquotient of $M$ is isomorphic to $L(\lam)$ for some $\lam\in\Omega$.
For $\lam\in\KR$ we set
\eq\label{def:o}
&&\BO[\lam]=\BO\{W(\lam)\circ\lam\},\qquad
\BO(\lam)=\BO\{W\circ\lam\}.
\endeq
By the definition, for any $\lam\in\KR$, we have
\eq
&&\BO[\lam]=\BO[w\circ\lam]\;
\mbox{
for any $w\in W(\lam)$,
}
\\
&&\BO(\lam)=\BO(w\circ\lam)\;
\mbox{for any $w\in W$.}
\endeq
By Proposition~\ref{D:mod:prop:[M:L]} we have
\eq
M(\lam)\in\Ob(\BO[\lam])\;
\mbox{
for any $\lam\in\KR$.
}
\endeq

By Lemma~\ref{KM:cat:dimHom} and Corollary~\ref{D:mod:cor:[M:L]} we have the
following lemma.
\begin{lemma}
\label{KM:cat:dimHom2}
Let $\lam\in\KRP$ and $w\in W$.
Assume that $M\in\Ob(\BO)$ 
satisfies the conditions
\eqn
&&[M:L(w\circ\lam)]\ne0,\\
&&\mbox{$[M:L(y\circ\lam)]=0$ for any $y\in wW(\lam)$ such that $y<w$.}
\endeqn
Then neither $\Hom_\Gg(M(w\circ\lam),M)$ nor
$\Hom_\Gg(M,M^*(w\circ\lam))$ vanishes.
\end{lemma}

\bigskip
We shall use later the following result of S. Kumar~\cite{Kumar} 
(a generalization of a result in
Deodhar-Gabber-Kac~\cite{DGK}).

\begin{theorem}
\label{thm:Kumar}
Any object $M$ of $\BO\{\KR\}$  decomposes uniquely into
\[
M=\bigoplus_{\lam\in\KRP}M^\lam
\quad(M^\lam\in\Ob(\BO[\lam])).
\]
\end{theorem}
In \cite{Kumar}, the theorem is proved for
$M$ with a semisimple action of $\Gh$.
However the same arguments can be applied in our situation.

For $\lam\in\KR$ we denote by
\eq
\label{eq:projection functor}
P_\lam:\BO\{\KR\}\to\BO[\lam]
\endeq
the projection functor.

We define a new abelian category
$\TBO$ by 
\eq\label{def:TBO}
&&\TBO=\prod_{\lam\in\KRP}\BO[\lam].
\endeq
We denote by the same symbol $P_\lam$
the projection functor
$P_\lam:\TBO\to \BO[\lam]$.
It is an exact functor.
By the definition we have
\[
\Hom_{\TBO}(M,N)=\prod_{\lam\in\KRP}\Hom_\Gg(P_\lam(M),P_\lam(N))
\quad\mbox{for $M$, $N\in\TBO$.}
\]
The category $\BO\{\KR\}$ can be regarded as a full subcategory
of $\TBO$.
For $M\in \TBO$ and $\lam\in\KR$ we set
\eqn&&
[M:L(\lam)]=[P_\lam(M):L(\lam)]\,.
\endeqn
For a subset $\Omega$ of $\KR$, we set
\eq
&&\TBO\{\Omega\}=
\prod_{\lam\in\KRP}\BO\{\Omega\cap(W(\lam)\circ\lam)\}\,,
\endeq
and for $\lam\in\KR$,
\eq
&&\TBO(\lam)=\TBO\{W\circ\lam\}\,.
\endeq

\newsection{Twisted $D$-modules}
We shall give a generalization of the theory of $D$-modules on
infinite-dimensional schemes developed in \cite{Kpos1} and \cite{KTpos2} to
that of twisted left $D$-modules (modules over a TDO-ring).
Since the arguments are analogous to the original non-twisted case, we
only state the results and omit proofs.

\subsection{Finite-dimensional case}
For a scheme $X$ we denote by $\CO_X$ the structure sheaf.
For a scheme $X$ smooth (in particular quasi-compact and separated)
over $\BC$, 
we denote by $\Omega_X$, $\Theta_X$ and $D_X$ 
the canonical sheaf, the sheaf of vector fields, and the sheaf of
rings of differential operators on $X$, respectively.

Let $X$ be a scheme smooth over $\BC$.
A TDO-ring on $X$ is by definition a sheaf $A$ of rings on $X$ containing
$\CO_X$ as a subring such that there exists an increasing filtration
$F=\{F_nA\}_{n\in\BZ}$ of the abelian sheaf $A$ satisfying the following
conditions.
\eq
&&F_nA=0\;\mbox{ for }n<0.\label{eq:TDO2}
\\
&&
F_nA\cdot F_mA\subset F_{n+m}A.\label{eq:TDO4}
\\
&&
[F_nA,F_mA]\subset F_{n+m-1}A.\label{eq:TDO5}
\\
&&
F_0A=\CO_X.\label{eq:TDO6}
\\
&&
\mbox{
The homomorphism
$\gr_1A\to\Theta_X\;
(P\:\mod\: F_0A\mapsto(\CO_X\ni a\mapsto[P,a]\in\CO_X))$ of
}
\\
&&
\mbox{
$\CO_X$-modules induced by (\ref{eq:TDO5}),  (\ref{eq:TDO6}) is an isomorphism.
}
\label{eq:TDO7}
\nn
\\
&&
\mbox{The homomorphism
$S_{\CO_X}(\gr_1A)\to\gr A$ of commutative
$\CO_X$-algebras is an}
\\
&&
\mbox{isomorphism.}
\nn
\endeq
Here we set $\gr_nA=F_nA/F_{n-1}A$, $\gr A=\bigoplus_n\gr_nA$, and
$S_{\CO_X}(\gr_1A)$ denotes the symmetric algebra of the locally free
$\CO_X$-module $\gr_1A$.
The filtration $F$ is uniquely determined by the above conditions, and it is
called the order filtration.
A TDO-ring is quasi-coherent over $\CO_X$ with respect to its left and right
$\CO_X$-module structures.

Let $A$ be a TDO-ring on a scheme $X$ smooth over $\BC$.
For a coherent (left) $A$-module $\CM$ we can define its characteristic
variety $\Ch(M)$ as a subvariety of the cotangent bundle $T^*X$ as in the
case $A=D_X$.
A coherent $A$-module $\CM$ is called holonomic if
$\dim\Ch(\CM)\leq\dim X$.
We denote by $\BM_h(A)$ the category of holonomic $A$-modules, and by
$D^b_h(A)$ the derived category consisting of bounded complexes of
quasi-coherent $A$-modules with holonomic cohomologies.
Set
\[
A^{-\sharp}=
\Omega_X^{\otimes-1}
\otimes_{\CO_X}
A^{\rm op}
\otimes_{\CO_X}
\Omega_X,
\]
where $A^{\rm op}$ denotes the opposite ring of $A$.
Then $A^\sharp$ is also a TDO-ring.
We define the duality functor
\[
\dual:D^b_h(A)\to D^b_h(A^{-\sharp})^{\rm op}
\]
by
\[
{\bf D}\CM=
\BR\sHom_A(\CM,A)
\otimes_{\CO_X}
\Omega_X^{\otimes-1}[\dim X].
\]

Let $f:X\to Y$ be a morphism of smooth schemes over $\BC$, and let $A$ be a
TDO-ring on $Y$.
Set
\[
A_{X\to Y}=\CO_X\otimes_{f^{-1}\CO_Y}f^{-1}A,\qquad
A_{Y\leftarrow X}
=f^{-1}A
\otimes_{f^{-1}\CO_Y}
f^{-1}\Omega_Y^{\otimes-1}
\otimes_{f^{-1}\CO_Y}
\Omega_X.
\]
Define the subring $f^\sharp A$ of $\sEnd_{f^{-1}A}(A_{X\to Y})$ by
\eq
&&
f^\sharp A=\bigcup_{n\in\BN}F_n(f^\sharp A),
\label{TDOpullback1}
\\
&&
F_n(f^\sharp A)=0
\mbox{ for }
n<0,
\label{TDOpullback2}
\\
&&
F_n(f^\sharp A)=
\{P\in\sEnd_{f^{-1}A}(A_{X\to Y})\,;\,
[P,\CO_X]\subset F_{n-1}(f^\sharp A)\}
\mbox{ for }
n\geq0.
\label{TDOpullback3}
\endeq
Then $f^\sharp A$ is a TDO-ring on $X$.
Moreover, $A_{X\to Y}$ 
has a structure of an $(f^\sharp A,f^{-1}A)$-bimodule, and
$A_{Y\leftarrow X}$ has a structure of an $(f^{-1}A,f^\sharp A)$-bimodule.
We have
\eq
f^\sharp(A^{-\sharp})=(f^\sharp A)^{-\sharp}
\endeq
for any TDO-ring $A$.
We define functors
\eq
\BD f^\BUL:D^b_h(A)\to D^b_h(f^\sharp A),
\endeq
\eq
\int_f:D^b_h(f^\sharp A)\to D^b_h(A),\qquad
\int_{f!}:D^b_h(f^\sharp A)\to D^b_h(A)
\endeq
by
\eq\label{def:pull}
&&
\begin{array}{l}
\BD f^\BUL(\CM)=A_{X\to Y}\otimes^{\BL}_{f^{-1}A}f^{-1}M,\quad
\int_f\CM=\BR f_*(A_{Y\leftarrow X}\otimes^{\BL}_{f^\sharp A}\CM),\\[5pt]
\int_{f!}=\dual\circ\int_f\circ\dual.
\end{array}
\endeq
We shall also use their variants
\eq
\BD f^*,\; \BD f^!:D^b_h(A)\to D^b_h(f^\sharp A),\qquad
\BD f_*,\; \BD f_!:D^b_h(f^\sharp A)\to D^b_h(A)
\endeq
given by
\eqn
&&
\begin{array}{ll}
\BD f^*=\dual\circ\BD f^\BUL\circ\dual,\qquad
&\BD f^!=\BD f^\BUL\ [2(\dim X-\dim Y)],
\\[5pt]
\BD f_*=\int_f\ [\dim Y-\dim X],\qquad&
\BD f_!=\int_{f!}\ [\dim Y-\dim X].
\end{array}
\endeqn

\subsection{Infinite-dimensional case}
Now we shall study TDO-rings on infinite-dimensional varieties.
We say that a scheme $X$ over $\BC$ satisfies the property (S) if
$X\simeq\limp_n S_n$
for some projective system $\{S_n\}_{n\in\BN}$ satisfying the following
conditions.
\eq
&&
\mbox{The scheme $S_n$ is smooth 
(in particular quasi-compact and separated)
}\label{property S1}\\
&&\mbox{over $\BC$ for any $n$.}\nn\\
&&\mbox{
The morphism $p_{nm}:S_m\to S_n$ is affine and smooth for any $m\geq n$.
}
\label{property S2}
\endeq
We call such $\{S_n\}_{n\in\BN}$ a {\em smooth projective system} for $X$.
For example, the infinite-dimensional affine space
\[
\BA^\infty=\limp_n\BA^n=\Spec \BC[x_i;i\in\BN]
\]
satisfies the property (S).

Let ${\cal S}$ denote the category whose objects are 
smooth $\BC$-schemes
and whose morphisms are affine and smooth morphisms.
Then the pro-object $``{\limp_n}$''$S_n$ in ${\cal S}$ 
depends only on $X$ and does not depend on the choice of 
a smooth projective system $\{S_n\}_{n\in\BN}$ (\cite{EGA}).
This follows from the fact
\eqn
&&\Hom(X,S)=\limi_n\Hom(S_n,S)
\endeqn
for any scheme $S$.

A $\BC$-scheme $X$ is called {\em pro-smooth} if it is covered by open subsets
satisfying (S).

Let $f:X\to Y$ be a morphism of $\BC$-schemes such that $X$ is pro-smooth
and $Y$ is smooth over $\BC$, and let $A$ be a TDO-ring on $Y$.
Set $A_{X\to Y}= \CO_X\otimes_{f^{-1}\CO_Y}f^{-1}A$, and define the subring
$f^\sharp A$ of $\sEnd_{f^{-1}A}(A_{X\to Y})$ by
(\ref{TDOpullback1})--(\ref{TDOpullback3}).
Then $A_{X\to Y}$ is an $(f^\sharp A,f^{-1}A)$-bimodule.
For an $A$-module $\CM$ we define the $f^\sharp A$-module $f^\BUL\CM$ by
\eq
f^\BUL\CM=A_{X\to Y}\otimes_{f^{-1}A}f^{-1}\CM.
\endeq
For a pro-smooth scheme $X$, a TDO-ring on $X$ is by definition a sheaf $A$
of rings on $X$ containing $\CO_X$ as a subring satisfying the following
condition.
\eq
\label{TDOonPROsmooth}
&&
\mbox{
\quad For any $x\in X$, 
there exist a morphism $f:U\to Y$ from an open neighborhood
}\\
&&
\mbox{\quad
$U$ of $x$ to a smooth $\BC$-scheme $Y$ and a TDO-ring $B$ on $Y$ such that
$A|_U=f^\sharp B$.
}
\nn
\endeq
By the definition, a TDO-ring $A$ on a pro-smooth scheme $X$ is locally of
the form $A=f^\sharp B$ where $B$ is a TDO-ring on a smooth $\BC$-scheme.
We can patch together $f^\sharp B^{-\sharp}$ and obtain a TDO-ring
$A^{-\sharp}$ on $X$.

For an invertible $\CO_X$-module $\CL$ on a pro-smooth scheme $X$, we have a
TDO-ring $D_X(\CL)$ given as follows.
\eq
&&
D_X(\CL)=\bigcup_nF_nD_X(\CL)\subset\sEnd_\BC\CL.
\\
&&
F_nD_X(\CL)=0\;\mbox{ for }n<0.
\\
&&
F_nD_X(\CL)=
\{P\in\sEnd_\BC\CL\,;\,
[P,a]\in F_{n-1}D_X(\CL)\;
\mbox{ for any }
a\in\CO_X\}
\;\mbox{ for }n\geq0.
\endeq
We set $D_X=D_X(\CO_X)$.
Then we have
$D_X(\CL)\simeq\CL\otimes_{\CO_X}D_X\otimes_{\CO_X}\CL^{\otimes-1}$.
More generally, for an invertible $\CO_X$-module $\CL$ and a scalar
$a\in\BC$ we can define a TDO-ring
$D_X(\CL^a)=\CL^a\otimes_{\CO_X}D_X\otimes_{\CO_X}\CL^{\otimes-a}$ by the
following patching procedure although $\CL^a$ does not necessarily exist.
A section of $D_X(\CL^a)$ is locally of the form $s^a\otimes P\otimes
s^{-a}$, where $s$ is a nowhere vanishing section of $\CL$ and $P$ is a
section of $D_X$, and we have
$s_1^a\otimes P_1\otimes s_1^{-a}=s_2^a\otimes P_2\otimes s_2^{-a}$ if and
only if $P_1=(s_2/s_1)^aP_2(s_2/s_1)^{-a}$ as sections of $D_X$.

Let $A$ be a TDO-ring on a pro-smooth scheme $X$.
We call a (left) $A$-module $\CM$ admissible if it satisfies the following
conditions.
\eq
&&
\ \mbox{
$\CM$ is quasi-coherent over $\CO_X$.
}
\\
&&
\ \mbox{
For any affine open subset $U$ of $X$ and any $s\in\Gamma(U;\CM)$, there
exists a finitely}
\\
&&
\ \mbox{
generated $\BC$-subalgebra $B$ of $\Gamma(U;\CO_X)$ such that $Ps=0$
for any $P\in\Gamma(U;A)$}
\nn\\
&&
\ \mbox{
satisfying $P(B)=0$.
}
\nn
\endeq
We denote the category of admissible $A$-modules by $\MA(A)$.
We call an admissible $A$-module $\CM$ {\em holonomic} if it satisfies the
following condition.
\eq
&&
\mbox{For any $x\in X$ there exist 
a morphism $f:U\to Y$, a TDO-ring $B$ on $Y$}
\\
&&
\mbox{as in
(\ref{TDOonPROsmooth}) and a holonomic $B$-module $\CM'$ 
such that $\CM|_U=f^\BUL\CM'$.}\nn
\endeq
We denote the category of holonomic $A$-modules by $\BM_{h}(A)$.

Let $A$ be a TDO-ring on a $\BC$-scheme $X$ satisfying the property (S).
Then we can take a smooth projective system $\{S_n\}_{n\in\BN}$ for $X$ and
TDO-rings $A_n$ on $S_n$ such that
\eq
p_{nm}^\sharp A_n\cong A_m
\mbox{ for any $m\geq n$, }\qquad
p_{n}^\sharp A_n\cong A,
\endeq
where $p_n:X\to S_n$ is the projection.
We call $\{(S_n,A_n)\}_{n\in\BZ}$ 
a smooth projective system for $(X,A)$.
Since $p_{nm}$ is smooth, the functor
$p_{nm}^\BUL:\BM_{h}(A_n)\to\BM_{h}(A_m)$ is exact, and we have
the equivalence of categories
\[
\BM_{h}(A)\cong\limi_n\BM_{h}(A_n).
\]

Let $X$ be a pro-smooth $\BC$-scheme and $A$ a TDO-ring on $X$.
Let $D(\MA(A))$ be the derived category of
$\MA(A)$.
Let us denote by $D^b_h(A)$ the full subcategory of $D^b(\MA(A))$
consisting of bounded complexes  whose cohomology groups
are holonomic.

If $X$ satisfies (S) and 
$\{(S_n,A_n)\}_{n\in\BZ}$ is a smooth projective system for $(X,A)$,
then we have an equivalence of categories:
\[
D^b_h(A)=\limi_nD^b_h(A_n).
\]
The duality functors $\dual:D^b_h(A_n)\to D^b_h(A_n^{-\sharp})^{\rm op}$
induce the duality functor
\eq
\dual:D^b_h(A)\to D^b_h(A^{-\sharp})^{\rm op}.
\endeq

For a morphism $f:X\to Y$ of pro-smooth $\BC$-schemes and a TDO-ring $A$
on $Y$, we define a TDO-ring $f^\sharp A$
on $X$ by the same formulas (\ref{TDOpullback1})--(\ref{TDOpullback3}).
The functor
\eq
\BD f^\BUL:D^b_h(A)\to D^b_h(f^\sharp A)
\endeq
is defined by the same formula as in (\ref{def:pull}).
It is well-defined as seen in the following.
The question being local, we may assume that $X$ and $Y$ satisfy (S).
Then we can take a smooth projective system $\{X_n\}$ for $X$ 
and a smooth projective system $\{(Y_n,A_n)\}_{n\in\BZ}$ for $(Y,A)$.
Let $p_{Xn}:X\to X_n$ and $p_{Yn}:Y\to Y_n$ be the projections.
We may assume further that there exists $\{f_n\}:\{X_n\}\to\{Y_n\}$ such that
$f=\limp_nf_n$.
For $\CM\in\Ob(D^b_h(A))$ there exist some $n$ and
$\CM_n\in\Ob(D^b_h(A_n))$ such that $\CM=p_{Yn}^\BUL\CM_n$.
Then we have
\[
\BD f^\BUL\CM=p_{Xn}^\BUL\BD f_n^\BUL\CM_n.
\]

Let $f:X\to Y$ be a morphism of pro-smooth schemes.
We assume that $f$ is of finite presentation type.
Let $A$ be a TDO-ring on $Y$.
We define a functor
\eq
\label{functor1}
\int_f:D^b_h(f^\sharp A)\to D^b_h(A)
\endeq
by the same formula as in (\ref{def:pull}).
We can see that it is well-defined as follows.
The question being local on $Y$, 
we can take a smooth projective system $\{X_n\}$ for $X$, a smooth projective
system $\{(Y_n,A_n)\}_{n\in\BZ}$ for $(Y,A)$, 
and $\{f_n\}:\{X_n\}\to\{Y_n\}$ such
that $f=\limp_nf_n$ and the following diagram is Cartesian for any $n$
(\cite{EGA}).
\[
\begin{array}{ccc}
X_n
&
\maprightu{f_n}
&
Y_n\phantom{\,.}
\\
\mapdownl{}
&
&
\mapdownl{}
\\
X_0
&
\maprightu{f_0}
&
Y_0\,.
\end{array}
\]
Let $p_{Xn}:X\to X_n$ and $p_{Yn}:Y\to Y_n$ be the projections.
Let $\CM\in\Ob(D^b_h(f^\sharp A))$.
There exists some $n$ and $\CM_n\in\Ob(D^b_h(f_n^\sharp A_n))$ such that
$\CM=p_{Xn}^\BUL\CM_n$.
Then we have
\[
\int_f\CM=p_{Yn}^\BUL\int_{f_n}\CM_n.
\]
%
Under the same assumption, we define the relative dimension
$d_f$ by $\dim X_0-\dim Y_0$.

\medskip
We shall also use the following functors
for a morphism $f:X\to Y$ of schemes satisfying (S):
\eq\label{functor2}
&&
\BD f^*,\ \BD f^!:D^b_h(A)\to D^b_h(f^\sharp A),\qquad
\int_{f!}\,,\ \BD f_*\,, \ \BD f_!:D^b_h(f^\sharp A)\to D^b_h(A)
\endeq
defined by
\[
\BD f^!=\BD f^\BUL[2d_f],\quad
\BD f^*=\dual\circ\BD f^\BUL\circ\dual,
\]
\[
\int_{f!}=\dual\circ\int_f\circ\dual,\quad
\BD f_*=\int_{f}[-d_f],\quad
\BD f_!=\int_{f!}[-d_f].
\]
Note that $\BD f^\bullet$ and $\BD f^*$ are defined for any morphism $f$ of
schemes satisfying (S), while other functors in (\ref{functor1}),
(\ref{functor2}) are defined only when $f$ is of finite presentation type.

For a morphism $f:X\to Y$ of pro-smooth schemes, a TDO-ring $A$ on $Y$, and
$k\in\BZ$ we can define functors
\eq
H^k\BD f^*:\BM_h(A)\to \BM_h(f^\sharp A)
\endeq
by patching together the locally defined object
$H^k(\BD f^*\CM)$ for $\CM\in\Ob(\BM_h(A))$.
Similarly, for a morphism $f:X\to Y$ of pro-smooth schemes which is
of finite presentation type, a TDO-ring $A$ on $Y$, and $k\in\BZ$ we can
define functors
\eq
&&
\begin{array}{l}
H^k\int_{f!}\,,\ H^k\BD f_!\,:\,\BM_h(f^\sharp A)\to
\BM_h(A).
\end{array}
\endeq

\subsection{Equivariant $D$-modules}\label{sec:equiv}
Let $G$ be an affine group scheme over $\BC$.
We assume that $\CO_G(G)$ is generated by countably many generators
as a $\BC$-algebra.
Then $G\simeq\limp_{n\in\BN}G_n$
for a projective systems of affine algebraic group over $\BC$,
and hence $G$ satisfies the condition (S). Let $\Gg$ be the Lie algebra of $G$.
Then 
$\Gg$ is the projective limit of the Lie algebras $\Gg_n$ of $G_n$.
\hb
Let $X$ be a pro-smooth $\BC$-scheme with an action of $G$.

Then we have the diagram
\eq
&&
\begin{array}{ccccc}
&\mapr{p_1}\\[-10pt]
&&&\mapr{\mu_X}\\[-4pt]
G\times G \times X&
\mapr{p_2}
&G\times X&\mapl{i}&X\,.\\[-4pt]
&&&\mapr{\pr_X}\\[-10pt]
&\mapr{p_3}
\end{array}
\endeq
Here $\mu_X:G\times X\to X$ is the action morphism,
$\pr_X: G\times X\to X$ the second projection,
and
\eqn
i(x)&=&(1,x),\\
p_1(g_1,g_2,x)&=&(g_1,g_2x),\\
p_2(g_1,g_2,x)&=&(g_1g_2,x),\\
p_3(g_1,g_2,x)&=&(g_2,x).
\endeqn
A $G$-equivariant TDO-ring $A$ on $X$
is a TDO-ring endowed with an isomorphism of TDO-rings
\[\alpha:\mu_X^\sharp A\iso\pr_X^\sharp A\]
with the cocycle condition (see \cite[\S 4.6]{Krep}),
i.e. the commutativity of the following diagram.
\eq
\begin{array}{ccc}
p_2^\sharp\mu_X^\sharp A&
\mathop{\hbox to 120pt{\rightarrowfill}}\limits^{p_2^\sharp\alpha}
&p_2^\sharp\pr_X^\sharp A\\
\downeq&&\downeq\\
p_1^\sharp\mu_X^\sharp A&\maprightu{p_1^\sharp\alpha}
p_1^\sharp\pr_X^\sharp A\cong p_3^\sharp\mu_X^\sharp A
\maprightu{p_3^\sharp\alpha}
&p_3^\sharp\pr_X^\sharp A
\end{array}
\endeq
Then the $G$-equivariance structure induces a ring homomorphism
\eqn
&&U(\Gg)\to \Gamma(X;A).
\endeqn
A $G$-equivariant module $M$ 
over a $G$-equivariant TDO-ring $A$ on $X$
is an $A$-module endowed with an isomorphism of $\pr_X^\sharp A$-modules
\[\mu_X^\BUL M\iso\pr_X^\BUL M\]
with a similar cocycle condition 
(see \cite[\S 4.7]{Krep} and (\ref{mod:equiv}) below).

We can generalize the notion of equivariance to that of
twisted equivariance.
Assume for the sake of simplicity that
$G$ is a finite-dimensional affine algebraic group with Lie algebra
$\Gg$.
Let $q_i:G\times G\to G$ ($i=1,2$) be the first and the second projection
and $\mu_G:G\times G\to G$ the multiplication morphism.
Let $i_G:\pt\to G$ be the identity.

Let $\lam\in\Gg^*$ be a $G$-invariant vector.
Let $\CT(\lam)$ be the free $\CO_G$-module generated by the symbol
$\e^\lam$. We define its $D_G$-module structure by
\eqn
&&R_A\e^\lam=\lam(A)\e^\lam
\quad\mbox{for any $A\in\Gg$,}
\endeqn
where $R_A$ is the left invariant vector field on $G$ corresponding to $A$.
Then we have a canonical isomorphism of
$D_G$-modules
\eq\label{mul.sys}
&&
\ba{l}
i_G^\BUL \CT(\lam)\cong \BC\,,\\[3pt]
m_\lam:\mu_G^\BUL \CT(\lam)
\iso q_1^\BUL \CT(\lam)\otimes q_2^\BUL \CT(\lam) \,,
\ea
\endeq
sending $\e^\lam$ to $1$ and $\e^\lam\otimes \e^\lam$,
respectively.
A twisted $G$-equivariant $A$-module $\CM$ with twist $\lam$
is an $A$-module
with an isomorphism of $\pr_X^\sharp A$-modules
\eq
&&\beta:\mu_X^\BUL \CM\iso q^\BUL \CT(\lam)\otimes \pr_X^\BUL \CM,
\endeq
with the cocycle condition.
Here $q:G\times X\to G$ is the first projection.
The cocycle condition means
the commutativity of the following 
diagram of $r_3^\sharp A$-modules on $G\times G\times X$.
\eq\label{mod:equiv}
&&
\begin{array}{ccc}
q_{12}^\BUL\mu_G^\BUL \CT(\lam)\ot r_3^\BUL\CM
& \mapr{m_\lam}
& q_{12}^\BUL(r_1^\BUL \CT(\lam)\ot r_2^\BUL \CT(\lam)) \ot r_3^\BUL\CM\\
\|&&\raise-6pt\hbox{$\Big\|$}\\
p_2^\BUL(q^\BUL \CT(\lam) \ot\pr_X^\BUL\CM)
& &\raise-2pt\hbox{$\Big\|$}\\
\Big\uparrow\rlap{\ $\scriptstyle \beta$}
&&
\raise0pt\hbox{$\Big\|$}
\\
p_2^\BUL\mu_X^\BUL\CM
&&\raise6pt\hbox{$\Big\|$}
\\
\|
&&
r_1^\BUL \CT(\lam)\ot r_2^\BUL \CT(\lam)\ot r_3^\BUL\CM
\\
p_1^\BUL\mu_X^\BUL\CM
&&\raise-6pt\hbox{$\Big\|$}
\\
\Big\downarrow\rlap{\ $\scriptstyle \beta$}
&&\hbox{$\Big\|$}
\\
p_1^\BUL(q^\BUL \CT(\lam) \ot\pr_X^\BUL\CM)
&&\hbox{$\Big\|$}\\
\|&&\raise6pt\hbox{$\Big\|$}\\
r_1^\BUL \CT(\lam)\ot p_3^\BUL\mu_X^\BUL\CM
& \mapr{\beta}
&
\ r_1^\BUL \CT(\lam) \ot p_3^\BUL(q^\BUL \CT(\lam) \ot \pr_X^\BUL\CM)\,.
\end{array}
\endeq
Here $q_{12}:G\times G\times X\to G\times G$ 
is the $(1,2)$-th projection,
and $r_i$ is the $i$-th projection from $G\times G\times X$.
%

Let $\psi:G\to\BC^\times$ be a character, and let
$\delta\psi\in\Gg^*$ be its differential.
Then $\e^{\lam+\delta\psi}\leftrightarrow \psi\e^\lam$ gives
a canonical isomorphism
\[\CT(\lam+\delta\psi)\cong\CT(\lam)\]
compatible with the multiplicative structure (\ref{mul.sys}).
Hence the twisted equivariance with twist $\lam$
is equivalent to  that with twist $\lam+\delta\psi$.

\newsection{$D$-modules on the flag manifold}
\subsection{Flag manifolds}
\label{subsection:flag manifolds}
We recall basic properties of the flag manifold for the Kac-Moody Lie
algebra $\Gg$ (Kashiwara~\cite{Kflag}).

Fix a $\BZ$-lattice $P$ of $\Gh^*$ satisfying
\eq
\alpha_i\in P,\quad \langle P,h_i\rangle\subset\BZ\quad
\mbox{ for any } i\in I.
\endeq
We define affine group schemes as follows:
\begin{eqnarray}
H&=&\Spec\BC[P],\\
N^{\pm}&=&\limp_k\exp(\Gn^\pm/\Gn^{\pm}_k),\\
B&=&\mbox{(the semi-direct product of } H \mbox{ and } N^+{\rm)},\\
B^-&=&\mbox{(the semi-direct product of } H \mbox{ and } N^-{\rm)}
\end{eqnarray}
(see (\ref{eq:nkpm}) for the definition of $\Gn^{\pm}_k$).
Here, for a finite-dimensional nilpotent Lie algebra $\Ga$ we denote the
corresponding unipotent algebraic group by $\exp(\Ga)$.
Then $N^\pm$ is an affine scheme isomorphic to $\Spec(S(\Gn^{\mp}))$.

For a subset $\Theta$ of $\Delta^\pm$ such that
$(\Theta+\Theta)\cap\Delta\subset\Theta$, we denote by $N(\Theta)$
the subgroup $\exp(\Gn(\Theta))$ of $N^\pm$.

In Kashiwara~\cite{Kflag}, a separated scheme $G$ is constructed
with a free right action of $B^-$ and a free left action of $B$.
The flag manifold $X$ is defined as the
quotient scheme $X=G/B^-$. The flag manifold is a separated scheme.
For $w\in W$, $U_w=wBB^-/B^-$ is an open subset of $X$.
A locally closed
subscheme $X_w=BwB^-/B^-$ of $X$ is called a Schubert cell.

\begin{proposition}\label{D:flag:prop:Bruhat}
\begin{tenumerate}
\item
$X=\bigcup_{w\in W}U_w=\bigsqcup_{w\in W}X_w$.
\item
For any $w\in W$, we have
$\overline{X}_w=\bigcup_{y\ge w}X_y$ and
$X_w\subset U_w\subset\bigcup_{x\le w}X_x$.
\item
We have an isomorphism
\[
N(w\Delta^+\cap\Delta^-)\times N(w\Delta^+\cap\Delta^+)\iso
U_w\,\quad((x,y)\mapsto xywB^-)
\]
of schemes.
Moreover, the subscheme $\{1\}\times
N(w\Delta^+\cap\Delta^+)$ is isomorphic to $X_w$
by this isomorphism.
\end{tenumerate}
\end{proposition}

In particular, $U_w$ and $X_w$ are isomorphic to
\[
\BA^\ell=\Spec\BC[x_k\,;\,0\le k<\ell]
\]
for some $\ell\in\BZ_{\ge0}\sqcup\{\infty\}$, and the codimension of $X_w$
in $X$ is the length $\ell(w)$ of $w\in W$.
Hence $X$ is pro-smooth.

We call a subset $\Phi$ of $W$ {\em admissible} if
\[
w\in\Phi, y\le w\Longrightarrow y\in\Phi.
\]
For an admissible subset $\Phi$ of $W$ we define an open subset $X_\Phi$ of
$X$ by $X_\Phi=\bigcup_{w\in \Phi}X_w$.

For a finite admissible subset $\Phi$, $X_\Phi$ is a quasi-compact scheme
with the condition (S). 
Indeed for $k\gge0$, the subgroup 
$\exp(\Gn^+_k)=\limp_{l\ge k}\exp(\Gn^+_k/\Gn^+_l)$ acts
freely on $X_\Phi$ and $\{X_\Phi/\exp(\Gn^+_k)\}_k$ is
a smooth projective system for $X_\Phi$ (see \cite{KTneg2}).
Note that, since $X_\Phi$ is separated over $\BC$,
$X_\Phi/\exp(\Gn^+_k)$ is separated for $k\gge0$ by the following lemma.

\Lemma
Let $\{X_n\}_{n\in\BN}$ be a projective system of quasi-compact 
and quasi-separated schemes.
Assume that the morphism $X_{n+1}\to X_n$ is an affine morphism
for any $n$.
Let $X_\infty$ be its projective limit.
If $X_\infty$ is separated, then $X_n$ is separated for $n\gge0$.
\enlemma
\proof
Let $f_{nm}:X_m\to X_n$ be the canonical projection ($0\le n\le m\le\infty$).
Since $X_0$ is quasi-compact, 
$X_0$ is covered by finitely many affine open subsets
$U_j^0$ ($j=1,\ldots,N$).
Then the inclusion $U_j^0\to X_0$ is of finite presentation.
Set $U_j^n=f_{0n}^{-1}(U_j^0)$ ($0\le n\le\infty$).
Since $X_\infty$ is separated,
$U_j^\infty\to X_\infty$ is an affine morphism.
Hence \cite[Theorem (8.10.5)]{EGA} 
implies that $U_j^n\to X_n$ is an affine morphism for $n\gge0$.
Hence we may assume from the beginning that
$U_j^0\cap U_k^0$ is affine for any $j,k=1,\ldots,N$.

The ring homomorphism $\CO_X(U_j^\infty)\otimes \CO_X(U_k^\infty)\to
\CO_X(U_j^\infty\cap U_k^\infty)$ is surjective
by the assumption that $X_\infty$ is separated,
Since $\CO_X(U_j^0\cap U_k^0)$ is a finitely generated algebra
over $\CO_X(U_j^0)$, the image of 
$\CO_X(U_j^0\cap U_k^0)\to\CO_X(U_j^n\cap U_k^n)$
is contained in the image of
$\CO_X(U_j^n)\otimes \CO_X(U_k^n)\to
\CO_X(U_j^n\cap U_k^n)$ for $n\gge0$.
On the other hand, we have 
\[
\CO_X(U_j^n\cap U_k^n)\cong
\CO_X(U_j^n)\otimes_{\CO_X(U_j^0)}\CO_X(U_j^0\cap U_k^0).
\]
Hence $\CO_X(U_j^n)\otimes \CO_X(U_k^n)\to
\CO_X(U_j^n\cap U_k^n)$ is surjective for $n\gge0$.
This shows that $X_n$ is separated for $n\gge0$.
\qed

Set $\TX=G/N^-$. 
We denote by $\xi:\TX\to X$ the canonical projection.
It is an $H$-principal bundle.
For $w\in W$ we set $\TX_w=\xi^{-1}X_w=BwN^-/N^-$, and for an admissible
subset $\Phi$ of $W$ we set
$\TX_\Phi=\xi^{-1}X_\Phi=\bigcup_{w\in \Phi}\TX_w$.
The scheme $\TX$ is also pro-smooth.

\subsection{Twisted $D$-modules on the flag manifold}
Let $p:G\to X$ be the projection.
For $\mu\in P$ we define the invertible $\CO_X$-module $\CO_X(\mu)$
as follows:
\eqn
\Gamma(U;\CO_X(\mu))
&=&\{\varphi\in\Gamma(p^{-1}U;\CO_G)\,;\,
\varphi(xg)=g^{-\mu}\varphi(x)
\mbox{ for $(x,g)\in p^{-1}U\times B^-$}\}\\
&\cong&\{\varphi\in\Gamma(\xi^{-1}U;\CO_\TX)\,;\,
\varphi(xh)=h^{-\mu}\varphi(x)
\mbox{ for $(x,h)\in\xi^{-1}U\times H$}\}
\endeqn
for any open subset $U$ of $X$.
Here $x\mapsto x^{-\mu}$ is the character of $B^-$ corresponding to 
the weight $-\mu$.
Twisting $D_X$ by $\CO_X(\mu)$ we obtain a TDO-ring
\[
D_\mu=\CO_X(\mu)\otimes_{\CO_X}D_X\otimes_{\CO_X}\CO_X(-\mu).
\]
This definition can be generalized to any $\mu\in\Gh^*$ 
and we can define an $N^+$-equivariant 
TDO-ring $D_\mu$ on $X$
(see Kashiwara~\cite{Krep} and Kashiwara-Tanisaki~\cite{KTneg2}).
Note that the pull-back $\xi^\sharp D_\mu$ of the TDO-ring $D_\mu$ under
$\xi:\TX\to X$ is canonically isomorphic to $D_{\TX}$.
Hence the pull-back $\xi^\BUL\CM$ of an admissible $D_\mu$-module $\CM$ is
naturally regarded as a $D_{\TX}$-module.
Moreover, by the functor $\xi^\BUL$, the category of
admissible $D_\mu$-modules is equivalent to the category of
admissible twisted $H$-equivariant $D_{\TX}$-modules with twist $\mu$.

\medskip
The infinitesimal action of $\Gg$ on $X$ lifts to an algebra homomorphism
\[
U(\Gg)\to \Gamma(X;D_\mu).
\]
In particular, $H^n(U;\CM)$ has a $\Gg$-module structure for
any open subset $U$ of $X$, any $D_\mu|_U$-module $\CM$ and $n\in\BZ$.

For $w\in W$, let $i_w:X_w\hookrightarrow X$ be the inclusion.
Then for any $\mu\in\Gh^*$, the $N^+$-equivariant
TDO-ring $i_w^\sharp D_\mu$ is 
canonically isomorphic to the $N^+$-equivariant TDO-ring $D_{X_w}$.
We define the $N^+$-equivariant holonomic $D_\mu$-modules
$\CB_w(\mu), \CM_w(\mu)$ by
\eq
\CB_w(\mu)=H^0\int_{i_{w}}\CO_{X_w}, \quad
\CM_w(\mu)=H^0\int_{i_{w}!}\CO_{X_w}.
\endeq
Note that $H^k\int_{i_{w}}\CO_{X_w}=H^k\int_{i_{w}!}\CO_{X_w}=0$ for any
$k\ne0$ because $i_w$ is an affine embedding.
By the definition we have
\eq
\Hom_{D_\mu}(\CM_w(\mu),\CM)&\iso&
\Hom_{D_{\mu}|_U}(\CM_w(\mu)|_U,\CM|_U),\label{eq:proporty of Mw}\\
\Hom_{D_\mu}(\CM,\CB_w(\mu))&\iso&
\Hom_{D_\mu|_U}(\CM|_U,\CB_w(\mu)|_U)\label{eq:proporty of Bw}
\endeq
for any open subset $U$ of $X$ containing $X_w$ and any holonomic
$D_\mu$-module $\CM$.
The isomorphism $\CM_w(\mu)|_U\cong\CB_w(\mu)|_U$
extends to a canonical
non-zero homomorphism $\CM_w(\mu)\to\CB_w(\mu)$.
We denote its image by $\CL_w(\mu)$.
It is a unique irreducible quotient of $\CM_w(\mu)$ and a unique irreducible
submodule of $\CB_w(\mu)$.
Note that $\Hom_{D_\mu}(\CM_w(\mu),\CB_w(\mu))=1$.

For $\mu\in\Gh^*$ and a finite admissible subset $\Phi$ of $W$, we denote by
$\BMD_\Phi(\mu)$ the category of $N^+$-equivariant holonomic
$D_\mu|_{X_\Phi}$-modules.
For $\mu\in\Gh^*$ we denote by $\BMD(\mu)$ the
category of $N^+$-equivariant holonomic $D_\mu$-modules.
Then we have obviously
\[
\BMD(\mu)=\limp_\Phi \BMD_\Phi(\mu),
\]
where $\Phi$ ranges over the set of finite admissible subset of $W$.


For any $w\in W$ and $\mu\in\Gh^*$ the $D_\mu$-modules $\CB_w(\mu)$,
$\CM_w(\mu)$ and $\CL_w(\mu)$ are objects of $\BMD(\mu)$.
%
%
Note that $\CL_w(\mu)$ is a simple object of $\BMD(\mu)$.

For any $N^+$-equivariant admissible $D_\mu$-module $\CM$,
we denote the support of $\CM$ by $\Supp(\CM)$.
It is an $N^+$-stable closed subset of $X$,
and hence it is also $B$-stable.

\begin{lemma}
\label{lem:K-equiv}
For $w\in W$, let $U$ be an open subset of $X$ which contains
$X_w$ as a closed subset.
For any $N^+$-equivariant admissible $D_\mu$-module $\CM$
such that $\Supp(\CM)\cap U\subset X_w$,
there exist some index set $J$ 
and isomorphisms
\[
\CM_w(\mu)^{\oplus J}|_U\iso\CM|_U\iso\CB_w(\mu)^{\oplus J}|_U\,.
\]
Furthermore if $\CM$ is holonomic, then
$J$ is a finite set and the above isomorphisms
can be extended to morphisms in $\BMD(\mu)$
\[
\CM_w(\mu)^{\oplus J}\to\CM\to\CB_w(\mu)^{\oplus J}.
\]
\end{lemma}
\proof
Let $i:X_w\to U$ be the closed embedding.
By the condition on $\CM$ there exists an $N^+$-equivariant
holonomic $D_{X_w}$-module
$\CN$ such that $\CM|_U\simeq\int_i\CN$.
Since $X_w$ is a homogeneous space of $N^+$ 
with a connected isotropy subgroup,
we see that $\CN$ is isomorphic to $\CO_{X_w}^{\oplus J}$ for some $J$.
If $\CM$ is holonomic then $\CN$ is holonomic, and hence $J$ is a finite set.
Thus we have
\[
\CM|_U\simeq \int_i\CO_{X_w}^{\oplus J}\simeq(\CB_w(\mu)|_U)^{\oplus
J}\simeq(\CM_w(\mu)|_U)^{\oplus J}.
\]
To see the last statement, it is enough to apply (\ref{eq:proporty of Mw})
and (\ref{eq:proporty of Bw}).
\qed

Let $\CM\in\Ob(\BMD(\mu))$.
By Lemma \ref{lem:K-equiv}, for any finite admissible subset $\Phi$ of $W$, 
$\CM|_{X_\Phi}\in \Ob(\BMD_\Phi(\mu))$ has finite length
and it has finite
composition series whose composition factors are isomorphic to
$\CL_w(\mu)|_{X_\Phi}$ for some $w\in\Phi$.
For $w\in W$ the multiplicity of $\CL_w(\mu)|_{X_\Phi}$ in the composition
series of $\CM|_{X_\Phi}$ does not depend on the choice of a finite admissible
subset $\Phi$ of $W$ such that $w\in \Phi$.
We denote it by $[\CM:\CL_w(\mu)]$.
Note that the multiplicity does not depend on the $N^+$-equivariance
structure.

\begin{lemma} \label{D:cat:lemma:mult}
We have $[\CM_w(\mu):\CL_y(\mu)]=[\CB_w(\mu):\CL_y(\mu)]$ for any $w,y \in W$.
\end{lemma}
%
Replacing the modules $\CM_w(\mu)$, $\CL_y(\mu)$ and $\CB_w(\mu)$
with their images by $\xi^\BUL$, 
this follows from the following general result.

\Prop[\cite{KKw}]
Let $j:X\to Y$ be an embedding of smooth $\BC$-schemes.
Then for any holonomic $D_X$-module $\CM$,
we have the equality
\[
\sum_{i\in\BZ}(-1)^i[H^i(\int_j\CM)]
=\sum_{i\in\BZ}(-1)^i[H^i(\int_{j!}\CM)]
\]
in the Grothendieck group of the category of
holonomic $D_Y$-modules.
\enprop
\proof
We can decompose this proposition into the closed 
embedding case and the open embedding case.
Since the first case is obvious, we may assume that $j$ is an open embedding.
Since the question being local on $Y$,
we can easily reduce to the case where $X$ is 
the complement of a hypersurface of $Y$.
Then $\int_j\CM$ and $\int_{j!}\CM$ are concentrated at degree $0$.
We may assume further
that $Y\setminus X$ is defined by $f=0$ for some $f\in\Gamma(Y;\CO_Y)$.
Let $\psi_f(\CM)$ be the near-by cycle of $\CM$,
which is a holonomic $D_Y$-module with support in $Y\setminus X$.
Let $\var:\psi_f(\CM)\to\psi_f(\CM)$ be the variation.
Then the kernel (resp. the cokernel) of
$\int_{j!}\CM\to\int_{j}\CM$ is isomorphic to 
the kernel (resp. the cokernel) of $\var:\psi_f(\CM)\to\psi_f(\CM)$.
Therefore we have
\[
[\int_{j}\CM]-[\int_{j!}\CM]
=[\Coker(\var)]-[\Ker(\var)]=0.\]
\qed

\subsection{Cohomologies of $\CB_w(\mu)$}

We first study the cohomology groups of $\CB_w(\mu)$.
\begin{proposition}
\label{D:mod:prop:gamma(B_w)}
Let $\mu\in\Gh^*$, $w\in W$ and let $\Phi$ be a finite admissible subset of $W$
containing $w$.
Then we have
\begin{tenumerate}
\item
$H^n(X_\Phi;\CB_w(\mu))=H^n(X;\CB_w(\mu))=0$ for any $n\ne0$.
\item
We have
$$\Gamma(X_\Phi;\CB_w(\mu))=
\Gamma(X;\CB_w(\mu))\simeq
U({}^w\Gb)\otimes_{U({}^w\Gb\cap\Gb)}(\Gamma(X_w;\CO_{X_w})
\otimes\BC_{w\circ\mu})$$
as a $U({}^w\Gb)$-module. 
Here ${}^w\Gb=\Gh\oplus(\bigoplus_{\alpha\in w\Delta^+}\Gg_\alpha)$ and
$\BC_{w\circ\mu}$ is the one-dimensional $U({}^w\Gb\cap\Gb)$-module
with $w\circ\mu$ as a weight.
\end{tenumerate}
\end{proposition}
\proof
By the definition of $\CB_w(\mu)$ we have
\[
H^n(X_\Phi;\CB_w(\mu))
\simeq
H^n(X;\CB_w(\mu))
\simeq
H^n(U_w;\CB_w(\mu)).
\]
Hence the assertion easily follows from
Proposition~\ref{D:flag:prop:Bruhat} (iii) (cf. \cite{Kpos1,KTneg2}).
\qed

Later, we shall see that $\Gamma(X;\CB_w(\mu))$ is isomorphic
to the dual Verma module under certain conditions.
The following corollary is a key of its proof.
\begin{corollary}\label{D:mod:cor:gamma(B_w)1}
Let $\mu\in\Gh^*$ and $w\in W$.
\begin{tenumerate}
\item
$\ch(\Gamma(X;\CB_w(\mu)))=\ch(M(w\circ\mu)).$
\item
If $\zeta\in\Gh^*$ and non-zero $m\in
\Gamma(X;\CB_w(\mu))_\zeta$ satisfy
$\Gn(w\Delta^+\cap\Delta^+)m=0$, 
then we have 
$\zeta\in w\circ\mu-
\sum_{\alpha\in\Delta^+\cap w\Delta^-}\BZ_{\geq0}\alpha$.
\item
If $\zeta\in\Gh^*$ and non-zero $m\in
\Gamma(X;\CB_w(\mu))^*_\zeta$ satisfy
$\Gn(w\Delta^-\cap\Delta^+)m=0$, then we have $\zeta\in
w\circ\mu-\sum_{\alpha\in\Delta^+\cap w\Delta^+}\BZ_{\geq0}\alpha$.
\end{tenumerate}
\end{corollary}
\proof
By Proposition~\ref{D:mod:prop:gamma(B_w)} we have
\eq\label{eq:26}
\begin{array}{rcl}
\Gamma(X;\CB_w(\mu))&\simeq&
U(\Gn(w\Delta^+\cap\Delta^-))\otimes\Gamma(X_w;\CO_{X_w})\otimes
\BC_{w\circ\mu}\\
&\simeq&
U(\Gn(w\Delta^+\cap\Delta^-))\otimes
S(\Gn(w\Delta^-\cap\Delta^-))\otimes\BC_{w\circ\mu},
\end{array}
\endeq
and hence we have
\eq
\ch(\Gamma(X;\CB_w(\mu)))&=&
\ch(U(\Gn(w\Delta^+\cap\Delta^-)))
\ch(S(\Gn(w\Delta^-\cap\Delta^-)))
\e^{w\circ\mu}\nn\\
&=&\prod_{\alpha\in w\Delta^+\cap\Delta^-}
(1-\e^{\alpha})^{-\dim\Gg_{\alpha}}
\prod_{\alpha\in w\Delta^-\cap\Delta^-}
(1-\e^{\alpha})^{-\dim\Gg_{\alpha}}\,
\e^{w\circ\mu}\nn\\
&=&\prod_{\alpha\in\Delta^-}
(1-\e^{\alpha})^{-\dim\Gg_{\alpha}}\,
\e^{w\circ\mu}\nn\\
&=&\ch(M(w\circ\mu)).\nn
\endeq
Thus (i) is proved.

Let us prove (iii).
Assume that a non-zero vector $m\in
\Gamma(X;\CB_w(\mu))^*_\zeta$ with $\zeta\in\Gh^*$ satisfies
$\Gn(w\Delta^-\cap\Delta^+)m=0$. Then we have
$\Bigl(\Gamma(X;\CB_w(\mu))/\Gn(w\Delta^+\cap\Delta^-)
\Gamma(X;\CB_w(\mu))\Bigr)_\zeta\ne0$.
Then (\ref{eq:26}) implies
\eqn
&&\ch\Bigl(\Gamma(X;\CB_w(\mu))/\Gn(w\Delta^+\cap\Delta^-)
\Gamma(X;\CB_w(\mu))\Bigr)\\
&=&
\ch(S(\Gn(w\Delta^-\cap\Delta^-)))
\e^{w\circ\mu}\\
&=&\prod_{\alpha\in w\Delta^-\cap\Delta^-}
(1-\e^{\alpha})^{-\dim\Gg_{\alpha}}\,
\e^{w\circ\mu},
\endeqn
and we obtain (iii).

Let us finally show (ii).
Define a filtration $\{F_\ell\}_{\ell\in\BZ}$ of
$U(\Gn(w\Delta^+\cap\Delta^-))$ by
\[
F_\ell=
\bigoplus_{\HT(\beta)\le\ell}
U(\Gn(w\Delta^+\cap\Delta^-))_{-\beta}.
\]
Then the subspace $F_\ell\otimes\Gamma(X_w;\CO_{X_w})\otimes\BC_{w\circ\mu}$
of $\Gamma(X;\CB_w(\mu))$  is stable under the action of
$\Gn(w\Delta^+\cap\Delta^+)$, and the action of
$\Gn(w\Delta^+\cap\Delta^+)$ on the quotient
\eq
&&(F_\ell\otimes\Gamma(X_w;\CO_{X_w})\otimes\BC_{w\circ\mu})/
(F_{\ell-1}\otimes\Gamma(X_w;\CO_{X_w})\otimes\BC_{w\circ\mu})\nn\\
&=&(F_\ell/F_{\ell-1})\otimes\Gamma(X_w;\CO_{X_w})\otimes\BC_{w\circ\mu}\nn
\endeq
is given by
\[
x(u\otimes n)=u\otimes xn\quad
\mbox{for $x\in\Gn(w\Delta^+\cap\Delta^+)$,
$u\in F_\ell/F_{\ell-1}$,
$n\in\Gamma(X_w;\CO_{X_w})\otimes\BC_{w\circ\mu}$.}
\]
Take the smallest $\ell$ such that $m\in
F_\ell\otimes\Gamma(X_w;\CO_{X_w})\otimes\BC_{w\circ\mu}$ and denote by
$\overline{m}$ the corresponding element of
$(F_\ell/F_{\ell-1})\otimes\Gamma(X_w;\CO_{X_w})\otimes\BC_{w\circ\mu}$.
Write $\overline{m}$ as
\[
\overline{m}=
\sum_{j=1}^r\overline{u}_j\otimes n_j,
\]
where $\overline{u}_j\,(j=1,\dots, r)$ are linearly independent elements of
$F_\ell/F_{\ell-1}$ and $n_j\,(j=1,\dots, r)$ are  elements of
$\Gamma(X_w;\CO_{X_w})\otimes\BC_{w\circ\mu}$.
By the assumption on $m$ we have $\Gn(w\Delta^+\cap\Delta^+)n_j=0$ for any
$j$.
Since $X_w$ is a homogeneous space of $N(w\Delta^+\cap\Delta^+)$, we have
$n_j\in\BC\otimes\BC_{w\circ\mu}
\subset\Gamma(X_w;\CO_{X_w})\otimes\BC_{w\circ\mu}$.
Thus
\[
\overline{m}\in(F_\ell/F_{\ell-1})\otimes\BC\otimes\BC_{w\circ\mu},
\]
and we obtain (ii).
\qed

\begin{proposition}
\label{D:mod:cor:gamma(B_w)2}
For $\lambda\in\KRP$ and $w\in W$, we have
\[
\Gamma(X;\CB_w(\lambda))\simeq
M^*(w\circ\lambda).
\]
\end{proposition}
\proof
Since $\ch(\Gamma(X;\CB_w(\lambda)))=
\ch(M^*(w\circ\lambda))$, it is sufficient to show that if there exists
$m\in\Gamma(X;\CB_w(\lambda))_\zeta\setminus\{0\}$ such that
$\Gn^+m=0$, then $\zeta=w\circ\lam$.

Since $[M(w\circ \lam):L(\zeta)]\ne0$,
Corollary \ref{D:mod:cor:[M:L]} implies $\zeta=y\circ\lam$ 
for $y\in wW(\lam)$.
Hence there exist some $\gamma_1, \cdots, \gamma_r\in\Delta^+(\lam)$ such that
\[
w^{-1}y=s_{\gamma_1}\cdots s_{\gamma_r},\qquad
s_{\gamma_1}\cdots s_{\gamma_{j-1}}\gamma_j\in\Delta^+.
\]
Then we have
\eqn
\lam+\rho-w^{-1}y(\lam+\rho)
&=&\lam+\rho-s_{\gamma_1}\cdots s_{\gamma_r}(\lam+\rho)\\
&=&\sum_{j=1}^r(s_{\gamma_1}\cdots
s_{\gamma_{j-1}}(\lam+\rho)-s_{\gamma_1}\cdots s_{\gamma_{j}}(\lam+\rho))
\\
&=&\sum_{j=1}^r(\lam+\rho,\gamma_j^\vee) s_{\gamma_1}\cdots
s_{\gamma_{j-1}}(\gamma_j).
\endeqn
Since $\lam\in\KRP$, we have $(\lam+\rho,\gamma_j^\vee)\in\BZ_{>0}$, and
hence $\lam+\rho-w^{-1}y(\lam+\rho)\in Q^+$.

On the other hand,
Corollary~\ref{D:mod:cor:gamma(B_w)1} (ii) implies
\[
\lam+\rho-w^{-1}y(\lam+\rho)
=w^{-1}(w\circ\lam-\zeta)
\in w^{-1}(\sum_{\alpha\in\Delta^+\cap w\Delta^-}\BZ_{\ge0}\alpha)
\subset-Q^+\,.
\]
Thus we obtain $\lam+\rho-w^{-1}y(\lam+\rho)=0$, 
and hence $\zeta=w\circ\lam$.
\qed
\Remark
This proposition \ref{D:mod:cor:gamma(B_w)2} can be also proved
by using the theory of the Radon transform in \S \ref{Radon}.
Indeed, Theorem \ref{thm:Radon:2} implies 
$\Gamma(X;\CB_w(\lam))\cong\Gamma(X;\CB_e(w\circ\lam))$.
The last module is isomorphic to
$\Gamma(X_e;\CO_X)\otimes \BC_{w\circ\lam}$,
and it has a unique highest weight vector,
because $X_e$ is a homogeneous space of $N^+$.

\subsection{Modified cohomology groups}
The cohomology group $H^n(X;\CM)$ itself may be too wild.
We shall replace it with a modified one easier to manipulate.

\begin{lemma}
\label{D:mod:lem:Hn(XPhi;M)}
Let $\Phi$ be a finite admissible subset of $W$.
For any $\mu\in\Gh^*$, $n\in\BZ$ and $\CM\in\Ob(\BMD(\mu))$,
we have the following.
\begin{tenumerate}
\item
$H^n(X_\Phi;\CM)$ is an object of $\BO$.
\item
If $[H^n(X_\Phi;\CM):L(\zeta)]\ne0$, then there exists some $w\in \Phi$ such
that $X_w\subset\Supp(\CM)$ and that $[M(w\circ\mu):L(\zeta)]\ne0$.
\item
For any admissible subset $\Psi$ of $W$ such that $\Psi\subset\Phi$,
let $N_1$ $($resp. $N_2$$)$ be the kernel $($resp. cokernel$)$
of the natural homomorphism
$H^n(X_\Phi;\CM)\to H^n(X_\Psi;\CM)$.
Then $N_i$ belongs to $\Ob(\BO)$ for $i=1,2$.
Moreover, if $[N_i:L(\zeta)]\ne0$ for $i=1$ or $2$, 
then $[M(x\circ\mu):L(\zeta)]\ne0$ for some $x\in\Phi\setminus\Psi$.
\end{tenumerate}
\end{lemma}
\proof
%
We first show (iii) 
by the induction of $\sharp(\Phi\setminus\Psi)$.

If $\sharp(\Phi\setminus\Psi)=0$, it is trivial.
In the case $\sharp(\Phi\setminus\Psi)=1$,
set $\Phi\setminus\Psi=\{x\}$.
Let $i:X_x\to X_{\Phi}$ and $j:X_{\Psi}\to X_{\Phi}$ be the inclusion.
By the assumption $i$ is a closed embedding and $j$ is an open embedding.
The distinguished triangle
\[
\BD i_*\BD i^!\CM\to \CM\to \BD j_*\BD j^\BUL M\maprightu{+1}
\]
induces an exact sequence
\eq
&&H^n(X_{\Phi};\BD i_*\BD i^!\CM)\to
H^n(X_{\Phi};\CM)\to
H^n(X_{\Psi};\CM)\to
H^{n+1}(X_{\Phi};\BD i_*\BD i^!\CM).
\endeq
%
Therefore the kernel $N_1$ is a quotient of 
$H^n(X_{\Phi};\BD i_*\BD i^!\CM)$,
and the cokernel $N_2$ is a submodule of 
$H^{n+1}(X_{\Phi};$$\BD i_*\BD i^!\CM)$.
By Lemma~\ref{lem:K-equiv}, the object $H^k(\BD i_*\BD i^!\CM)$ 
in $\BMD_\Phi(\mu)$ is isomorphic to a
direct sum of finitely many copies of $\CB_x(\mu)|_{X_\Phi}$.
Hence Proposition \ref{D:mod:prop:gamma(B_w)} (i)
implies
\[
H^k(X_{\Phi};\BD i_*\BD i^!\CM)=
\Gamma(X_{\Phi};H^k(\BD i_*\BD i^!\CM))\,,\]
and its character is a constant multiple of the character of
the Verma module $M(x\circ\mu)$ for any $k$
by Corollary~\ref{D:mod:cor:gamma(B_w)1}.
This shows (iii) in the case $\sharp(\Phi\setminus\Psi)=1$.

Assume $\sharp(\Phi\setminus\Psi)>1$.
Taking a maximal element $x$ of $\Phi\setminus\Psi$,
set $\Phi'=\Phi\setminus\{x\}$.
Then $\Phi'$ is an admissible subset such that
$\Psi\subset\Phi'\subset\Phi$.
Consider the diagram
\eqn&&
H^n(X_\Phi;\CM)\maprightu{\alpha}H^n(X_{\Phi'};\CM)
\maprightu{\beta}H^n(X_\Psi;\CM).
\endeqn
Then we have an exact sequence
\eqn
0\to \Ker \alpha\to\Ker (\beta\circ\alpha)\to\Ker \beta\\
\Cok \alpha\to\Cok (\beta\circ\alpha)\to\Cok \beta\to0.
\endeqn
By the induction hypothesis,
$\Ker \alpha$ and $\Ker \beta$ belong to $\BO$.
Hence $N_1=\Ker (\beta\circ\alpha)$ belongs to $\BO$.
If $[N_1:L(\zeta)]\not=0$, then
$[\Ker \alpha:L(\zeta)]\not=0$ or $[\Ker \beta:L(\zeta)]\not=0$.
The induction hypothesis implies
$[M(x\circ\mu):L(\zeta)]\not=0$ in the first case
and $[M(w\circ\mu):L(\zeta)]\not=0$ for some $w\in\Phi'\setminus\Psi$
in the second case. This shows the assertion for $N_1$.
The assertion for $N_2$ is similarly proved.

\medskip
We  obtain (i) and (ii) from (iii)
by taking $\Psi=\emptyset$ or
$\Psi=\{w\in\Phi\,;\,X_w\cap\Supp \CM=\emptyset\}$.
\qed

By Lemma~\ref{D:mod:lem:Hn(XPhi;M)}, Proposition~\ref{D:mod:prop:[M:L]}
and the $W$-invariance of $\KR$, we have the following corollary.

\begin{corollary}\label{cor:coh}
For $\lam\in\KR$ and $\CM\in\Ob(\BMD(\lam))$,
we have
\[
H^n(X_\Phi;\CM)\in\Ob(\BO(\lam))
\]
for any finite admissible subset $\Phi$ of $W$ and any $n\in\BZ$.
\end{corollary}

\begin{lemma}
\label{D:mod:eq:bij}
Let $\lambda, \mu\in\KR$.
Then for any $\zeta\in\Gh^*$ there exists a finite admissible subset $\Phi$
of $W$ such that the restriction homomorphism
\[(P_\mu(H^n(X_{\Phi'};\CM)))_{\zeta'}\to
(P_\mu(H^n(X_{\Phi};\CM)))_{\zeta'}\]
is bijective for any finite admissible
subset $\Phi'$ of $W$ containing $\Phi$, $\zeta'\in\zeta+Q^+$,
$\CM\in\Ob(\BMD(\lambda))$ and $n\geq0$ 
$($see $(\ref{eq:projection functor})$ for the definition of $P_\mu$$)$.
\end{lemma}
\proof
We may assume $\mu\in\KRP$.
Then $W(\mu)\circ\mu\cap(\zeta-\mu+Q^+)$ is a finite set.
Since $\{w\in W;w\circ\lam=\lam\}=\{1\}$,
there exist only finitely many $w\in W$ satisfying $w\circ\lam\in
W(\mu)\circ\mu$ and $w\circ\lam-\zeta\in Q^+$.
Thus we conclude that there exists a finite
admissible subset $\Phi$ of $W$ satisfying
\[
w\in W,\; w\circ\lam\in W(\mu)\circ\mu, \;w\circ\lam-\zeta\in Q^+
\Longrightarrow
w\in\Phi.
\]
Then the assertion follows from Lemma~\ref{D:mod:lem:Hn(XPhi;M)} (iii)
and the assumption on $\Phi$.
\qed

For $\lam,\mu\in\KR$, $\CM\in\BMD(\lam)$ and $n\in\BZ_{\ge0}$ we set
\eq
\tH_\mu^n(\CM)=\bigoplus_{\xi\in\Gh^*}\limp_\Phi
\left(P_\mu(H^n(X_\Phi;\CM))\right)_\xi,
\endeq
where $\Phi$ is running over finite admissible subsets of $W$.

\Lemma
For $\lam$, $\mu\in\KR$,
$\tH_\mu^n$ is an additive functor
from $\BMD(\lam)$ to $\BO[\mu]$ for any integer $n$.
\enlemma
\proof
Let $\CM\in\Ob(\BMD(\lam))$.
Take any $\zeta\in\Gh^*$.
Let $\Phi$ be as in Lemma~\ref{D:mod:eq:bij}.
Then we have
\[
\sum_{\xi\in\zeta+Q^+}\dim\tH_\mu^n(\CM)_\xi
=\sum_{\xi\in\zeta+Q^+}\dim (P_\mu(H^n(X_\Phi;\CM)))_\xi
<\infty
\]
by Lemma~\ref{D:mod:lem:Hn(XPhi;M)} (i).
Thus $\tH^n_\mu(\CM)\in\Ob(\BO)$.
\qed

Now for $\lam\in\KR$ and $\CM\in\BMD(\lam)$ we define the object of $\TBO$ 
(see (\ref{def:TBO})) as follows:
\eq
&&\tH^n(\CM)=\prod_{\mu\in\KRP}\tH^n_\mu(\CM).
\endeq
Then $\tH^n(\CM)$ is the projective limit of
$H^n(X_\Phi;\CM)$ in the category $\TBO$.
\hb
We write $\tG(\CM)$ for $\tH^0(\CM)$.

\begin{proposition}\label{D:mod:prop:tH1}
Let $\lambda\in\KR$.
\begin{tenumerate}
\item
$\tH^n$ is a functor from $\BMD(\lam)$ to $\TBO(\lam)$.
\item
For any short exact sequence
\[
0\to \CM_1\to \CM_2\to \CM_3\to0
\]
in $\BMD(\lambda)$, we have a functorial long exact sequence
\[
\begin{array}{ccccccc}
0&\to&\tH^0(\CM_1)&\to&\tH^0(\CM_2)&\to&\tH^0(\CM_3)\\
&\to&\tH^1(\CM_1)&\to&\tH^1(\CM_2)&\to&\tH^1(\CM_3)\\
&\to&\cdots
\end{array}
\]
in $\TBO$.
\item
Assume that $[\tH^n(\CM):L(\zeta)]\ne0$ for $\CM\in\Ob(\BMD(\lambda))$,
$\zeta\in\Gh^*$, $n\in\BZ_{\ge0}$.
Then there exists some $w\in \Phi$ such
that $X_w\subset\Supp(\CM)$ and that $[M(w\circ\mu):L(\zeta)]\ne0$.
\item
If $\lam\in\KRP$ and $\CM\in\BMD(\lam)$, then
$\tH^n(\CM)$ belongs to $\TBO\{S(\CM)\circ\lam\}$,
where $S(\CM)=\{w\in W\,;\,X_w\subset\Supp(\CM)\}$.
\end{tenumerate}
\end{proposition}
\proof
(i), (ii) and (iii) easily follow from the definition along with
Lemma~\ref{D:mod:lem:Hn(XPhi;M)} and Corollary \ref{cor:coh}.
Let us prove (iv).
Assume $[\tH^n(\CM):L(\zeta)]\ne0$.
By (iii) there exists some $x\in \Phi$ such
that $X_x\subset\Supp(\CM)$ and that $[M(x\circ\mu):L(\zeta)]\ne0$.
Then by Corollary \ref{D:mod:cor:[M:L]},
we have $\zeta=w\circ\lam$ for
some $w\in W$ with $w\geq x$.
Hence $X_w\subset\ol X_x\subset \Supp(\CM)$.
Therefore $w\in S(\CM)$.
%
\qed
%
\begin{lemma}
\label{lem:Gamma}
Let $\lam\in\KR$.
For $M\in\Ob(\BO\{\KR\})$ and $\CM\in\Ob(\BMD(\lam))$ we have
\[
\Hom_\TBO(M,\tG(\CM))\simeq
\Hom_\Gg(M,\Gamma(X;\CM)).
\]
\end{lemma}
\proof
Note that $\tG(\CM)$ is the projective limit of
$\{\Gamma(X_\Phi;\CM)\}_\Phi$ in $\TBO$,
where $\Phi$ ranges over the set of finite admissible subsets of $W$,
while $\Gamma(X;\CM)$ is the projective limit of
$\{\Gamma(X_\Phi;\CM)\}_\Phi$ in the category of $\Gg$-modules.
Hence we have
\eqn
\Hom_\TBO(M,\tG(\CM))&\simeq&
\limp_\Phi\Hom_\TBO(M,\Gamma(X_\Phi;\CM))\\
&\simeq&
\limp_\Phi\Hom_\Gg(M,\Gamma(X_\Phi;\CM))\\
&\simeq&
\Hom_\Gg(M,\Gamma(X;\CM)).
\endeqn
\qed

By  Corollary~\ref{D:mod:cor:gamma(B_w)1} and
Proposition~\ref{D:mod:cor:gamma(B_w)2} we have the following proposition.
\begin{proposition}\label{D:mod:prop:tH2}
Let $\lambda\in\KR$.
\begin{tenumerate}
\item
$\tH^n(\CB_w(\lambda))=0$ for $n\ne0$.
\item
$\tG(\CB_w(\lambda))=\Gamma(X;\CB_w(\lambda))$ and
$\ch(\tG(\CB_w(\lambda)))=\ch(M(w\circ\lam))$.
\item
$\tG(\CB_w(\lambda))\simeq
M^*(w\circ\lambda)$ if $\lam\in\KRP$.
\end{tenumerate}
\end{proposition}
\begin{proposition}
\label{cohom:prop:Hom}
For $\lambda\in\KRP$ and $w\in W$ we have
\eq
&&
\Hom_\TBO(M(w\circ\lambda),\tG(\CM_w(\lambda)))
\iso
\Hom_\TBO(M(w\circ\lambda),\tG(\CB_w(\lambda)))
\cong\BC.
\label{cohom:prop:Hom:eq2}
\endeq
\end{proposition}
\proof
Set $\CN_1=\CB_w(\lam)/\CL_w(\lam), \CN_2=\Ker(\CM_w(\lam)\to\CL_w(\lam))$.
Then we have $\Supp(\CN_i)\subset\overline{X}_w\setminus X_w$ for $i=1,2$.
Hence Proposition~\ref{D:mod:prop:tH1} (iv) implies
$[\tH^n(\CN_i):L(w\circ\lam)]=0$.
Taking the cohomologies of the short exact sequences
\eqn
&&0\to \CN_2\to\CM_w(\lam)\to\CL_w(\lam)\to0,\\[5pt]
&&0\to \CL_w(\lam)\to\CB_w(\lam)\to\CN_1\to0,
\endeqn
we obtain exact sequences
\eqn
&&0\to
\tG(\CN_2)\to\tG(\CM_w(\lam))\to\tG(\CL_w(\lam))\to\tH^1(\CN_2),\\[5pt]
&&0\to
\tG(\CL_w(\lam))\to\tG(\CB_w(\lam))\to\tG(\CN_1).
\endeqn
Therefor we have
\eqn
[\tG(\CM_w(\lambda)):L(w\circ\lam)]&=&
[\tG(\CL_w(\lambda)):L(w\circ\lam)]\\
&=&[\tG(\CB_w(\lambda)):L(w\circ\lam)]\\
&=&[M^*(w\circ\lam):L(w\circ\lam)]\\
&=&1\,.
\endeqn
Here the third equality follows from (iii) in
the preceding proposition.
%
By Proposition~\ref{D:mod:prop:tH1} (iv),
we have
\[[\tG(\CM_w(\lam)):L(y\circ\lam)]=0
\quad\mbox{for $y<w$.}\]
Hence Lemma~\ref{KM:cat:dimHom2} implies
that $\Hom_\TBO(M(w\circ\lambda),\tG(\CM_w(\lambda)))$ does not vanish.
By the exact sequence
\[0\to\tG(\CN_2)\to\tG(\CM_w(\lam))\to\tG(\CB_w(\lam))\]
and $\Hom_\TBO(M(w\circ\lam),\tG(\CN_2))=0$,
the homomorphism
\[\Hom_\TBO(M(w\circ\lam),\tG(\CM_w(\lam)))
\to
\Hom_\TBO(M(w\circ\lam),\tG(\CB_w(\lam)))\cong\BC\]
is injective, which implies the desired result.
\qed

\subsection{Modified localization functor}

For $\mu\in\Gh^*$ there exists a unique additive functor
\eq
\label{eq:localization}
D_\mu\htensor\ \bullet:
\BMA
\to
\MA(D_\mu)
\qquad
(M\mapsto D_\mu\htensor M),
\endeq
called the modified localization functor, such that
\eq
\label{eq:localization2}
&&\Hom_\Gg(M,\Gamma(X;\CM))=
\Hom_{D_\mu}(D_\mu\htensor M,\CM)\\
&&
\qquad\qquad
\mbox{ for }
M\in\Ob(\BMA),\;
\CM\in\Ob(\MA(D_\mu)).
\nn
\endeq
In \cite{Kpos1} it is constructed in the case where $\mu$ is integral.
Since the construction in the general case is completely similar, we do not
repeat it here.
As in \cite{Kpos1} we have the following proposition.

\begin{proposition}
Let $\mu\in\Gh^*$.
\begin{tenumerate}
\item
The functor {\rm (\ref{eq:localization})} is right exact, and commutes with
the inductive limit.
\item
For any $M\in\BO$, 
$D_\mu\htensor M$ is an $N^+$-equivariant admissible $D_\mu$-module.
\end{tenumerate}
\end{proposition}

In particular, for $M\in\Ob(\BO)$, the support $\Supp(D_\mu\htensor M)$ of
$D_\mu\htensor M$ is a $B$-stable closed subset of $X$.

By Lemma~\ref{lem:Gamma} we have the following lemma.
\begin{lemma}
\label{lem:htensor:hom}
For $\lam\in\KR$,
$M\in\Ob(\BO\{\KR\})$ and $\CM\in\Ob(\BMD(\lam))$ we have
\[
\Hom_\TBO(M,\tG(\CM))=
\Hom_{D_\lam}(D_\lam\htensor M,\CM)\,.
\]
\end{lemma}

\begin{proposition}
\label{D:loc:prop:supp}
Let $\lam\in\KRP$, $M\in\Ob(\BO)$.
\begin{tenumerate}
\item
Assume that $X_w$ is open in $\ol X_w\bigcup\Supp(D_{\lam}\htensor M)$.
Then $D_\lam\htensor M|_{U_w}$
is isomorphic to the direct sum of 
$\dim \Hom_\Gg(M,M^*(w\circ\lam))$ copies of
$\CB_w(\lam)|_{U_w}$.
\item
$\Supp(D_{\lam}\htensor M)$ is the union of
$\overline{X}_w$ such that $\Hom_\Gg(M,M^*(w\circ\lam))\ne0$.
\item
If $\lam\in\KRP$, then
$\Supp(D_{\lam}\htensor M)$ is the union of
$\ol{X}_w$ such that $[M:L(w\circ\lam)]\ne0$.
\end{tenumerate}
\end{proposition}
\proof
Let us first show (i).
Assume that $X_w$ is open in $\ol X_w\cup\Supp(D_{\lam}\htensor M)$.
Lemma~\ref{lem:K-equiv} implies that
$D_{\lam}\htensor M|_{U_w}$ is isomorphic to
$\CB_w(\lam)^{\oplus J}|_{U_w}$ for some index set $J$.
Hence by (\ref{eq:proporty of Bw}), (\ref{eq:localization2}), and
Proposition~\ref{D:mod:cor:gamma(B_w)2} we have
\eqn
\Hom_\Gg(M,M^*(w\circ\lam))
&\simeq&
\Hom_{D_\lam}(D_\lam\htensor M,\CB_w(\lam))\\
&\simeq&
\Hom_{D_\lam|_{U_w}}(D_\lam\htensor M|_{U_w},\CB_w(\lam)|_{U_w})
\\
&\simeq&\Hom_\BC(\BC^{\oplus J},\BC).
\endeqn
Hence $\sharp J=\dim\Hom_\Gg(M,M^*(w\circ\lam))<\infty$.
Let us show (ii). By (i) it is obvious that
$\Supp(D_\lam\htensor M)$ in contained in the union.
Conversely assume $\Hom_\Gg(M,M^*(w\circ\lam))\ne0$.
If $X_w$ is not contained in $\Supp(D_\lam\htensor M)$,
then $X_w$ is open in $X_w\cup\Supp(D_\lam\htensor M)$,
and hence (i) implies that $\Supp(D_\lam\htensor M)$ contains $X_w$.

Let us show (iii).
By (ii), it is enough to show that $[M:L(w\circ\lam)]\ne0$ implies
$X_w\subset\Supp(D_\lam\htensor M)$.
Let us take $x\in W$ such that $x\le w$, $[M:L(x\circ\lam)]\ne0$
and $[M:L(y\circ\lam)]=0$ for any $y<x$.
Then Lemma \ref{KM:cat:dimHom2} implies $\Hom_\Gg(M,M^*(x\circ\lam))\ne0$.
Hence (ii) implies $\Supp(D_\lam\htensor M)\supset \ol{X}_x\supset X_w$.
\qed

\Prop
For $\lam\in\KRP$ the functor {\rm (\ref{eq:localization})} 
induces the functor
\eq
D_\mu\htensor\ \bullet:
\BO
\to
\BMD(\lam).
\endeq
\enprop
\proof
Set $M_0=U(\Gg)/U(\Gg)\Gn^+$.
We shall first show that
$D_\lam\htensor M_0$ is holonomic.
Let $\Phi$ be a finite admissible subset of $W$.
Set $Y_k=X_{\Phi}/\exp(\Gn_k^+)$, and let
$p_k:X_{\Phi}\to Y_k$ be the projection.
Then $Y_k$ is a smooth $\BC$-scheme for $k\gge0$.
For $k_0\gge0$, let $\{(Y_k,A_k)\}_{k\ge k_0}$ be 
a smooth projective system of $(X_\Phi,D_\lam|_{X_\Phi})$.
Since $Y_k$ has finitely many orbits by the action of
$\exp(\Gn^+/\Gn^+_k)$,
the $A_k$-module $A_k\otimes_{U(\Gn^+/\Gn^+_k)}\BC$ 
is a holonomic $A_k$-module.
Since $D_\lam\htensor M_0|_{X_\Phi}=
p_k^\BUL (A_k\otimes_{U(\Gn^+/\Gn^+_k)}\BC)$
for $k\gge0$ (see \cite{Kpos1}),
it is also holonomic.
Thus we have proved that
$D_\lam\htensor M_0$ is holonomic.
Then we see that
$D_\lam\htensor(U(\Gg)\otimes_{U(\Gn^+)}V)$ is also holonomic
for any finite-dimensional $\Gb$-module $V$.
Indeed $V$ has a finite filtration whose graduation is isomorphic to $\BC$ 
as an $\Gn^+$-module.
More generally for any $\Gg$-module $M$ which is generated
by a finite-dimensional $\Gb$-module $V$,
$D_\lam\htensor M$ is holonomic.

Now let us show that
$D_\lam\htensor M$ is holonomic for any $M\in\BO$.
Let $\Phi$ be a finite admissible subset of $W$.
Let $M_1$ be the $\Gg$-submodule of $M$ generated by
$\bigoplus_{w\in\Phi}M_{w\circ\lam}$,
and set $M_2=M/M_1$.
Then we have $[M_2:L(w\circ\lam)]=0$ for any $w\in\Phi$.
Hence Proposition \ref{D:loc:prop:supp} (iii) implies that
$D_\lam\htensor M_2|_{X_\Phi}=0$.
By the exact sequence
\[D_\lam\htensor M_1\to D_\lam\htensor M\to D_\lam\htensor M_2\to0
\]
$D_\lam\htensor M_1|_{X_\Phi}\to D_\lam\htensor M|_{X_\Phi}$
is surjective.
Since $D_\lam\htensor M_1$ is holonomic,
$D_\lam\htensor M|_{X_\Phi}$ is holonomic.
\qed

\Cor\label{cor:tensor's support}
For $\lam\in\KRP$ and $\mu=w\circ\lam\in\KRP$ with $w\in W$,
$\Supp(D_\lam\htensor M)\subset \overline{X}_w$ for any $M\in\BO[\mu]$.
\encor
\proof
For $x\in W$, if $[M:L(x\circ\lam)]\ne0$ then $x\circ\lam\in W(\mu)\circ\mu$.
Hence $x\in W(\mu)w$ and
Lemma \ref{lemma:ord} implies $x\ge w$. 
Then the desired result follows from 
Proposition \ref{D:loc:prop:supp} (iii).
\qed

\begin{corollary}\label{cor:tens}
Let $\lam\in\KRP$ and $\mu\in\KR$.
For $M\in\Ob(\BO[\mu])$ we have $D_\lam\htensor M=0$ unless $\mu\in
W\circ\lam$.
\end{corollary}
\proof
Assume $D_\lam\htensor M\ne0$.
Take $w\in W$ such that $X_w$ is open in $\Supp(D_\lam\htensor M)$.
Then Proposition~\ref{D:loc:prop:supp} implies $[M:L(w\circ\lam)]\ne0$.
Hence $w\circ\lam\in W\circ\mu$.
\qed

\bigskip
We define $D_\lam\htensor M$ for $\lam\in\KRP$ and  $M\in\TBO$
by
\eqn
&&D_\lam\htensor M=
\bigoplus_{\mu\in\KRP}
D_\lam\htensor P_\mu(M)
=\bigoplus_{\mu\in W\circ\lam\cap\KRP}
D_\lam\htensor P_\mu(M)\,.
\endeqn
Here the last equality follows from Corollary \ref{cor:tens}.
Then Corollary \ref{cor:tensor's support} implies that,
for any finite admissible subset $\Phi$, 
$D_\lam\htensor P_\mu(M)|_{X_\Phi}=0$ except finitely many 
$\mu\in W\circ\lam\cap\KRP$. Hence $D_\lam\htensor\,\BUL$
is a right exact functor from $\TBO$ to $\BMD(\lam)$.
The composition of functors
\[\BO\{\KR\}\longrightarrow\TBO\TO\limits^{D_\lam\htensor}\BMD(\lam)\]
coincides with the restriction of the original functor $D_\lam\htensor\ \BUL$
by Theorem \ref{thm:Kumar}.
Moreover we have the following property 
analogous to (\ref{eq:localization2}).

\Prop
For $\lam\in\KRP$, the functor
$D_\lam\htensor\ \BUL:\TBO\to \BMD(\lam)$ is a left adjoint functor of
$\,\tG:\BMD(\lam)\to\TBO$.
Namely we have an isomorphism functorial in
$\CM\in\BMD(\lam)$ and $M\in\TBO\,$$:$
\eqn
&&\Hom_{\BMD(\lam)}(D_\lam\htensor M,\CM)\cong
\Hom_\TBO(M,\tG(\CM))\,.
\endeqn
\enprop

\noindent
In particular we have a morphism functorial in $\CM\in\BMD(\lam)$
\eqn
&&
D_\lam\htensor\tG(\CM)\to\CM\,.
\endeqn

\begin{proposition}
\label{D:main:prop:tensor}
Let $\lam\in\KRP$, and $\Phi$ a finite admissible subset of $W$.
For $M\in\Ob(\TBO)$ 
we have $D_{\lam}\htensor M|_{X_\Phi}=0$ if and only
if $M\in\Ob(\TBO\{(\KR\setminus(\Phi\circ\lam)\})$.
\end{proposition}
\proof
We may assume $M\in\Ob(\BO)$.
Then it is an immediate consequence of Proposition \ref{D:loc:prop:supp}.
\qed
\begin{theorem}\label{D:loc:thm:DtensorMw}
For any $\lam\in\KRP$ and $w\in W$, there is a canonical 
isomorphism 
\[D_{\lam}\htensor
M(w\circ\lam)\iso \CM_w(\lam).\]
\end{theorem}
\proof
We have $\dim\Hom_\Gg(M(w\circ\lam),M^*(y\circ\lam))=1$ or $0$ according to
whether $w=y$ or $w\ne y$.
Hence by Proposition~\ref{D:loc:prop:supp} we have $\Supp(D_{\lam}\htensor
M(w\circ\lam))=\overline{X}_w$, and
\[
D_{\lam}\htensor M(w\circ\lam)|_{U_w}\simeq \CB_w(\lam)|_{U_w}
\simeq \CM_w(\lam)|_{U_w}.
\]
By (\ref{eq:proporty of Mw}) the isomorphism $\CM_w(\lam)|_{U_w}\simeq
D_{\lam}\htensor M(w\circ\lam)|_{U_w}$ is uniquely extended to a homomorphism
$\varphi:\CM_w(\lam)\to D_{\lam}\htensor M(w\circ\lam)$.

On the other hand, by Proposition~\ref{cohom:prop:Hom} we have a unique
non-zero homomorphism $M(w\circ\lam)\to\tG(\CM_w(\lambda))$.
Let $\psi:D_{\lam}\htensor M(w\circ\lam)\to\CM_w(\lam)$ be the corresponding
homomorphism.
By the same proposition, the composition
$D_{\lam}\htensor M(w\circ\lam)\to\CM_w(\lam)\to\CB_w(\lam)$
is non-zero.
Hence $\varphi|_{U_w}$ and $\psi|_{U_w}$ are inverse to each other
up to a non-zero constant multiple.
In particular we have $\psi\circ\varphi|_{U_w}$ is the identity endomorphism of
$\CM_w(\lam)|_{U_w}$.
Thus $\psi\circ\varphi=\id$ by (\ref{eq:proporty of Mw}).

It remains to show that $\varphi$ is an epimorphism.
Note that $\Supp(\Cok(\varphi))\subset\overline{X}_w\setminus X_w$.
Assume that $\Cok(\varphi)\ne0$.
By Lemma~\ref{lem:K-equiv}, there exists some $y\in W$ such that
$y>w$ and that $\Hom_{D_{\lam}}(\Cok(\varphi),\CB_y(\lam))\ne0$.
Since $\Cok(\varphi)$ is a quotient of $D_{\lam}\htensor M(w\circ\lam)$, we
obtain
\[
\Hom_{\Gg}(M(w\circ\lam),M^*(y\circ\lam))\simeq
\Hom_{D_{\lam}}(D_{\lam}\htensor M(w\circ\lam),\CB_y(\lam))\ne0.
\]
This implies $w=y$, which contradicts $y>w$.
Thus $\Cok(\varphi)=0$ and hence $\varphi$ is an epimorphism.
\qed

For $\lam\in\KRP$, a canonical morphism $D_{\lam}\htensor
M^*(w\circ\lam)\to \CB_w(\lam)$ is defined 
by Proposition~\ref{D:mod:cor:gamma(B_w)2} and (\ref{eq:proporty of Bw}).


\begin{lemma}\label{D:loc:lemp:DtensorM*w}
For any $\lam\in\KRP$ and $w\in W$, the canonical morphism $D_{\lam}\htensor
M^*(w\circ\lam)\to \CB_w(\lam)$ is surjective.
\end{lemma}
\proof
By Proposition \ref{D:mod:cor:gamma(B_w)2}, we have
$\Gamma(X;\CB_w(\lam))\cong M^*(w\circ\lam)$.
Since the open embedding $i:U_w\to X$ is an affine morphism
and $\CB_w(\lam)$ is a quasi-coherent $\CO_X$-module,
the natural homomorphism
$\CO_X\otimes M^*(w\circ\lam)\cong
\CO_X\otimes\Gamma(U_w;\CB_w(\lam))\to i_*i^{-1}\CB_w(\lam)\cong\CB_w(\lam)$
is surjective.
Hence the assertion follows from the fact that
$\CO_X\otimes M^*(w\circ\lam)\to\CB_w(\lam)$ decomposes into
$\CO_X\otimes M^*(w\circ\lam)\to
D_{\lam}\htensor M^*(w\circ\lam)\to \CB_w(\lam)$.
\qed

\subsection{Radon transforms}\label{Radon}

For $i\in I$, $\mu\in\Gh^*$ and $n\in\BZ$, we shall construct functors
\eq
S_{i\,*}^n:\BMD(\mu)\to\BMD(s_i\circ\mu),\qquad
S_{i\,!}^n:\BMD(\mu)\to\BMD(s_i\circ\mu),
\endeq
called the Radon transforms, and investigate their properties.
We use results in \cite{KTneg2} without giving proofs.

Fix $i\in I$.
Let $P^-$ be the algebraic group containing $B^-$ with Lie algebra
$\Gb^-\oplus\Gg_{\alpha_i}$.
We have a natural free right action of $P^-$ on $G$ (see
\cite{Kflag}).
Set $Y=G/P^-$, and let $\pi:X\to Y$ be the projection.
Then $Y$ is a separated pro-smooth scheme,
and $\pi$ is a $\BP^1$-bundle.
Set
\[
Z=X\times_Y X,\qquad Z_0=Z\setminus\Delta(X),
\]
where $\Delta:X\hookrightarrow Z$ denotes the diagonal embedding.
Let $p_r:Z_0\to X$ ($r=1, 2$) be the first and the second projections.

For a holonomic $D_\mu$-module $\CM$ we set
\eq
S^n_{i\,*}\CM=H^n(\int_{p_{1}}p_2^\BUL\CM), \qquad
S^n_{i\,!}\CM=H^n(\int_{p_{1}!}p_2^\BUL\CM).
\endeq
Since $\Omega_\pi=\CO_X(-\alpha_i)$, we have $p_1^\sharp
D_{s_i\circ\mu}=p_2^\sharp D_\mu$ (see Lemma~1.3.3 of \cite{KTneg2}), and
hence $S^n_{i\,*}$ and  $S^n_{i\,!}$ are well-defined.

By Lemma 1.4.3 and Theorem 1.5.1 of \cite{KTneg2},
we have the following proposition.

\begin{proposition}
\label{prop:Radon:1}
Let $\mu\in\Gh^*$, and let $i\in I$ and $w\in W$ such that $ws_i<w$.
\begin{tenumerate}
\item
$S^0_{i\,*}\CB_{w}(\mu)=\CB_{ws_i}(s_i\circ\mu)$ and $S^n_{i\,*}\CB_{w}(\mu)=0$
for $n\ne0$.
\item
$S^0_{i\,!}\CM_{w}(\mu)=\CM_{ws_i}(s_i\circ\mu)$ and $S^n_{i\,!}\CM_{w}(\mu)=0$
for $n\ne0$.
\item
$S^0_{i\,!}\CB_{ws_i}(\mu)=\CB_{w}(s_i\circ\mu)$
and $S^n_{i\,!}\CB_{ws_i}(\mu)=0$ for $n\ne0$.
\item
$S^0_{i\,*}\CM_{ws_i}(\mu)=\CM_{w}(s_i\circ\mu)$ and
$S^n_{i\,*}\CM_{ws_i}(\mu)=0$ for $n\ne0$.
\end{tenumerate}
\end{proposition}
\begin{theorem}
\label{thm:Radon:2}
Let $\lam\in\KR, w\in W, i\in I$ such that $ws_i<w$ and $\langle
h_i,\lam+\rho\rangle\notin\BZ_{>0}$.
Then we have
\[
\tH^n(\CB_{ws_i}(\lam))\simeq\tH^n(\CB_{w}(s_i\circ\lam)),\qquad
\tH^n(\CM_{w}(\lam))\simeq\tH^n(\CM_{ws_i}(s_i\circ\lam))
\]
for any $n\in\BZ$.
\end{theorem}
\proof
Take a finite admissible subset $\Phi$ of $W$ such that $\Phi s_i=\Phi$.
By Proposition~\ref{prop:Radon:1}, Corollary 1.6.2 of \cite{KTneg2}, and by
$\langle h_i,\lam+\rho\rangle\notin\BZ_{>0}$ we have
\[
H^n(X_\Phi;\CB_{ws_i}(\lam))\simeq
H^n(X_\Phi;\CB_{w}(s_i\circ\lam)),\quad
H^n(X_\Phi;\CM_{w}(\lam))\simeq
H^n(X_\Phi;\CM_{ws_i}(s_i\circ\lam))
\]
for any $n\in\BZ$.
Hence the desired results follow by taking the projective limit with 
respect to $\Phi$.
\qed

\begin{corollary}
\label{Radon:cor}
Let $\lambda\in\KRP, z\in W$.
Assume that $z\circ\lambda\in\KRP$.
Then we have
\[
\tH^n(\CB_w(\lambda))\simeq\tH^n(\CB_{wz^{-1}}(z\circ\lambda)), \qquad
\tH^n(\CM_w(\lambda))\simeq\tH^n(\CM_{wz^{-1}}(z\circ\lambda))
\]
for any $w\in W$ and $n\in\BZ$.
\end{corollary}
\proof
Take a reduced expression $z=s_{i_1}s_{i_2}\cdots s_{i_p}$ of $z\in W$.
It is sufficient to show
\eq
&&\tH^n(\CB_{ws_{i_p}\cdots s_{i_{j+1}}}(s_{i_{j+1}}\cdots
s_{i_p}\circ\lam))\simeq
\tH^n(\CB_{ws_{i_p}\cdots s_{i_{j}}}(s_{i_{j}}\cdots
s_{i_p}\circ\lam))\nn\\
&&\tH^n(\CM_{ws_{i_p}\cdots s_{i_{j+1}}}(s_{i_{j+1}}\cdots
s_{i_p}\circ\lam))\simeq
\tH^n(\CM_{ws_{i_p}\cdots s_{i_{j}}}(s_{i_{j}}\cdots s_{i_p}\circ\lam))\nn
\endeq
for any $j$.
By Theorem~\ref{thm:Radon:2} we have only to show
\[
(s_{i_{j+1}}\cdots
s_{i_p}\circ\lam+\rho,\alpha^\vee_{i_j})\notin
\BZ_{>0}\cup\BZ_{<0}.
\]
Since $\lam\in\KRP$ and $s_{i_p}\cdots
s_{i_{j+1}}(\alpha_{i_j})\in\Delta^+_{\re}$, we have
\[
(s_{i_{j+1}}\cdots s_{i_p}\circ\lam+\rho,\alpha^\vee_{i_j})
=(\lam+\rho,s_{i_p}\cdots s_{i_{j+1}}(\alpha_{i_j}^\vee))
\notin\BZ_{<0}.
\]
On the other hand, since $z\circ \lam\in\KRP$ and $s_{i_1}\cdots
s_{i_{j-1}}(\alpha_{i_j})\in\Delta^+_{\re}$, we have
\[
(s_{i_{j+1}}\cdots s_{i_p}\circ\lam+\rho,\alpha^\vee_{i_j})
=-(z\circ\lam+\rho,s_{i_1}\cdots s_{i_{j-1}}(\alpha_{i_j}^\vee))
\notin\BZ_{>0}.
\]
The proof is completed.
\qed

\subsection{Global sections of $\CM_w(\lambda)$}
For $\lambda\in\KRP$ and $w\in W$, we denote by
\eq
\varphi^\lam_w:M(w\circ\lam)\to\tG(\CM_w(\lam))
\endeq
a non-zero morphism in $\TBO$
(see Proposition~\ref{cohom:prop:Hom}).
Note that $\varphi^\lam_w$ is unique up to a non-zero constant multiple.
The aim of this section is to prove that it is a monomorphism.

Let $i\in I$.
Let $\pi:X\to Y$ be the $\BP^1$-bundle as in \S \ref{Radon}.
Assume that $w\in W$ satisfy $ws_i>w$.
Set $Y_w=\pi(X_w)$. Then $Y_w$ is an affine scheme.
$X_{iw}=\pi^{-1}(Y_w)=X_w\sqcup X_{ws_i}$ and
$X_{ws_i}\to Y_w$
is an isomorphism.
Hence $X_{ws_i}$ is a closed hypersurface of $X_{iw}$.
Let $j:X_{iw}\to X$ be the inclusion.
Let $\lam$ be an element of $\Gh^*$ satisfying
$(\lam+\rho,\alpha_i^\vee)\in\BZ$.
Then there exists an $N^+$-equivariant line bundle $L$ on $X_{iw}$
such that $j^\sharp D_\lam\simeq D_{X_{iw}}(L)$.
Let $j_0:X_w\to X_{iw}$ be the open embedding and $i_0:X_{ws_i}\to X_{iw}$
the closed embedding.
We have the exact sequences of holonomic $D_{X_{iw}}$-modules
\eqn
&0\to\int_{j_0\,!}\CO_{X_{ws_i}}\to\int_{i_0\,!}\CO_{X_{w}}
\to\CO_{X_{iw}}\to0,&\\
&0\to\CO_{X_{iw}}\to\int_{i_0}\CO_{X_{w}}\to\int_{i_0}\CO_{X_{ws_i}}\to0.&
\endeqn
Since $j$ is an affine morphism,
$H^0\int_{j\,!}$ and $H^0\int_{j}$ are exact functors.
Tensoring $L$ to the exact sequences above,
and applying the exact functors $H^0\int_{j\,!}$
and $H^0\int_{j}$, we obtain exact sequences in $\BMD(\lam)$
\eq\label{pr:em}
&&
\begin{array}{c}
0\longrightarrow\CM_{ws_i}(\lam)\buildrel\iota\over\longrightarrow
\CM_{w}(\lam)\longrightarrow\CL\longrightarrow0\,,\\[5pt]
0\longrightarrow\CL'\longrightarrow\CB_w(\lam)
\longrightarrow\CB_{ws_i}(\lam)\longrightarrow0\,.
\end{array}
\endeq
where $\CL=H^0\int_{j\,!}L$ and $\CL=H^0\int_{j}L$.
We have $\CL|_U\simeq\CL'|_U$ for any open set containing $X_{iw}$
as a closed subset.

\Lemma\label{D:Mw:comm}
Assume that $\lam\in\KRP$, $i\in I$ and $w\in W$ satisfy $ws_i>w$ and
$\alpha_i\in\Delta(\lam)$.
Let
$\iota_*:\tG(\CM_{ws_i}(\lam))\to\tG(\CM_{w}(\lam))$
be the monomorphism in $\TBO$ induced by the monomorphism
$\iota:\CM_{ws_i}(\lam)\to\CM_{w}(\lam)$ in $(\ref{pr:em})$.
Let $j:M(ws_i\circ\lam)\to M(w\circ\lam)$ be the injective homomorphism of
$\Gg$-modules {\rm(}see {\rm Proposition~\ref{D:mod:prop:[M:L]})}.
Then the following diagram in $\TBO$
is commutative up to a non-zero constant multiple.
\[
\begin{array}{ccc}
M(ws_i\circ\lam)&\maprightu{\varphi^\lam_{ws_i}}&
\tG(\CM_{ws_i}(\lam))\\
\mapdownl{j}&&\mapdownr{\iota_*}\\
M(w\circ\lam)&\maprightu{\varphi^\lam_{w}}&
\tG(\CM_{w}(\lam))\,.
\end{array}
\]
\enlemma
\proof
It is sufficient to show the following two statements.
\eq
&&\dim\Hom_\TBO(M(ws_i\circ\lambda),\tG(\CM_w(\lambda)))
=1.\label{aa1}\\
&&\varphi_{w}^\lam\circ j\ne0.\label{aa2}
\endeq

We first show (\ref{aa1}).
The first exact sequence in (\ref{pr:em})
\[
0\to \CM_{ws_i}(\lam)\to\CM_{w}(\lam)\to\CL\to 0
\]
induces exact sequences
\eq
&&
0\to \Hom_\TBO(M(ws_i\circ\lam),\tG(\CM_{ws_i}(\lam)))
\to\Hom_\TBO(M(ws_i\circ\lam),\tG(\CM_{w}(\lam)))\\
&&\qquad\qquad\qquad\qquad\qquad\qquad\qquad\qquad
\qquad\qquad\to\Hom_\TBO(M(ws_i\circ\lam),\tG(\CL)).\nonumber
\endeq
Since $\dim\Hom_\TBO(M(ws_i\circ\lam),\tG(\CM_{ws_i}(\lam)))=1$
by Proposition~\ref{cohom:prop:Hom},
it is sufficient to show
$\Hom_\TBO(M(ws_i\circ\lam),\tG(\CL))=0$.
Since 
$\CM_{ws_i}(\lam)\cong\CD_\lam\htensor M(ws_i\circ\lam)$ 
by Theorem \ref{D:loc:thm:DtensorMw},
we reduce this to
$\Hom(\CM_{ws_i}(\lam),\CL)=0$.
Set $\Phi=\{x\in W\,;\,x\leq ws_i\}$
and $\Psi=\Phi\setminus\{ws_i\}$.
They are finite admissible subsets of $W$, and 
$X_\Phi\cap \overline{X}_w=X_{iw}=X_w\sqcup X_{ws_i}$.
Then (\ref{pr:em}) induces
an exact sequence:
\[
0\to\CL|_{X_\Phi}\to\CB_w(\lam)|_{X_\Phi}
\to\CB_{ws_i}(\lam)|_{X_\Phi}\to0.
\]
Since $\CM_{ws_i}(\lam)|_{X_\Psi}=0$,
we have by (\ref{eq:proporty of Mw}) and (\ref{eq:proporty of Bw}) 
\eqn
\Hom(\CM_{ws_i}(\lam),\CL)&\simeq&
\Hom(\CM_{ws_i}(\lam)|_{X_\Phi},\CL|_{X_\Phi})\\
&\subset&\Hom(\CM_{ws_i}(\lam)|_{X_\Phi},\CB_w(\lam)|_{X_\Phi})\\
&\simeq&\Hom(\CM_{ws_i}(\lam)|_{X_\Psi},\CB_w(\lam)|_{X_\Psi})
=0.
\endeqn
Thus the proof of (\ref{aa1}) is completed.

\medskip
Let us prove (\ref{aa2}).
Consider the chain of morphisms
\[
\psi:M(w\circ\lam)\to\tG(\CM_w(\lam))\iso\tG(\CM_{ws_i}(s_i\circ\lam)
)\to
\tG(\CB_{ws_i}(s_i\circ\lam)).
\]
Here the middle isomorphism follows from Theorem \ref{thm:Radon:2},
because
\[(s_i\circ\lam+\rho,\alpha_i^\vee)=
-(\lam+\rho,\alpha_i^\vee)\not\in\BZ_{>0}.\]
In order to prove (\ref{aa2}),
it is sufficient to show that the composition of
\[
M(ws_i\circ\lam)\maprightu{j}M(w\circ\lam)
\maprightu{\psi}\tG(\CB_{ws_i}(s_i\circ\lam))
\]
is non-zero.
Assuming that $\psi$ is a non-zero homomorphism for a while,
we shall finish the proof.
Since $\ch(\tG(\CB_{ws_i}(s_i\circ\lam)))=\ch(M(w\circ\lam))$,
we have $[\Ker(\psi):L(\zeta)]=[\Cok(\psi):L(\zeta)]$ for any
$\zeta$. 
Assume that $\psi\circ j=0$.
Then $\Ker(\psi)$ contains $M(ws_i\circ\lam)$ as a submodule.
Hence we have
\[[\Cok(\psi):L(ws_i\circ\lam)]=[\Ker(\psi):L(ws_i\circ\lam)]>0,\]
Since $\Cok(\psi)_{w\circ\lam}=0$,
$[\Cok(\psi):L(\zeta)]\ne0$ only if $\zeta=y\circ\lam$ for some
$y\in W$ such that $y>w$.
If $[M(\zeta):L(ws_i\circ\lam)]\ne0$ and $[\Cok(\psi):L(\zeta)]\ne0$, then
we have $\zeta=y\circ\lam$ with $ws_i\ge y>w$ 
(by Corollary \ref{D:mod:cor:[M:L]}), and hence
$\zeta=ws_i\circ\lam$.
Hence Lemma~\ref{KM:cat:dimHom} implies
\[
\dim\Hom_\Gg(M(ws_i\circ\lam),\Cok(\psi)^*)>0.
\]
Thus we obtain
\[
\dim\Hom_\Gg(M(ws_i\circ\lam),\tG(\CB_{ws_i}(s_i\circ\lam))^*)\ge
\dim\Hom_\Gg(M(ws_i\circ\lam),\Cok(\psi)^*)>0.
\]
Applying Corollary~\ref{D:mod:cor:gamma(B_w)1} (iii), we have
\[
ws_i\circ\lam\in w\circ\lam-\sum_{\alpha\in\Delta^+\cap
ws_i\Delta^+}\BZ_{\ge0}\alpha,
\]
and hence
\eqn
(\lam+\rho,\alpha_i^\vee)\alpha_i=w^{-1}(w\circ\lam-ws_i\circ\lam)
&\in&w^{-1}\sum_{\alpha\in\Delta^+\cap ws_i\Delta^+}\BZ_{\ge0}\alpha
=\sum_{\alpha\in w^{-1}\Delta^+\cap s_i\Delta^+}\BZ_{\ge0}\alpha\\
&\subset&
\sum_{\alpha\in \Delta^+\setminus\{\alpha_i\}}\BZ_{\ge0}\alpha.
\endeqn
This is a contradiction.
Hence $\psi\circ j\ne0$.

It remains to prove that $\psi$ does not vanish.
In order to see this, it is sufficient to show the injectivity of
the homomorphism
\[
\Hom_\TBO(M(w\circ\lam),\tG(\CM_{ws_i}(s_i\circ\lam)))\to
\Hom_\TBO(M(w\circ\lam),\tG(\CB_{ws_i}(s_i\circ\lam)))\,.
\]
Let $\CN$ be the kernel of
the morphism
$\CM_{ws_i}(s_i\circ\lam)\to\CB_{ws_i}(s_i\circ\lam)$.
By the exact sequence
\eqn
&&0\to
\Hom_\TBO(M(w\circ\lam),\tG(\CN))\to
\Hom_\TBO(M(w\circ\lam),\tG(\CM_{ws_i}(s_i\circ\lam)))\\
&&\hspace{180pt}
\to
\Hom_\TBO(M(w\circ\lam),\tG(\CB_{ws_i}(s_i\circ\lam))),
\endeqn
we can reduce the assertion to
$[\tG(\CN):L(w\circ\lam)]=0$.
Assume the contrary.
Then by Proposition \ref{D:mod:prop:tH1} (iii),
there exists some $x\in W$ such that
\eq
[M(xs_i\circ\lam):L(w\circ\lam)]\ne0\quad\mbox{and}\label{eq:1}\\
X_x\subset\Supp\CN\subset\ol{X_{ws_i}}\setminus X_{ws_i}.\label{eq:2}
\endeq
Then (\ref{eq:1}) implies $xs_i\le w$ by Corollary \ref{D:mod:cor:[M:L]},
and (\ref{eq:2}) implies $ws_i<x$.
Hence we obtain
$xs_i\le w<ws_i<x$, which implies
$l(xs_i)\le l(w)<l(ws_i)<l(x)\le l(xs_i)+1$.
This is a contradiction.
Hence we have proved the non-vanishing of $\psi$.
\qed

Now we are ready to prove the following result.
\begin{proposition}\label{pro:inj}
Let $\lambda\in\KRP$, $w\in W$.
Then the morphism
\[
\varphi_w^\lam:M(w\circ\lam)\to\tG(\CM_w(\lam))
\]
is a monomorphism.
\end{proposition}
\proof
Take $x\in W(w\circ\lam)$ such that $\lam'=x^{-1}w\circ\lam\in \KRP$.
Then $w\circ\lam=x\circ\lam'$ and $x\in W(\lam')$.
By Corollary~\ref{Radon:cor} we have
$\tG(\CM_w(\lam))\simeq\tG(\CM_x(\lam'))$.
The composition of
\[
M(w\circ\lam) =
M(x\circ\lam')\maprightu{{\varphi_x^{\lam'}}}\tG(\CM_x(\lam')
)\simeq\tG(\CM_w(\lam))
\]
coincides with $\varphi_w^\lam$ up to a non-zero scalar multiple by
Proposition~\ref{cohom:prop:Hom}.
Hence we may assume from the beginning that $w\in W(\lam)$.

For $y, x\in W(\lam)$ such that $y\ge_\lam x$, we denote by
$f_x^y:M(y\circ\lam)\to M(x\circ\lam)$ a non-zero homomorphism of
$\Gg$-modules.

Assume that for any $w\in W(\lam)$ there exists a monomorphism
\[
F_w:\tG(\CM_w(\lam))\to\tG(\CM_e(\lam))
\]
in $\TBO$
such that
the diagram
\eq\label{M_w:eq:e}
\begin{array}{ccc}
M(w\circ\lam)
&\maprightu{\varphi_w^\lam}
&\tG(\CM_w(\lam))\\
\mapdownl{f_e^w}
&
&\mapdownr{F_w}\\
M(\lam)
&\maprightd{\varphi_e^\lam}
&\tG(\CM_e(\lam))
\end{array}
\endeq
is commutative.
If $\Ker \varphi_e^\lam\ne0$, then there exists some $w\in W$ such that
$\Image f_e^w\subset\Ker \varphi_e^\lam$.
Since $\varphi_w^\lam\ne0$,
this contradicts the injectivity of $F_w$.
Hence $\varphi_e^\lam$ is injective.
Thus $\varphi_w^\lam$ is a monomorphism 
by the commutativity of (\ref{M_w:eq:e}).

Therefore it remains to show that for any $w\in W(\lam)$ there exists a
monomorphism
$F_w:\tG(\CM_w(\lam))\to\tG(\CM_e(\lam))$
in $\TBO$
such that
the diagram
(\ref{M_w:eq:e})
is commutative.

For $w\in W(\lam)$, take a reduced expression
$w=s_{\beta_1}\cdots s_{\beta_r}\,(\beta_k\in\Pi(\lam))$.
Set $w_k=s_{\beta_1}\cdots s_{\beta_{k}}$.
Then $w_k=w_{k-1} s_{\beta_{k}}$, and $w_{k-1}\beta_{k}\in\Delta^+$.

Then it is sufficient to show that for any $k$ there exists a
monomorphism
\[
F_k:\tG(\CM_{w_k}(\lam))\to\tG(\CM_{w_{k-1}}(\lam))
\]
in $\TBO$
such that
the diagram
\eq
\begin{array}{ccc}
M(w_k\circ\lam)
&\maprightu{\varphi_{w_k}^\lam}
&\tG(\CM_{w_k}(\lam))\\
\mapdownl{
f_{w_{k-1}}^{w_k}
}
&
&\mapdownr{F_k}\\
M(w_{k-1}\circ\lam)
&\maprightd{\varphi_{w_{k-1}}^\lam}
&\tG(\CM_{w_{k-1}}(\lam))
\end{array}
\endeq
is commutative (see Theorem \ref{KM:emb:thm:embed}).
By Lemma~\ref{km:int:lem:8} we can take $x\in W$ such that
$\lam'=x\circ\lam\in\KRP, x(\beta_k)=\alpha_i\in\Pi$.
We have $\alpha_i\in\Delta(\lam')$.
Set $y=w_{k-1}x^{-1}$.
Then $y\alpha_i=w_{k-1}\beta_k\in \Delta^+$ and hence $ys_i>y$.
We have $w_{k-1}\circ\lam=y\circ\lam',  w_{k}\circ\lam=ys_i\circ\lam'$,
and hence $M(w_{k-1}\circ\lam)=M(y\circ\lam'),
M(w_{k}\circ\lam)=M(ys_i\circ\lam')$.
Moreover, by Corollary~\ref{Radon:cor}, we have
$\tG(\CM_{w_{k-1}}(\lam))=\tG(\CM_{y}(\lam')),
\tG(\CM_{w_{k}}(\lam))=\tG(\CM_{ys_i}(\lam'))$.
Hence the desired result follows from Lemma~\ref{D:Mw:comm}.
\qed

\subsection{Correspondence of $\GGg$-modules and $D$-modules}
We shall prove the following theorem
on a partial correspondence between $\BMD(\lam)$ and
$\TBO(\lam)$.
This theorem says in particular that the composition
$\BMD(\lam)\TO\limits^{\tG}\TBO(\lam)
\TO\limits^{D_\lam\htensor}\BMD(\lam)$
is an equivalence of categories.
Hence $\BMD(\lam)$ is equivalent to a full subcategory of
$\TBO(\lam)$.
Moreover,
$\CM_w(\lam)$, $\CB_w(\lam)$ and $\CL_w(\lam)$ in the category
$\BMD(\lam)$ correspond to
$M(w\circ\lam)$, $M^*(w\circ\lam)$ and $L(w\circ\lam)$
in $\TBO(\lam)$, respectively.
\begin{theorem}\label{Main theorem}
Let $\lam\in\KRP$.
\begin{tenumerate}
\item
$\tH^n(\CM)=0$ for any $\CM\in\Ob(\BMD(\lam))$ and $n\ne0$.
\item
$\tG(\CB_w(\lam))\simeq M^*(w\circ\lam)$ and
$\tG(\CM_w(\lam))\simeq M(w\circ\lam)$ for any $w\in W$.
\item\label{prop:tens}
$\tG(\CL_w(\lam))\simeq L(w\circ\lam)$ for any $w\in W$.
\item $D_\lam\htensor \tG(\CM)\iso\CM$ for any $\CM\in\Ob(\BMD(\lam))$.
\end{tenumerate}
\end{theorem}

\proof
Note that the first statement in (ii) is already proved
(Proposition \ref{D:mod:prop:tH2}).
We first show (i), (ii) and (iv).
It is sufficient to show the following statements ({\rm a}), ({\rm b})
and ({\rm c}) for
any finite admissible subset $\Phi$ of $W$.

\newenvironment{alphenumerate}{
  \begin{enumerate}
  \renewcommand{\labelenumi}
  {\rm
  {(\alph{enumi})}
  }
  }{\end{enumerate}}

\begin{alphenumerate}
\item
For $\CM\in\Ob(\BMD(\lam))$, we have
$\tH^n(\CM)\in\Ob(\TBO\{(W\setminus\Phi)\circ\lam\})$ for any $n\ne0$.
\item
For $w\in \Phi$,
the cokernel of the monomorphism 
$\varphi_w^\lam:M(w\circ\lam)\to\tG(\CM_w(\lam))$
belongs to $\TBO\{(W\setminus\Phi)\circ\lam\}$.
\item
For $\CM\in\Ob(\BMD(\lam))$ and a morphism $\varphi:M\to\tG(\CM)$
in $\TBO$, assume that $\Ker(\varphi)$ and $\Cok(\varphi)$ belong to
$\TBO\{(W\setminus\Phi)\circ\lam\}$.
Then the canonical homomorphism
$D_\lam\htensor M|_{X_\Phi}\to \CM|_{X_\Phi}$
is an isomorphism.
\end{alphenumerate}
Fixing $\Phi$, we shall show the following statements 
(a$)_{\Psi}$, (b$)_{\Psi}$,
(c$)_{\Psi}$ for a finite admissible subset $\Psi$ of $W$ such that
$\Psi\subset\Phi$ by induction on $\sharp(\Phi\setminus\Psi)$.
Note that (a)=(a$)_\emptyset$, (b)=(b$)_\emptyset$ and (c)=(c$)_\emptyset$.

\newenvironment{penumerate}{
  \begin{enumerate}
  \renewcommand{\labelenumi}
  {\rm
  {(\alph{enumi}$)_{\Psi}$}
  }
  }{\end{enumerate}}

\begin{penumerate}
\item
Let $\CM\in\Ob(\BMD(\lam))$ such that $\Supp(\CM)\cap X_\Psi=\emptyset$.
Then we have
$\tH^n(\CM)\in\Ob(\TBO\{(W\setminus\Phi)\circ\lam\})$ for any $n>0$.
\item
For $w\in \Phi\setminus\Psi$,
the cokernel of the monomorphism 
$\varphi_w^\lam:M(w\circ\lam)\to\tG(\CM_w(\lam))$
belongs to $\TBO\{(W\setminus\Phi)\circ\lam\}$.
\item
Let $\CM\in\Ob(\BMD(\lam))$ such that $\Supp(\CM)\cap X_\Psi=\emptyset$.
Assume that a morphism $\varphi:M\to\tG(\CM)$ in $\TBO$ satisfies
$\Ker(\varphi)$, $\Cok(\varphi)\in\Ob(\TBO\{(W\setminus\Phi)\circ\lam\})$.
Then the canonical homomorphism
$
D_\lam\htensor M|_{X_\Phi}\to \CM|_{X_\Phi}
$
is an isomorphism.
\end{penumerate}

In the case $\Phi=\Psi$, (a$)_\Phi$ follows from
Proposition~~\ref{D:mod:prop:tH1} (iv), (b$)_\Phi$ is trivial, and
(c$)_\Phi$ is a consequence of Proposition~\ref{D:mod:prop:tH1} (iv) and
Proposition~\ref{D:main:prop:tensor}.

Assume $\sharp(\Phi\setminus\Psi)>0$.
Take a minimal element $y$ of the set $\Phi\setminus\Psi$, and set
$\Psi'=\Psi\cup\{y\}$.
Then $\Psi'$ is a finite admissible subset of $W$ satisfying
$\Psi\subset\Psi'\subset\Phi$ and
$\sharp(\Phi\setminus\Psi')=\sharp(\Phi\setminus\Psi)-1$.
Hence we may assume the statements (a$)_{\Psi'}$, (b$)_{\Psi'}$,
(c$)_{\Psi'}$ by the hypothesis of induction.

Set $\CL=\CB_y(\lam)/\CL_y(\lam)$.
Since $\Supp(\CL)\cap X_{\Psi'}=\emptyset$, (a$)_{\Psi'}$ implies
$\tH^n(\CL)\in\Ob(\TBO\{(W\setminus\Phi)\circ\lam\})$ for $n\ne0$.
By considering the long exact sequence associated to the short exact sequence
\[
0\to\CL_y(\lam)\to\CB_y(\lam)\to\CL\to0\,,
\]
we obtain
$\tH^n(\CL_y(\lam))\in\Ob(\TBO\{(W\setminus\Phi)\circ\lam\})$ for any
$n\ge2$ and an exact sequence
\eq\label{mainthm:eq:1}
\tG(\CB_y(\lam))\to\tG(\CL)\to\tH^1(\CL_y(\lam))\to0.
\endeq
Consider the following natural commutative diagram.
\[
\begin{array}{ccc}
D_{\lam}\htensor\tG(\CB_y(\lam))|_{X_\Phi}
&\maprightu{f_1}
&D_{\lam}\htensor\tG(\CL)|_{X_\Phi}\\
\mapdownl{h_1}
&
&\mapdownr{h_2}\\
\CB_y(\lam)|_{X_\Phi}
&
\maprightu{f_2}
&\CL|_{X_\Phi}\,.
\end{array}
\]
Then  $h_2$ is an isomorphism by (c$)_{\Psi'}$, $h_1$ is surjective by
Lemma~\ref{D:loc:lemp:DtensorM*w}, and $f_2$ is obviously surjective.
Thus $D_{\lam}\htensor\tG(\CB_y(\lam))|_{X_\Phi}
\to D_{\lam}\htensor\tG(\CL)|_{X_\Phi}$ is surjective.
Hence we have $D_{\lam}\htensor\tH^1(\CL_y(\lam))|_{X_\Phi}=0$ 
by (\ref{mainthm:eq:1}).
Then we obtain
$\tH^1(\CL_y(\lam))\in\Ob(\TBO\{(W\setminus\Phi)\circ\lam\})$ by
Proposition~\ref{D:main:prop:tensor}.
We have thus proved that
$\tH^n(\CL_y(\lam))\in\Ob(\TBO\{(W\setminus\Phi)\circ\lam\})$ for 
$n>0$.

Set $\CL'=\Ker(\CM_y(\lam)\to\CL_y(\lam))$.
Since $\Supp(\CL')\cap X_{\Psi'}=\emptyset$, (a$)_{\Psi'}$ implies
$\tH^n(\CL')\in\Ob(\TBO\{(W\setminus\Phi)\circ\lam\})$ for any
$n>0$.
By considering the long exact sequence associated to the short exact sequence
\[
0\to\CL'\to\CM_y(\lam)\to\CL_y(\lam)\to0\,,
\]
we obtain
\eq
\label{eq:HnMy}
\tH^n(\CM_y(\lam))\in\Ob(\TBO\{(W\setminus\Phi)\circ\lam\})\quad
\mbox{ for any $n>0$.}
\endeq

Let us show (a$)_{\Psi}$.
By Lemma~\ref{lem:K-equiv} there exists 
a morphism $f:\CM_y(\lam)^{\oplus r}\to\CM$
whose restriction
$\CM_y(\lam)^{\oplus r}|_{X_{\Psi'}}\iso\CM|_{X_{\Psi'}}$
is an isomorphism.
Setting $\CN=\Image(f)$, $\CN_1=\Ker(f)$ and $\CN_2=\Cok(f)$, 
we obtain exact
sequences
\eq
&0\to\CN_1\to\CM_y(\lam)^{\oplus r}\to\CN\to0,&\nn\\
&0\to\CN\to\CM\to\CN_2\to0.&\nn
\endeq
Note that $\Supp(\CN_i)\cap X_{\Psi'}=\emptyset$ for $i=1$, $2$.
Let $n>0$.
By (a$)_{\Psi'}$ and (\ref{eq:HnMy}), the objects
$\tH^{n+1}(\CN_1)$, $\tH^n(\CN_2)$ and $\tH^n(\CM_y(\lam))$
belong to $\TBO\{(W\setminus\Phi)\circ\lam\}$.
Hence $\tH^n(\CM)$ also belongs to $\TBO\{(W\setminus\Phi)\circ\lam\}$ by
the exact sequences
\eq
&\tH^n(\CM_y(\lam))^{\oplus
r}\to\tH^n(\CN)\to\tH^{n+1}(\CN_1),&\nn\\
&\tH^n(\CN)\to\tH^n(\CM)\to\tH^{n}(\CN_2).&\nn
\endeq
The statement (a$)_{\Psi}$ is proved.

We next show (b$)_{\Psi}$.
By (b$)_{\Psi'}$ we have only to deal with the case $w=y$.
By Lemma~\ref{D:cat:lemma:mult} we have, for any $\mu\in \KR$,
\eqn
\sum_n(-1)^n\ch(P_\mu\tH^n(\CM_y(\lam)))
&=&\sum_n(-1)^n\ch(P_\mu\tH^n(\CB_y(\lam)))
=\ch(P_\mu\tG(\CB_y(\lam)))\\
&=&
\ch\big(P_\mu(M(y\circ\lam))\big)\,.
\endeqn
Since $\tH^n(\CM_y(\lam))$ belongs to
$\TBO\{(W\setminus\Phi)\circ\lam\}$ for any $n>0$ by (a$)_{\Psi}$, 
we see that
$\ch(P_\mu(\Cok \varphi^\lam_y))=
\ch(P_\mu\tG(\CM_y(\lam)))-\ch(P_\mu(M(y\circ\lam)))$ is in
$\sum_{x\in W\setminus\Phi}\BZ\ch(L(x\circ\lam))$.
The statement (b$)_{\Psi}$ is proved.

Let us show (c$)_{\Psi}$.
Take  $r, f, \CN, \CN_1, \CN_2$ for $\CM$ as in the proof of (a$)_{\Psi}$.
Set $N=\varphi^{-1}(\Image(\tG(\CN)\to\tG(\CM)))$ and
$N_2=M/N$.
Then we obtain the following commutative diagrams whose rows are exact.
\eq
\begin{array}{ccccccccc}
&&&&&&0&&\\
&&&&&&\mapdownl{}&&\\
0
&
\to
&
N
&
\maprightu{}
&
M
&
\maprightu{}
&
N_2
&
\maprightu{}
&
0\\
&&\mapdownl{\psi}
&&\mapdownl{\varphi}
&&\mapdownr{\varphi'}
&&\mapdownr{}
\\
0
&
\to
&
\tG(\CN)
&
\maprightu{}
&
\tG(\CM)
&
\maprightu{h}
&
\tG(\CN_2)
&
\maprightu{}
&
\tH^1(\CN)
\end{array}\nn
\endeq
\eq
\begin{array}{ccccccccc}
&&
D_{\lam}\htensor N
&
\maprightu{}
&
D_{\lam}\htensor M
&
\maprightu{}
&
D_{\lam}\htensor N_2
&
\to
&
0\phantom{.}\\
&&\mapdownl{}
&&\mapdownl{}
&&\mapdownr{}
\\
0&
\to
&
\CN
&
\maprightu{}
&
\CM
&
\maprightu{}
&
\CN_2
&
\to
&
0.
\end{array}\nn
\endeq
By the definition of $N_2$, $\varphi'$ is a monomorphism.
By (a$)_{\Psi}$, we have
$\tH^1(\CN)\in\Ob(\TBO\{(W\setminus\Phi)\circ\lam\})$.
Hence by the exact sequences
\eq
&&
0\to\Ker\psi\to\Ker\varphi\to0\to
\Cok\psi\to\Cok\varphi\to\Cok\varphi'\to
\tH^1(\CN),
\nn
\endeq
and by the assumption on $\varphi$, we obtain
\eq
\label{eq:Maintheorem:a}
\Ker(\psi), \;\Cok(\psi),\;
\Cok(\varphi')\in\Ob(\TBO\{(W\setminus\Phi)\circ\lam\}).
\endeq
Since $\Supp(\CN_2)\cap X_{\Psi'}=\emptyset$, the morphism
$D_{\lam}\htensor N_2|_{X_\Phi}\to\CN_2|_{X_\Phi}$ is an
isomorphism by (c$)_{\Psi'}$ and (\ref{eq:Maintheorem:a}).
Hence it is sufficient to show that 
$D_{\lam}\htensor N|_{X_\Phi}\to\CN|_{X_\Phi}$ is an isomorphism.

Set $N_0=\psi^{-1}(\Image(\tG(\CM_y(\lam)^{\oplus r})\to\tG(\CN)))$ 
and let $\psi_{0}:N_0\to\tG(\CN)$ be the
restriction of $\psi$.
Since $\tH^1(\CN_1)\in\Ob(\TBO\{(W\setminus\Phi)\circ\lam\})$ by
(a$)_{\Psi'}$, we have
\eq
\label{eq:Maintheorem:b}
N/N_0,\;
\Ker\psi_{0},\;
\Cok\psi_{0}
\in
\Ob(\TBO\{(W\setminus\Phi)\circ\lam\})
\endeq
by (\ref{eq:Maintheorem:a}) and the exact sequences
\eq
&&0\to N/N_0\to\tH^1(\CN_1),
\qquad
0\to\Ker\psi_{0}\to\Ker\psi,
\nn
\\
&&\qquad\qquad
N/N_0\to\Cok\psi_{0}\to\Cok\psi\to0.
\nn
\endeq
Proposition~\ref{D:main:prop:tensor} implies $D_{\lam}\htensor
(N/N_0)|_{X_\Phi}=0$, and hence $D_{\lam}\htensor N_0|_{X_\Phi}\to
D_{\lam}\htensor N|_{X_\Phi}$ is surjective.
Since $D_{\lam}\htensor N_0|_{X_\Phi}\to\CN|_{X_\Phi}$ decomposes into
$D_{\lam}\htensor N_0|_{X_\Phi}\to D_{\lam}\htensor
N|_{X_\Phi}\to\CN|_{X_\Phi}$, it is sufficient to show that
$D_{\lam}\htensor N_0|_{X_\Phi}\to\CN|_{X_\Phi}$ is an
isomorphism.

Let $R$ be the fiber product
of $\tG(\CM_y(\lam)^{\oplus r})$ and $N_0$ over $\tG(\CN)$, 
and consider the following commutative diagrams whose rows are
exact:
\eq
\begin{array}{ccccccccc}
0
&
\to
&
\tG(\CN_1)
&
\maprightu{}
&
R
&
\maprightu{}
&
N_0
&
\maprightu{}
&
0\\
&&\mapdownl{\id}
&&\mapdownl{\psi'}
&&\mapdownr{\psi_{0}}
&&
\mapdownl{}
\\
0
&
\to
&
\tG(\CN_1)
&
\maprightu{}
&
\tG(\CM_y(\lam)^{\oplus r})
&
\maprightu{k}
&
\tG(\CN)
&
\maprightu{}
&
\tH^1(\CN_1)\,,
\end{array}\nn
\endeq
\eq
\begin{array}{ccccccccccc}
&&
D_{\lam}\htensor\tG(\CN_1)
&
\maprightu{}
&
D_{\lam}\htensor R
&
\maprightu{}
&
D_{\lam}\htensor N_0
&
\to
&
0\phantom{\ .}\\
&&\mapdownl{}
&&\mapdownl{}
&&\mapdownr{}
\\
0&\maprightu{}
&
\CN_1
&
\maprightu{}
&
\CM_y(\lam)^{\oplus r}
&
\maprightu{}
&
\CN
&
\to
&
0\ .
\end{array}\nn
\endeq
By (c$)_{\Psi'}$ the morphism
$D_{\lam}\htensor\tG(\CN_1)|_{X_\Phi}\to\CN_1|_{X_\Phi}$
is an isomorphism.
Hence it is sufficient to show that the morphism
$D_{\lam}\htensor R|_{X_\Phi}\to\CM_y(\lam)^{\oplus
r}|_{X_\Phi}$ is an isomorphism.

Since $\Ker\psi'$ and $\Cok\psi'$ are isomorphic to subobjects of
$\Ker\psi_{0}$ and $\Cok\psi_{0}$ respectively, we have
\eq
\label{eq:Maintheorem:c}
\Ker(\psi'),\; \Cok(\psi')
\in\Ob(\TBO\{(W\setminus\Phi)\circ\lam\})
\endeq
by (\ref{eq:Maintheorem:b}).
Thus we have reduced (c$)_\Psi$ to the case $\CM=\CM_y(\lam)^{\oplus r}$.

The statement (b$)_{\Psi}$ implies that
$\varphi^\lam_y:M(y\circ\lam)\to\tG(\CM_y(\lam))$ is a
monomorphism such that
$\Cok\varphi^\lam_y\in\Ob(\TBO\{(W\setminus\Phi)\circ\lam\})$.
Let $R'$ be the Cartesian product of
$M(y\circ\lam)^{\oplus r}$ and $R$ over
$\tG(\CM_y(\lam)^{\oplus r})\,$:
\eqn
&&\begin{array}{ccc}
R'&\mathop{\longhookrightarrow}\limits^\eta&R\\
\mapdownl{\psi''}&&\mapdownl{\psi'}\\
M(y\circ\lam)^{\oplus r}&
\mathop{\longhookrightarrow}\limits^{\varphi^\lam_y{}^{\oplus r}}
&\tG(\CM_y(\lam)^{\oplus r}).
\end{array}
\endeqn
Then $\Cok \eta$ belongs to
$\TBO\{(W\setminus\Phi)\circ\lam\}$, and hence
$D_\lam\htensor \Cok \eta|_{X_\Phi}=0$
by Proposition \ref{D:main:prop:tensor}.
Hence we obtain
\eq\label{pro:sur}
&&
\mbox{$D_\lam\htensor R'|_{X_\Phi}\to D_\lam\htensor R|_{X_\Phi}$ 
is surjective.}
\endeq
On the other hand,
the cokernel of $\psi''$ belongs to $\TBO\{(W\setminus\Phi)\circ\lam\}$,
which implies that $\psi''$ is surjective.
Since the kernel of $\psi''$ also belongs to 
$\TBO\{(W\setminus\Phi)\circ\lam\}$, 
$D_\lam\htensor \Ker \psi''|_{X_\Phi}=0$. Hence we have
\eqn
&&
D_\lam\htensor R'|_{X_\Phi}
\iso D_\lam\htensor M(y\circ\lam)^{\oplus r}|_{X_\Phi}
\cong\CM_y(\lam)^{\oplus r}|_{X_\Phi}.
\endeqn
Since the isomorphism 
$D_\lam\htensor R'|_{X_\Phi}
\iso \CM_y(\lam)^{\oplus r}|_{X_\Phi}$
factors through $D_\lam\htensor R|_{X_\Phi}$,
(\ref{pro:sur}) implies that
$D_\lam\htensor R|_{X_\Phi}
\to \CM_y(\lam)^{\oplus r}|_{X_\Phi}$
is an isomorphism.
The statement (c$)_{\Psi}$ is proved.

The proof of (i), (ii)  and (iv) is now completed.

Let us finally show (iii).
By (i) the functor $\tG:\BMD(\lam)\to\TBO(\lam)$ is exact.
Hence the exact sequences
\[
\CM_w(\lam)\to\CL_w(\lam)\to 0,\qquad
0\to\CL_w(\lam)\to\CB_w(\lam)
\]
induce exact sequences
\[
\tG(\CM_w(\lam))\to\tG(\CL_w(\lam))\to 0,\qquad
0\to\tG(\CL_w(\lam))\to\tG(\CB_w(\lam))\,.
\]
Since we have already seen $\tG(\CM_w(\lam))\simeq M(w\circ\lam)$ and
$\tG(\CB_w(\lam))\simeq M^*(w\circ\lam)$, we have only to show
that the morphism
$\tG(\CM_w(\lam))\to\tG(\CB_w(\lam))$ induced by
the canonical morphism $\CM_w(\lam)\to\CB_w(\lam)$ is non-zero.
This follows from (iv).
The statement (iii) is proved.
\qed

\def\TS{\mathop{\tilde{S}}\nolimits}
\def\TV{\mathop{\tilde{V}}\nolimits}
\newsection{Twisted intersection cohomology groups}
\subsection{Combinatorics}
We first recall a result of Lusztig~\cite{L}.

Set
\[
\Gh^*_\BQ=\BQ\otimes_\BZ P,\qquad
\Gamma=\Gh^*_\BQ/P,
\]
where $P$ is as in \S\ref{subsection:flag manifolds}.
Note that the Weyl group $W$ naturally acts on $\Gamma$.
For $\lam\in\Gamma$ let $M^\lam$ be the free $\BZ[q,q^{-1}]$-module with
basis $\{A_w^\lam\}_{w\in W}$.

For $i\in I$ we define $\theta_{i*},
\theta_{i!}\in\Hom_{\BZ[q,q^{-1}]}(M^\lam,M^{s_i\lam})$
by the following.
\eq
&&\theta_{i*}(A_w^\lam)=
\left\{
\begin{array}{ll}
q^{-1}A_{ws_i}^{s_i\lam}\quad
&\mbox{if $\lan\lam,h_i\ran\notin\BZ$, $ws_i>w$},\\
A_{ws_i}^{s_i\lam}\quad
&\mbox{if $\lan\lam,h_i\ran\notin\BZ$,  $ws_i<w$},\\
q^{-1}A_{ws_i}^{\lam}+(q^{-1}-1)A_w^\lam\quad
&\mbox{if $\lan\lam,h_i\ran\in\BZ$, $ws_i>w$},\\
A_{ws_i}^{\lam}\quad
&\mbox{if $\lan\lam,h_i\ran\in\BZ$,  $ws_i<w$},
\end{array}\right.\nn\\
&&\theta_{i!}(A_w^\lam)=
\left\{
\begin{array}{ll}
A_{ws_i}^{s_i\lam}\quad
&\mbox{if $\lan\lam,h_i\ran\notin\BZ$, $ws_i>w$},\\
qA_{ws_i}^{s_i\lam}\quad
&\mbox{if $\lan\lam,h_i\ran\notin\BZ$,  $ws_i<w$},\\
A_{ws_i}^{\lam}\quad
&\mbox{if $\lan\lam,h_i\ran\in\BZ$, $ws_i>w$},\\
qA_{ws_i}^{\lam}+(q-1)A_w^\lam\qquad\quad
&\mbox{if $\lan\lam,h_i\ran\in\BZ$,  $ws_i<w$}.
\end{array}\right.\nn
\endeq

Then $\theta_{i!}:M^\lam\to M^{s_i\lam}$
and $\theta_{i*}:M^{s_i\lam}\to M^{\lam}$ are inverse to each other.

\begin{lemma}[Lusztig~\cite{L}]\label{Lemma:Lusztig1}
\begin{tenumerate}
\item
There exists a unique endomorphism
$m\mapsto\overline{m}$ of the abelian group $M^\lam$ satisfying
\[
\overline{A^\lam_e}=A^\lam_e,\quad
\overline{qm}=q^{-1}\overline{m},\quad
\overline{\theta_{i*}(m)}=\theta_{i!}(\overline{m})
\]
for any $m\in M^\lam$.

\item
We have
$\overline{\overline{m}}=m$ for any $m\in M^\lam$,
and
\[
\overline{A_w^\lam}
\in
q^{-\ell(w)}A_w^\lam
+\sum_{y<w}\BZ[q,q^{-1}]A_y^\lam
\]
for any $w\in W$.
\end{tenumerate}
\end{lemma}
\begin{proposition}[Lusztig~\cite{L}]\label{Prop:Lusztig2}
\begin{tenumerate}
\item
For $w\in W$ and $\lam\in\Gamma$ there exists a unique $C_w^\lam\in M^\lam$
satisfying
\eq
\label{eq1:Prop:Lusztig2}
&&C_w^\lam
\in
A_w^\lam
+
\sum_{y<w}
\Bigl(q^{(\ell(w)-\ell(y)-1)/2}\BZ[q^{-1/2}]\cap\BZ[q,q^{-1}]\Bigr)
A_y^\lam,\\
\label{eq2:Prop:Lusztig2}
&&\overline{C_w^\lam}=q^{-\ell(w)}C_w^\lam.
\endeq
\item
If $w$ is the element of $wW(\lam)$ with minimal length, then for any
$x\in W(\lam)$ we have
\[
C_{wx}^\lam
=
\sum_{y\in W(\lam), y\leq_\lam x}
(-1)^{\ell_\lam(x)-\ell_\lam(y)}
q^{c(y,x)}
P^\lam_{y,x}(q)A^\lam_{wy},
\]
where 
$c(y,x)=\Bigl((\ell(x)-\ell(y))-(\ell_\lam(x)-\ell_\lam(y))\Bigr)/2$, and
$P^\lam_{y,x}(q)\in\BZ[q]$ denotes the Kazhdan-Lusztig polynomial for the
Coxeter group $W(\lam)$ {\rm(}Kazhdan-Lusztig~{\rm \cite{KL1}}{\rm)}.
\end{tenumerate}
\end{proposition}

We define $\kappa\in\Hom_{\BZ[q,q^{-1}]}(M^\lam,M^{-\lam})$ by
$\kappa(A_w^\lam)=A_w^{-\lam}$.

\begin{lemma}
\label{Lemma:Lusztig3}
\begin{tenumerate}
\item
For any $i\in I$ and $\lam\in\Gamma$ we have
$\theta_{i*}\circ\kappa=\kappa\circ\theta_{i*}$ and
$
\theta_{i!}\circ\kappa=\kappa\circ\theta_{i!}$
on $M^\lam$.
\item
We have $\kappa(\overline{m})=\overline{\kappa(m)}$ for any $m\in M^\lam$.
\item
We have $\kappa(C_w^\lam)=C_w^{-\lam}$ for any $w\in W$.
\end{tenumerate}
\end{lemma}
\proof
The statements (i) and (ii) follow 
from the definition of $\theta_{i*}$,
$\theta_{i!}$ and $\bar{}$.

Applying $\kappa$ to (\ref{eq1:Prop:Lusztig2}) we have
\[
\kappa(C_w^\lam)
\in
A_w^{-\lam}
+
\sum_{y<w}
(q^{(\ell(w)-\ell(y)-1)/2}\BZ[q^{-1/2}]\cap\BZ[q,q^{-1}])
A_y^{-\lam}.
\]
By (ii) and (\ref{eq2:Prop:Lusztig2}) we have
\[
\overline{\kappa(C_w^\lam)}
=\kappa(\overline{C_w^\lam})
=\kappa(q^{-\ell(w)}C_w^\lam)
=q^{-\ell(w)}\kappa(C_w^\lam).
\]
Thus we obtain (iii) by Proposition~\ref{Prop:Lusztig2}.
\qed

Lusztig~\cite{L} used Proposition~\ref{Prop:Lusztig2} to compute the twisted
intersection cohomology groups of the finite-dimensional Schubert varieties.
In order to compute that of the finite-{\it codimensional} Schubert
varieties, we need its dual version.

Set
\[
N^\lam=\Hom_{\BZ[q,q^{-1}]}(M^\lam,\BZ[q,q^{-1}]),
\]
and define $B^\lam_w\in N^\lam$ for $w\in W$ by
$\lan B_w^\lam,A_y^\lam\ran=\delta_{y,w}$.
Then any element of $N^\lam$ is uniquely written as a formal infinite sum
$\sum_{w\in W}a_wB_w^\lam$ with $a_w\in\BZ[q,q^{-1}]$.

For $i\in I$ we define $\theta_{i*},
\theta_{i!}\in\Hom_{\BZ[q,q^{-1}]}(N^\lam,N^{s_i\lam})$
by
\eq
\lan\theta_{i*}n,m\ran=\lan n,\theta_{i*}m\ran,\quad
\lan\theta_{i!}n,m\ran=\lan n,\theta_{i!}m\ran\qquad
\mbox{
for $n\in N^\lam, m\in M^\lam$}.\nn
\endeq
By the definition we have the following.
\eq
&&\theta_{i*}(B_w^\lam)=
\left\{
\begin{array}{ll}
B_{ws_i}^{s_i\lam}\quad
&\mbox{if $\lan\lam,h_i\ran\notin\BZ$, $ws_i>w$},\\
q^{-1}B_{ws_i}^{s_i\lam}\quad
&\mbox{if $\lan\lam,h_i\ran\notin\BZ$,  $ws_i<w$},\\
B_{ws_i}^{\lam}\quad
&\mbox{if $\lan\lam,h_i\ran\in\BZ$, $ws_i>w$},\\
q^{-1}B_{ws_i}^{\lam}+(q^{-1}-1)B_w^\lam\quad
&\mbox{if $\lan\lam,h_i\ran\in\BZ$,  $ws_i<w$},
\end{array}\right.\nn\\
&&\theta_{i!}(B_w^\lam)=
\left\{
\begin{array}{ll}
qB_{ws_i}^{s_i\lam}\quad
&\mbox{if $\lan\lam,h_i\ran\notin\BZ$, $ws_i>w$},\\
B_{ws_i}^{s_i\lam}\quad
&\mbox{if $\lan\lam,h_i\ran\notin\BZ$,  $ws_i<w$},\\
qB_{ws_i}^{\lam}+(q-1)B_w^\lam\qquad\quad
&\mbox{if $\lan\lam,h_i\ran\in\BZ$, $ws_i>w$},\\
B_{ws_i}^{\lam}\quad
&\mbox{if $\lan\lam,h_i\ran\in\BZ$,  $ws_i<w$}.
\end{array}\right.\nn
\endeq
Let $a\mapsto\overline{a}$ be the endomorphism of the ring $\BZ[q,q^{-1}]$
given by $\overline{q}=q^{-1}$.
We define an endomorphism $n\mapsto\overline{n}$ of the abelian group
$N^\lam$ by
\[
\lan\overline{n},m\ran=\overline{\lan n,\overline{m}\ran}\qquad
\mbox{ for $n\in N^\lam, m\in M^\lam$}.
\]
By Lemma~\ref{Lemma:Lusztig1} we have $\overline{\overline{n}}=n$ for any
$n\in N^\lam$, and
\[
\overline{B_w^\lam}
\in
q^{\ell(w)}B_w^\lam
+\sum_{y>w}\BZ[q,q^{-1}]B_y^\lam
\]
for any $w\in W$.
Define $D_w^\lam\in N^\lam$ by $\lan D_w^\lam,C_y^\lam\ran=\delta_{y,w}$.
By Proposition~\ref{Prop:Lusztig2} we have the following.
\begin{proposition}\label{Prop:Dw1}
\begin{tenumerate}
\item
For $w\in W$ and $\lam\in\Gamma$ we have
\eq
&&D_w^\lam
\in
B_w^\lam
+
\sum_{y>w}
(q^{(\ell(y)-\ell(w)-1)/2}\BZ[q^{-1/2}]\cap\BZ[q,q^{-1}])
B_y^\lam,\label{eq:Dw:unique1}\\
&&\overline{D_w^\lam}=q^{\ell(w)}D_w^\lam.\label{eq:Dw:unique2}
\endeq
\item
If $w$ is the element of $wW(\lam)$ with minimal length, then for any
$x\in W(\lam)$ we have
\[
D_{wx}^\lam
=
\sum_{y\in W(\lam), y\geq_\lam x}
q^{c(x,y)}
Q^\lam_{x,y}(q)B^\lam_{wy},
\]
where 
$c(x,y)=\Bigl((\ell(y)-\ell(x))-(\ell_\lam(y)-\ell_\lam(x))\Bigr)/2$
and
$Q^\lam_{x,y}(q)\in\BZ[q]$ 
denotes the inverse Kazhdan-Lusztig polynomial for
the Coxeter group $W(\lam)$ given by
\eq
\sum_{x\leq_\lam y\leq_\lam z}
(-1)^{\ell_\lam(y)-\ell_\lam(x)}
Q^\lam_{x,y}(q)P^\lam_{y,z}(q)
=\delta_{x,z}.
\label{eq:inverse KL}
\endeq
\end{tenumerate}
\end{proposition}

Define $\kappa\in\Hom_{\BZ[q,q^{-1}]}(N^\lam, N^{-\lam})$ by
\[
\lan\kappa(n),m\ran=\lan n,\kappa(m)\ran\qquad
\mbox{ for $n\in N^\lam, m\in M^{-\lam}$}.
\]

By Lemma~\ref{Lemma:Lusztig3} we obtain the following.
\begin{lemma}
\label{Lemma:Dw:2}
\begin{tenumerate}
\item
We have $\kappa(B_w^\lam)=B_w^{-\lam}$ for any $w\in W$.
\item
For any $i\in I$ and $\lam\in\Gamma$ we have
$\theta_{i*}\circ\kappa=\kappa\circ\theta_{i*}$ and
$
\theta_{i!}\circ\kappa=\kappa\circ\theta_{i!}$ on $N^\lam$.
\item
We have $\kappa(\overline{n})=\overline{\kappa(n)}$ for any $n\in N^\lam$.
\item
We have $\kappa(D_w^\lam)=D_w^{-\lam}$ for any $w\in W$.
\end{tenumerate}
\end{lemma}

Let $R$ be a ring containing $\BZ[q,q^{-1}]$ as a subring.
Assume that we are given an involutive automorphism $r\mapsto\overline{r}$
of the ring $R$ and a family of $\BZ$-submodules $\{R_i\}_{i\in\BZ}$ of $R$
such that
\eq
R=\bigoplus_{i\in\BZ}R_i,\quad
R_iR_j\subset R_{i+j},\quad
1\in R_0,\quad
q\in R_2,\quad
\overline{q}=q^{-1},\quad
\overline{R_i}=R_{-i}.
\label{eq:property of R}
\endeq
Define the endomorphism $m\mapsto \ol m$ of the abelian group
$R\otimes_{\BZ[q,q^{-1}]}M^\lam$ by
\[\ol{a\otimes m}=\ol a\otimes \ol m\quad 
\mbox{for $a\in R$ and $m\in M^\lam$.}
\]
Set $N_R^\lam=\Hom_R(R\otimes_{\BZ[q,q^{-1}]}M^\lam,R)$.
We can naturally regard $N^\lam$ as a $\BZ[q,q^{-1}]$-submodule of $N_R^\lam$.
Define an endomorphism $n\mapsto\overline{n}$ of the abelian group
$N_R^\lam$ by
\[
\lan\overline{n},m\ran=\overline{\lan n,\overline{m}\ran}\qquad
\mbox{ for $n\in N_R^\lam, m\in R\otimes_{\BZ[q,q^{-1}]}M^\lam$},
\]
and a homomorphism $\kappa:N_R^\lam\to N_R^{-\lam}$ of $R$-modules by
\[
\lan\kappa(n),m\ran=\lan n,\kappa(m)\ran\qquad
\mbox{ for $n\in N_R^\lam, m\in R\otimes_{\BZ[q,q^{-1}]}M^\lam$}.
\]

Then we have the following characterization of $D_w^\lam$.
\begin{proposition}
\label{prop:characterization of D_w}
Let $w\in W$ and $\lam\in\Gamma$.
Assume that $D^+\in N_R^\lam$ and $D^-\in N_R^{-\lam}$ satisfy the following
properties:
\eq
&&D^\pm
\in
B_w^{\pm\lam}
+
\sum_{y>w}
(\bigoplus_{i\leq\ell(y)-\ell(w)-1}R_i)
B_y^{\pm\lam},
\label{eq:Dw:unique1s}
\\
&&\overline{\kappa(D^+)}=q^{\ell(w)}D^-.
\label{eq:Dw:unique2s}
\endeq
Then we have $D^\pm=D_w^{\pm\lam}$.
\end{proposition}
\proof
Since (\ref{eq:Dw:unique1s}) and (\ref{eq:Dw:unique2s}) are satisfied for
$D^\pm=D^{\pm\lam}_w$, it is sufficient to show that there exist unique
$D^\pm\in N_R^\lam$ satisfying (\ref{eq:Dw:unique1s}) and
(\ref{eq:Dw:unique2s}).

By (\ref{eq:Dw:unique1s}) we have
\[
D^\pm=\sum_{y\geq w}F^\pm_yB_y^{\pm\lam}\qquad
\mbox{
with $F^\pm_w=1$ and $F^\pm_y\in\bigoplus_{i\leq\ell(y)-\ell(w)-1}R_i$ for
$y>w$.
}
\]
We have to show that $F_y^\pm$ are uniquely determined by the condition
(\ref{eq:Dw:unique2s}).
Write
\[
\overline{B_y^\lam}=\sum_{z\geq y}G_{y,z}B_z^\lam\qquad
\mbox{
with $G_{y,z}\in R, G_{y,y}=q^{\ell(w)}$
}.
\]
Then we have
\[
q^{-\ell(w)}\overline{\kappa(D^+)}=
q^{-\ell(w)}
\sum_{y\geq w}\overline{F^+_y}(\sum_{z\geq y}\overline{G_{y,z}}B_z^{-\lam})
=\sum_{z\geq w}
(\sum_{z\geq y\geq w}
q^{-\ell(w)}\overline{F^+_yG_{y,z}})B_z^{-\lam},
\]
and hence
\[
\sum_{z\geq y\geq w}
q^{-\ell(w)}\overline{F^+_yG_{y,z}}
=F^-_z
\]
for any $z\geq w$.
Thus
\[
F^-_z-q^{-\ell(w)+\ell(z)}\overline{F^+_z}
=\sum_{z> y\geq w}
q^{-\ell(w)}\overline{F^+_yG_{y,z}}
\]
for any $z>w$.
By the assumption we have
\[
F^-_z
\in
\bigoplus_{i\leq\ell(z)-\ell(w)-1}
R_i,
\qquad
q^{-\ell(w)+\ell(z)}\overline{F^+_z}
\in
\bigoplus_{i\geq\ell(z)-\ell(w)+1}
R_i.
\]
Therefore $F^\pm_z$ are uniquely determined inductively.
\qed

\subsection{Character formula}
For $\lam\in\Gh^*$ and a finite addmissible subset of $W$, 
let $K(\BMD_\Phi(\lam))$ be the Grothendieck group of the category
$\BMD_\Phi(\lam)$. It is a module with
$\Bigl\{[\CM_w(\lam)]\Bigr\}_{w\in\Phi}$ as a basis.
Let $K(\BMD(\lam))$ be the projective limit of
$K(\BMD_\Phi(\lam))$, where $\Phi$ ranges over the set of
finite admissible subsets of $W$.
Then $\bigl\{[\CM_w(\lam)]\bigr\}_{w\in W}$ as well as
$\bigl\{[\CL_w(\lam)]\bigr\}_{w\in W}$ is a formal basis of
$K(\BMD(\lam))$.

The aim of this section is to prove the following result.
\begin{theorem}\label{thm:perverse:mult}
Let $\lam\in\Gh_\BQ^*$, and let $w\in W$ such that
$\ell(z)>\ell(w)$ for any $z\in wW(\lam)\setminus\{w\}$.
Then for any $x\in W(\lam)$ we have
\begin{eqnarray}
&&[\CL_{wx}(\lam)]
=
\sum_{y\geq_\lam x}
(-1)^{l(y)-l(x)}Q^\lam_{x,y}(1)
[\CM_{wy}(\lam)],\label{eq:perverse:mult1}
\\
&&[\CM_{wx}(\lam)]
=
\sum_{y\geq_\lam x}
P^\lam_{x,y}(1)
[\CL_{wy}(\lam)].
\label{eq:perverse:mult2}
\end{eqnarray}
\end{theorem}

The proof of this theorem will be given in the next subsection.
The corresponding result for finite-dimensional Schubert varieties was
proved by Lusztig (see \cite{L0}, \cite{L}).

Note that (\ref{eq:perverse:mult1}) and (\ref{eq:perverse:mult2}) are
equivalent by (\ref{eq:inverse KL}).

By Theorem~\ref{Main theorem}, Proposition~\ref{D:mod:prop:tH1} and
Theorem~\ref{thm:perverse:mult}, 
we obtain the following main result of this paper.

\begin{theorem}\label{thm:MAIN THEOREM}
Assume that $\lam\in\Gh^*$ satisfies the following conditions.
\eq
&&\mbox{
$2(\alpha,\lam+\rho)\ne(\alpha,\alpha)$ 
for any positive imaginary root $\alpha$.
}\\
&&\mbox{
$(\alpha^\vee,\lam+\rho)\notin\BZ_{\leq0}$ 
for any positive real root $\alpha$.
}\\
&&\mbox{
If $w\in W$ satisfies $w\circ\lam=\lam$, then $w=1$.
}
\label{main:iso}\\
&&\mbox{
$(\alpha^\vee,\lam)\in\BQ$ for any real root $\alpha$.
}
\endeq
Then for any $w\in W(\lam)$ we have 
\eq
&&\label{eq:ch1}
\ch(M(w\circ\lam))
=\sum_{y\ge_\lam w}
P^\lam_{w,y}(1)
\ch(L(y\circ\lam)),
\\
&&
\ch(L(w\circ\lam))
=\sum_{y\ge_\lam w}
(-1)^{\ell_\lam(y)-\ell_\lam(w)}
Q^\lam_{w,y}(1)
\ch(M(y\circ\lam)).\label{eq:ch2}
\endeq
\end{theorem}

As a special case, we obtain the following result.

\begin{theorem}\label{thm:MAIN THEOREM:affine}
Assume that $\Gg$ is finite-dimensional or affine and
$\lam\in\Gh^*$ satisfies
\eq
&&(\beta^\vee,\lam+\rho)\in\BQ\setminus\BZ_{\le0} 
\quad\mbox{for any $\beta\in\Delta_\re^+$.}\\
&&(\delta,\lam+\rho)\not=0\quad
\mbox{if $\Gg$ is affine. Here $\delta$ is an imaginary root.}
\label{root:}
\endeq
Then $(\ref{eq:ch1})$ and $(\ref{eq:ch2})$ hold
for any $w\in W(\lam)$.
\end{theorem}
In the affine case, the condition (\ref{main:iso})
on the triviality of the isotropy subgroup of $\lam$
follows from the following well-known lemma.
\Lemma
The isotropy subgroup 
$\{w\in W;w\circ\lam=\lam\}$ is generated by 
$\{s_\beta;\beta\in\Delta^+_\re,\ (\beta,\lam+\rho)=0\}$
whenever $\Gg$ is affine and $\lam\in\Gh^*$ satisfies $(\ref{root:})$.
\enlemma

In the affine case,
we can derive the following result on the non-regular highest weight case
from the regular highest weight case above  
by using the translation functors (we omit the proof).

\Theorem
Let $\Gg$ be an affine Lie algebra, and assume that
$\lam\in\Gh^*$ satisfies
\eq
&&\mbox{
$(\delta,\lam+\rho)\not=0$, 
}\\
&&\mbox{
$(\alpha^\vee,\lam+\rho)\in\BQ\setminus\BZ_{<0}$ 
for any positive real root $\alpha$.
}
\endeq
Then $W_0(\lam)=\{w\in W;w\circ\lam=\lam\}$ is a finite group.
Let $w$ be an element of $W(\lam)$ which is
the longest element of $wW_0(\lam)$.
Then we have
\eqn
\ch(L(w\circ\lam))
=\sum_{y\ge_\lam w}
(-1)^{\ell_\lam(y)-\ell_\lam(w)}
Q^\lam_{w,y}(1)
\ch(M(y\circ\lam)).
\endeqn
\end{theorem}
\subsection{Hodge modules on flag manifolds}
\label{subsection:Hodge}
Let $R$ (resp.\ $R_i$ for $i\in\BZ$) denote the Grothendick group of the
category of mixed Hodge structures (resp. pure Hodge structures with weight
$i$) over $\QC$.
Then we have $R=\bigoplus_{i\in\BZ}R_i.$
Let $\QC^H(k)$ be the Hodge structure of Tate with weight $-2k$.
Set $q=[\QC^H(-1)]\in R_2$, and let $r\mapsto\overline{r}$ denote the
endomorphism of the ring $R$ induced by the duality operation on the mixed
Hodge structures.
Then the condition (\ref{eq:property of R}) is satisfied for the above $R$.

For a smooth $\BC$-scheme $S$, 
let $\MH(S)$ denote the category of mixed Hodge modules on $S$ 
(see Saito\cite{Saito}).
Here we use the convention that the perversity
is stabel under the smooth inverse image.

For a scheme $S$ satisfying (S)
with a smooth projective system $\{S_n\}_{n\in\BN}$,
let us denote by $\MH(S)$ the inductive limit of
$\MH(S_n)$. It is an abelian category
and there is an exact functor
\[\MH(S)\to \DM_h(\CD_S)\]
We call an object of $\MH(S)$ a mixed Hodge module over $S$.

For $\lam\in\Gamma=\Gh_\BQ^*/P$ we denote by $T^H(\lam)$ the Hodge module on
$H$ corresponding to $\CT(\lam)$ (see \S \ref{sec:equiv}) of weight $0$.
By the assumption on $\lam$, the monodromies of the corresponding local system
are roots of unity, and hence it has a structure of 
variation of polarizable Hodge structure. 
Hence $T^H(\lam)$ is defined as a Hodge module on $H$ of weight $0$.

For a $\BC$-scheme $S$ satisfying (S) with an action of $H$,
we can define the twisted $H$-equivariance of mixed Hodge module on $S$
as in \S \ref{sec:equiv} by the aid of $T^H(\lam)$.
We denote by $\MH(S,\lam)$ the category
of twisted $H$-equivariant mixed Hodge modules on $S$ with twist $\lam$.
It depends only on the image of $\lam$ in $\Gamma=\Gh_\BQ^*/P$.

Recall that $\TX=G/N^-$ and $\xi:\TX\to X$ is the natural projection.
Then $B\times H$ acts on $\TX$ by
$(b,h)\circ(gN^-)=bgh^{-1}N^-$.
By the action of $H$ on $\TX$,
$\xi:\TX\to X$ is a principal $H$-bundle.
For a finite admissible subset $\Phi$ of $W$ and $\lam\in\Gamma$ we denote
by $\MH_\Phi(\lam)$ the category of twisted $(N\times H)$-equivariant
mixed Hodge modules with twist $\lam$.

Set
\[
\MH(\lam)=\limp_{\Phi}\MH_\Phi(\lam).
\]

Since $\xi^\sharp D_{X,\lam}=D_\TX$  and since $\xi$ is a smooth morphism,
we have an equivalence $\xi^\bullet$ from the category of holonomic
$D_{\lam}$-modules to the category of $H$-equivariant
holonomic $D_{\TX}$-modules with twist $\lam$.
Hence we have an exact functor
\eq\label{hodg:D}
&&\MH(\lam)\to\BMD(\lam).
\endeq

For $w\in W$,
set $\TX_w=BwN^-/N^-=\xi^{-1}(X_w)$ and let $\Ti:\TX_w\hookrightarrow\TX$
be the embedding.
By the isomorphism of schemes
\[
H\times N(\Delta^+\cap w\Delta^+)\iso\TX_w\qquad
((h,u)\mapsto uwh^{-1}N^-),
\]
we can define a morphism $p_w:\TX_w\to H$ by $p_w(uwh^{-1}N^-)=h$ for $u\in
N(\Delta^+\cap w\Delta^+)$ and $h\in H$.
We define a Hodge module $F^H_w(\lam)$ on $\TX_w$ 
by $F^H_w(\lam)=p_w^*T(\lam)$.
We denote by the same letter $F^H_w(\lam)$ the object 
$\Ti_!F^H_w(\lam)$ in the derived category of
the category of mixed Hodge modules.
Let us denote by ${}^\pi F^H_w(\lam)[-\ell(w)]$
the minimal extension of $F^H_w(\lam)[-\ell(w)]$.
Then $F^H_w(\lam)[-\ell(w)]$ and ${}^\pi F^H_w(\lam)[-\ell(w)]$ 
are objects of $\MH(\lam)$.
By the functor (\ref{hodg:D}),
$F_w(\lam)[-\ell(w)]$ and ${}^\pi F_w(\lam)[-\ell(w)]$
correspond to the objects $\CM_w(\lam)$ and $\CL_w(\lam)$
of $\BMD(\lam)$.

\medskip
For a finite admissible set $\Phi$,
$\TX_\Phi$ has the $N\times H$-orbit
decomposition $\TX_\Phi=\bigsqcup_{w\in \Phi}\TX_w$.
Hence the irreducible objects of $\MH_\Phi(\lam)$ is
of the form $H\otimes{}^\pi F_w^H(\lam)[-\ell(w)]|\TX_\Phi$ 
for some $w\in\Phi$ and some irreducible Hodge structure $H$.
We denote the Grothendieck group of $\MH_\Phi(\lam)$
by $K(\MH_\Phi(\lam))$.
This has a structure of $R$-module.
We set
\[
K(\MH(\lam))=\limp_{\Phi}K(\MH_\Phi(\lam)).
\]
For $F\in\MH(\lam)$ we denote by $[F]$ the element of $K(\MH(\lam))$
corresponding to $F$, and $\Big[F[n]\Big]=(-1)^n[F]$.
Any $m\in K(\MH(\lam))$ can be written uniquely as
\[
m=\sum_{w\in W}a_w[{}^\pi F^H_w(\lam)]=\sum_{w\in
W}b_w[F_w^H(\lam)]\qquad(a_w, b_w\in R).
\]
Define an isomorphism
\[
\varphi_\lam:K(\MH(\lam))\iso N_R^\lam
\]
of $R$-modules by $\varphi_\lam([F_w^H(\lam)])=B^\lam_w$.

\medskip
For $i\in I$ we shall define
\[
\TS_{i*},\;\TS_{i!}\in\Hom_R\Bigl(K(\MH(\lam)),K(\MH(s_i\lam))\Bigr).
\]
The definition is analogous to \S \ref{Radon},
and we use the notations in \S \ref{Radon}.

Set $N^-_i=\exp(\Gn(\Delta^-\setminus\{-\alpha_i\}))\subset N^-$,
and $\TZ_0=G/N^-_i$.
The group $B\times H$ acts on $\TZ_0$ by $(b,h)\circ(gN^-_i)=gh^{-1}N^-_i$.
Let $\tp_i:\TZ_0\to \TX$ ($i=1,2$)
be the morphism defined by
\eqn
&&p_1(gN^-_i)=gN^-\quad\mbox{and}\quad
p_2(gN^-_i)=gs_iN^-.
\endeqn
We have the commutative diagram
\eq\label{daigram1}
\begin{array}{ccccc}
\TX&\mapleftu{\tp_1}&\TZ_0&\maprightu{\tp_2}&\TX\\
\mapdownl{\xi}&&\mapdownl{\xi'}&&\mapdownl{\xi}\\
X&\mapleftu{p_1}&Z_0&\maprightu{p_2}&X.
\end{array}
\endeq
Then $\xi':\TZ_0\to Z$ is a principal $H$-bundle.
The morphisms $\tp_1$ and $\tp_2$ are $B$-equivariant, 
and they satisfy the following
relation with the action of $H$.
\eq\label{ac:H}
&
\ba{l}
\tp_1(hz)=h\tp_1(z)\\[2pt]
\tp_2(hz)=s_i(h)\tp_2(z)
\ea
&\quad
\mbox{for $h\in H$ and $z\in\TZ_0$.}
\endeq
Here $s_i$ is the group automorphism of $H$
corresponding to the simple reflection $s_i\in\Aut(\Gh)$.

For $F\in\MH(\TX,\lam)$ we set
\eq
\TS_{i!}(F)=\BR\tp_2{}_{!}\tp_1{}^*F,\qquad
\TS_{i*}(F)=\BR\tp_2{}_{*}\tp_1{}^!F.
\endeq
Then $H^k(\TS_{i!}(F))$ and $H^k(\TS_{i!}(F))$ are objects of
$\MH(s_i\lam)$ by (\ref{ac:H}).

We define
$\TS_{i*}, \; \TS_{i!}\in\Hom_R\Bigl(K(\MH(\lam)),K(\MH(s_i\lam))\Bigr)$
by
\[
\TS_{i!}([F])=\sum_{k}(-1)^k[H^k(\TS_{i!}F)],\qquad
\TS_{i*}([F])=\sum_{k}(-1)^k[H^k(\TS_{i*}F)].
\]

\begin{proposition}\label{prop:S=theta}
We have
\[
\varphi_{s_i\lam}\circ \TS_{i*}=\theta_{i*}\circ\varphi_\lam,\qquad
\varphi_{s_i\lam}\circ \TS_{i!}=\theta_{i!}\circ\varphi_\lam.
\]
\end{proposition}
\proof
Fix $w\in W$ such that $ws_i>w$.
It is sufficient to show the following:
\eq
\TS_{i*}[F_{w}^H(\lam)]&=&
\left\{
\begin{array}{ll}
[F_{ws_i}^H(s_i\lam)]\qquad\qquad \qquad \qquad \qquad\quad
&\mbox{if $\lan\lam,h_i\ran\notin\BZ$,}\\
{}[F_{ws_i}^H(\lam)]\quad
&\mbox{if $\lan\lam,h_i\ran\in\BZ$},
\end{array}\right.\label{eq:S:1}\\
\TS_{i*}[F_{ws_i}^H(\lam)]&=&
\left\{
\begin{array}{ll}
q^{-1}[F_{w}^H(s_i\lam)]\quad
&\mbox{if $\lan\lam,h_i\ran\notin\BZ$,}\\
q^{-1}[F_{w}^H(\lam)]+(q^{-1}-1)[F_{ws_i}^H(\lam)]\quad
&\mbox{if $\lan\lam,h_i\ran\in\BZ$},
\end{array}
\right.\label{eq:S:2}\\
\TS_{i!}[F_{w}^H(\lam)]&=&
\left\{
\begin{array}{ll}
q[F_{ws_i}^H(s_i\lam)]\quad
&\mbox{if $\lan\lam,h_i\ran\notin\BZ$},\\
q[F_{ws_i}^H(\lam)]+(q-1)[F_{w}^H(\lam)]\qquad\quad
&\mbox{if $\lan\lam,h_i\ran\in\BZ$}.
\end{array}\right.\label{eq:S:3}\\
\TS_{i!}[F_{ws_i}^H(\lam)]&=&
\left\{
\begin{array}{ll}
[F_{w}^H(s_i\lam)]\qquad\qquad \qquad \qquad \quad\quad
&\mbox{if $\lan\lam,h_i\ran\notin\BZ$},\\
{}[F_{w}^H(\lam)]\quad
&\mbox{if $\lan\lam,h_i\ran\in\BZ$}.
\end{array}\right.\label{eq:S:4}
\endeq

Set $Y=\TX_w\sqcup \TX_{ws_i}$ and
let $j:Y\to \TX$ be the embedding.

As in Lemma 1.4.1 and Corollary 1.5.2 in \cite{KTneg2},
$\TS_{i!}$  and $\TS_{i*}$ commute with $j_!$.
Hence we can reduce these statements to the case 
$\Gg=\Gsl_2$ where we can check them directly.
Details are omitted.
\qed

The duality functor for Hodge modules induces a contravariant exact functor
\[
{\bf D}:\MH(\lam)\to\MH(-\lam).
\]
We also denote by
\[
{\bf D}:K(\MH(\lam))\to K(\MH(-\lam)).
\]
the induced homomorphism of abelian groups.
By the definition we have the following result.
\begin{lemma}
\label{lemma:duality1}
\begin{tenumerate}
\item
We have
\[
{\bf D}(rn)=\overline{r}{\bf D}(n)
\]
for any $r\in R$ and $n\in K(\MH(\lam))$.
\item
We have
\[
{\bf D}\circ\TS_{i*}=\TS_{i!}\circ{\bf D},\qquad
{\bf D}\circ\TS_{i!}=\TS_{i*}\circ{\bf D}
\]
on $K(\MH(\lam))$.
\end{tenumerate}
\end{lemma}

\begin{proposition}
\label{Prop:duality2}
We have
\[
\varphi_{-\lam}({\bf D}n)
=\kappa(\overline{\varphi_\lam(n)})
\]
for any $n\in K(\MH(\lam))$.
\end{proposition}
\proof
It is sufficient to show
\[
\lan\varphi_{-\lam}({\bf D}n),A_w^{-\lam}\ran=
\lan\kappa(\overline{\varphi_\lam(n)}),A_w^{-\lam}\ran
\]
for any $w\in W$.
By Lemma~\ref{Lemma:Dw:2} the right side coincides with
$\overline{
\lan
\varphi_\lam(n)
,
\overline{A_w^\lam}
\ran
}$
and hence we have to show
\eq
\label{eq:Prop:duality}
\lan\varphi_{-\lam}({\bf D}n),A_w^{-\lam}\ran=
\overline{\lan\varphi_\lam(n),\overline{A_w^\lam}\ran}.
\endeq

We first consider the case $w=e$.
In this case we have $\overline{A_e^\lam}=A_e^\lam$.
We may assume that $n=[F^H_x(\lam)]$ for $x\in W$.
Since ${\bf D}F^H_x(\lam)\cong F^H_x(\lam)(-\ell(x))[-2\ell(x)]$
on a neighborhood of $\TX_w$, we have
\[
\varphi_{-\lam}({\bf D}[F^H_x(\lam)])\in
q^{\ell(x)}B_x^{-\lam}+\sum_{y>x}RB_y^{-\lam}.
\]
Thus the both sides of (\ref{eq:Prop:duality}) are equal to $\delta_{x,e}$.

For general $w\in W$, take a reduced expression $w=s_{i_1}\cdots s_{i_r}$.
By the definition of $\theta_{i!}$ we have
$A_w^{\pm\lam}=\theta_{i_r!}\cdots\theta_{i_1!}A_e^{\pm w\lam}$.
Thus we have
\[
\overline{A_w^{\lam}}=
\overline{\theta_{i_r!}\cdots\theta_{i_1!}A_e^{w\lam}}=
\theta_{i_r*}\cdots\theta_{i_1*}\overline{A_e^{w\lam}}=
\theta_{i_r*}\cdots\theta_{i_1*}A_e^{w\lam}.
\]
Hence by Lemma~\ref{lemma:duality1} we obtain
\eq
\lan\varphi_{-\lam}({\bf D}n),A_{w}^{-\lam}\ran
&=&
\lan\varphi_{-\lam}({\bf
D}n),\theta_{i_r!}\cdots\theta_{i_1!}A_e^{-w\lam}\ran\nn\\
&=&
\lan\theta_{i_1!}\cdots\theta_{i_r!}\varphi_{-\lam}({\bf
D}n),A_e^{-w\lam}\ran\nn\\
&=&
\lan\varphi_{-w\lam}(\TS_{i_1!}\cdots \TS_{i_r!}{\bf
D}n),A_e^{-w\lam}\ran,\nn\\
&=&
\lan\varphi_{-w\lam}(
{\bf D}
\TS_{i_1*}\cdots \TS_{i_r*}n)
,A_e^{-w\lam}\ran,\nn\\
&=&
\overline{
\lan\varphi_{w\lam}(\TS_{i_1*}\cdots \TS_{i_r*}n),A_e^{w\lam}\ran
}
\nn\\
&=&
\overline{
\lan\theta_{i_1*}\cdots\theta_{i_r*}\varphi_\lam(n),A_e^{w\lam}\ran
}
\nn\\
&=&
\overline{\lan\varphi_{\lam}(n),\theta_{i_r*}\cdots\theta_{i_1*}A_e^{w\lam}
\ran}\nn\\
&=&
\overline{\lan\varphi_\lam(n),\overline{A_w^\lam}\ran}
\nn
\endeq
\qed

\begin{theorem}\label{thm:Hodge:mult}
We have $\varphi_\lam([{}^\pi F_w^H(\lam)])=D^\lam_w$ for any $w\in W$.
\end{theorem}

Note that Theorem~\ref{thm:perverse:mult} is a consequence of
Theorem~\ref{thm:Hodge:mult}.
In fact, if $w\in W$ satisfies $\ell(z)>\ell(w)$ for any $z\in
wW(\lam)\setminus\{w\}$, then we have
\eq\label{eq:Hodge:mult}
[{}^\pi F_{wx}^H(\lam)]
=
\sum_{y\in W(\lam), y\geq_\lam x}
q^{c(x,y)}
Q^\lam_{x,y}(q)[F_{wy}^H(\lam)]
\endeq
in $K(\MH(\lam))$ for any $x\in W(\lam)$ by Proposition~\ref{Prop:Dw1} and
Theorem~\ref{thm:Hodge:mult}.
Applying the canonical homomorphism $K(\MH(\lam))\to K(\BMD(\lam))$ to
(\ref{eq:Hodge:mult}) we obtain (\ref{eq:perverse:mult1}).

\medskip
\noindent
{\sl Proof of Theorem~\ref{thm:Hodge:mult}}.\quad
Set $D^\pm=\varphi_{\pm\lam}([{}^\pi F_w^H(\pm\lam)])$.
It is sufficient to show that $D^\pm$ satisfy the conditions
(\ref{eq:Dw:unique1s}) and (\ref{eq:Dw:unique2s}) in
Proposition~\ref{prop:characterization of D_w}.

By Proposition~\ref{Prop:duality2} we have
\eq
\overline{\kappa(D^+)}
&=&
\overline{\kappa\varphi_\lam([{}^\pi F_w^H(\lam)])}\nn\\
&=&
\varphi_{-\lam}({\bf D}([{}^\pi F_w^H(\lam)]))\nn\\
&=&
\varphi_{-\lam}([{\bf D}{}^\pi F_w^H(\lam)])\nn\\
&=&
\varphi_{-\lam}(\Bigl[{}^\pi F_w^H(-\lam)[-2\ell(w)](-\ell(w))\Bigr])\nn\\
&=&
q^{\ell(w)}\varphi_{-\lam}([{}^\pi F_w^H(-\lam)])\nn\\
&=&
q^{\ell(w)}D^-\nn,
\endeq
and hence (\ref{eq:Dw:unique2s}) holds.

For $x\in W$, let $\Ti_x:\TX_x\hookrightarrow \TX$
be the embedding and
let $\MH^x(\pm\lam)$ denote the abelian category of 
twisted $N\times H$-equivariant mixed Hodge
modules on $\TX_x$ with twist $\pm\lam$.
Since $\TX_w$ is an orbit of
$N\times H$, any object of  $\MH^x(\pm\lam)$ has 
a form $H\otimes F_x^H(\pm\lam)$ for some mixed Hodge module $H$.
Hence its Grothendieck group $K(\MH^x(\pm\lam))$ 
is a free $R$-module generated by $[F_x^H(\pm\lam)]$.
The inverse image functor $\Ti_x^*$ induces a homomorphism
\[
\iota_x:K(\MH(\pm\lam))\to R\qquad
([\Ti_x^*F]=\iota_x([F])[F_x^H(\pm\lam)]).
\]
of $R$-modules.
Then we have
\[
\varphi_{\pm\lam}(n)=\sum_{x\in W}\iota_x(n)B_w^{\pm\lam}
\]
for any $n\in\MH(\pm\lam)$ because this formula obviously holds for
$n=[F_y^H(\pm\lam)]$ with $y\in W$.
We have obviously $\iota_w([([{}^\pi F_w^H(\pm\lam)])=1$.
Let $y>w$.
Since ${}^\pi F_w^H(\pm\lam)$ is pure of weight 0,
$H^j(\Ti_y^*({}^\pi F_w^H(\pm\lam)))$ is a mixed Hodge module of weight
$\leq j$ for any $j$.
On the other hand the perversity property of ${}^\pi F_w^H(\pm\lam)$ implies
$H^j(\Ti_y^*({}^\pi F_w^H(\pm\lam)))=0$ for $j\geq\ell(y)-\ell(w)$.
Thus we obtain $\iota_y([([{}^\pi
F_w^H(\pm\lam)])\in\sum_{j\leq\ell(y)-\ell(w)-1}R_j$.
Hence the condition (\ref{eq:Dw:unique2s}) also holds.
\qed
By using a $\BC^\times$-action, we can prove that, for any $j$,
$H^j(\Ti_y^*({}^\pi F_w^H(\pm\lam)))$ is a pure Hodge module of 
weight $j$ as in
Kazhdan-Lusztig~\cite{KL1} and Kashiwara-Tanisaki~\cite{KTpos2}.
This gives the following stronger version of 
Theorem~\ref{thm:perverse:mult}.
Since this result is not used in this paper, the details are omitted.

\begin{theorem}\label{thm:perverse:mult-strong}
Let $\lam\in\Gamma=\Gh_\BQ^*/P$, and let $w\in W$ such that
$\ell(z)>\ell(w)$ for any $z\in wW(\lam)\setminus\{w\}$.
Let $x, y\in W(\lam)$ such that $y\geq x$ and write
$Q^\lam_{x,y}(q)=\sum_jc_jq^j$.
\begin{tenumerate}
\item
$H^{2j+1}(\Ti_y^*({}^\pi F^H_{wx}(\lam))=0$ for any $j\in\BZ$.
\item
$H^{2j}(\Ti_y^*({}^\pi F^H_{wx}(\lam))
\simeq
F^H_{wy}(\lam)(-j)^{\oplus c_j}$
for any $j\in\BZ$.
\end{tenumerate}
\end{theorem}


\bibliographystyle{unsrt}
\def\same{\,$\raise3pt\hbox to 25pt{\hrulefill}\,$}

\end{document}